\documentclass{article}[11]
\usepackage{amsmath}
\usepackage{graphicx}

\begin{document}

\title{How and why non smooth solutions \\
of the 3D Navier-Stokes equations \\ could possibly develop}

\author{Daniele Funaro}

%\date{}
\maketitle

\centerline{Dipartimento di Scienze Chimiche e Geologiche} \centerline{Universit\`a di
Modena e Reggio Emilia}
\centerline{Via Campi 103, 41125
Modena (Italy)} \centerline{daniele.funaro@unimore.it}

\begin{abstract} 
Fluid configurations in three-dimensions, displaying a plausible decay of regularity in a 
finite time, are suitably built and examined.
Vortex rings are the primary ingredients in this study. The full Navier-Stokes system
is converted into a 3D scalar problem, where appropriate numerical
methods are implemented in order to figure out the behavior of the solutions.
Further simplifications in 2D and 1D provide interesting toy problems, that may be used
as a starting platform for a better understanding of blowup phenomena.
\end{abstract}

\vspace{.8cm}
\par\smallskip

\noindent{AMS: 35Q30, 76N10.}

\numberwithin{equation}{section}

\par\medskip
\setcounter{equation}{0}
\section{Six collapsing rings}

The aim of this paper is to propose a way to build special explicit solutions of the entire set of time-dependent 
incompressible Navier-Stokes equations.
The model mainly consists of the law of momentum conservation, given by the vector equation: 
\begin{equation}\label{ns}
\frac{\partial {\bf v}}{\partial t} -\nu \bar\Delta {\bf v} + ({\bf v}\cdot \bar\nabla){\bf v}=-\bar\nabla p +{\bf f}
\end{equation}
where the velocity field ${\bf v}$ is required to be divergence-free, i.e.:
${\rm div}{\bf v}=0$. The last relation guarantees mass conservation.
The time $t$ belongs to the finite interval $[0,T]$. 
As customary, $\nu >0$ denotes the viscosity parameter.
The potential $p$ plays the role of pressure and ${\bf f}$ is a given force field. 
The equations are required to be satisfied in the whole three-dimensional space ${\bf R}^3$.
The symbol $\bar \Delta$ denotes the 3D vector Laplacian. Later, we will introduce another symbol
$\Delta$ (with no over bar) with a slightly different meaning.
\smallskip

Specifically, we will refer to those phenomena known as vortex rings (see, e.g.: \cite{akhme}, \cite{wu}). 
According to Fig.1, the fluid follows
a rotatory motion where the stream-lines revolve around the major circumferences of a doughnut.
As a consequence of diffusion, the movement of the particles is accompanied by a drifting of the
ring as indicated by the arrows. At the same time, a progressive reduction of the energy is also expected, depending
on the magnitude of $\nu$. We would like to see what happens when the ring is constrained inside an infinite cone, and in
particular to examine its behavior in an appropriate neighborhood of the vertex. There, the sections of the
ring that, in normal circumstances, tend to be approximated by circles, assume unusual shapes (see Fig.12).

\vspace{.1cm}
\begin{center}
\begin{figure}[h!]
\centerline{{\hspace{.3cm}\includegraphics[width=6.cm,height=4.7cm]{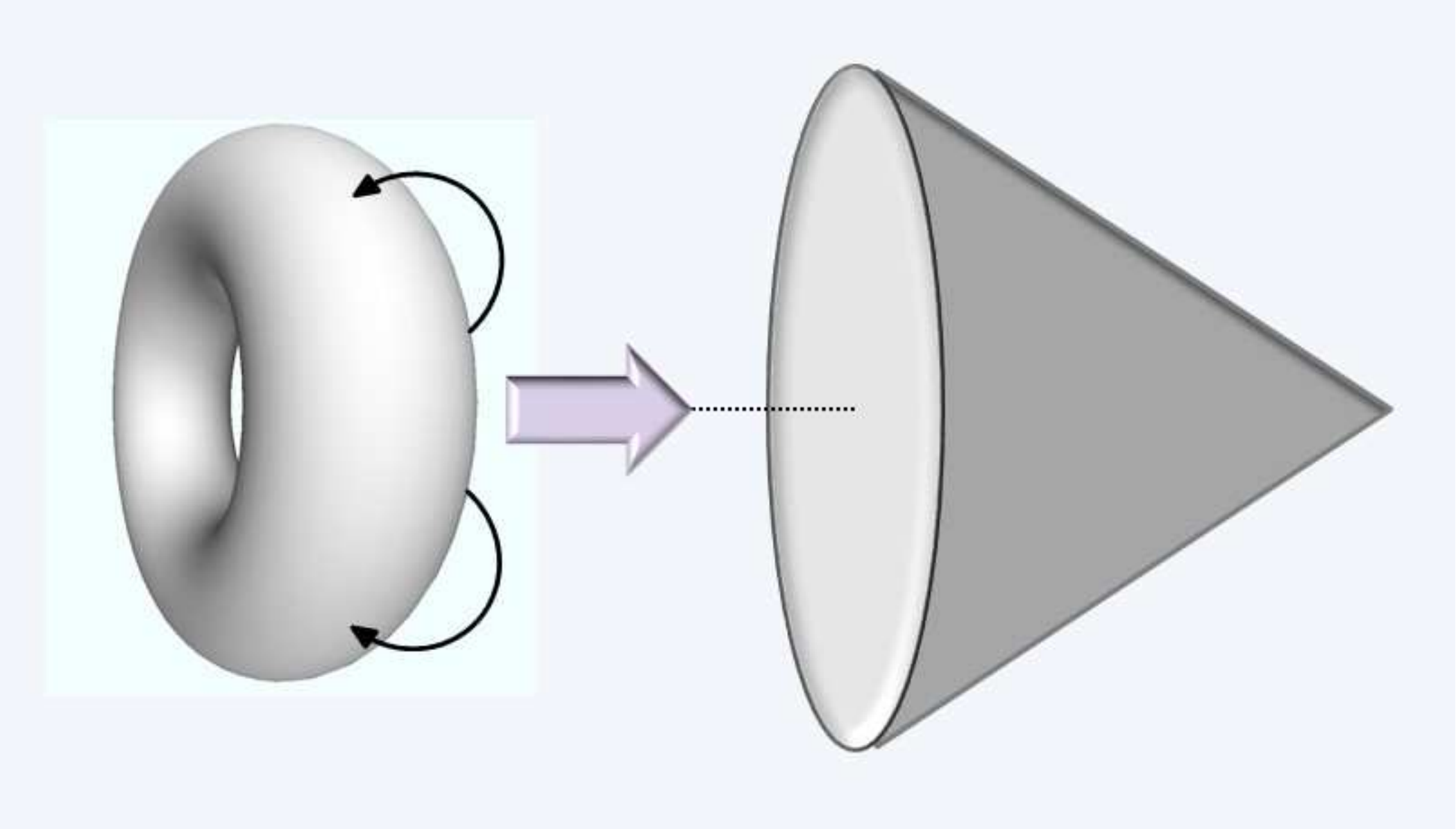}}
\hspace{.9cm}
{\hspace{-.3cm}\includegraphics[width=5.6cm,height=4.7cm]{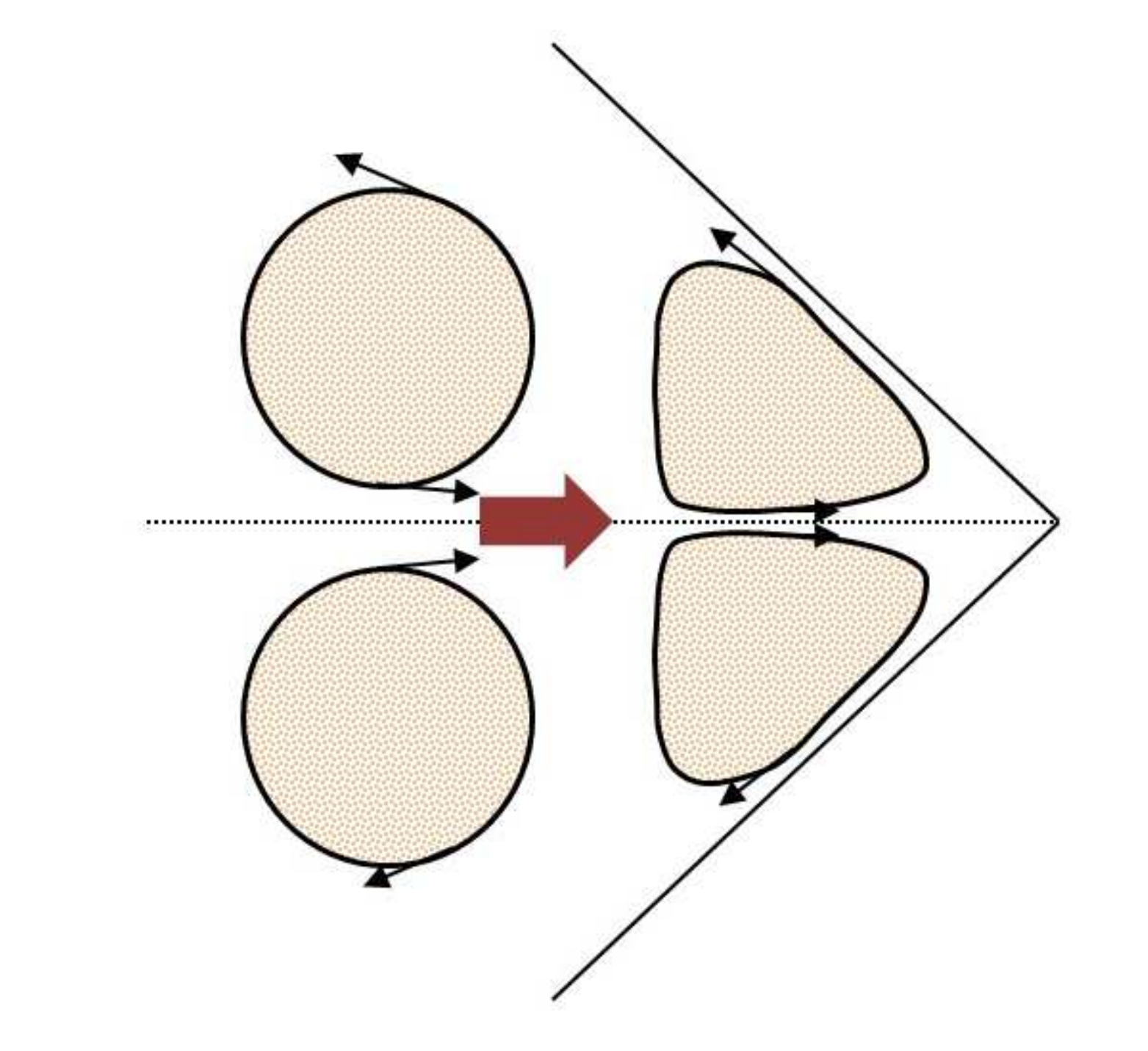}}}
%\vspace{-.5cm}
\begin{caption}{\small A vortex ring is constrained inside a cone. The sections
tend to change their shape when approaching the vertex.}
\end{caption}
\end{figure}
\end{center}
\vspace{-.2cm}

Since we want to avoid boundary conditions and have a solution defined on the whole space ${\bf R}^3$, we skip the
idea of the cones and we
divide instead the space into six virtual pyramidal regions as suggested by Fig.2. Each pyramid has an aperture of 90 degrees,
spanned by two independent angles $\theta$ and $\phi$. Six identical vortex rings
(or six sequences of them, lined up one after the other) are assembled along the six Cartesian semi-axes, so that they progress 
by maintaining a global symmetry and exerting reciprocal
constraints, without mixing each other. Such a congestion near the vertexes of the contiguous pyramids, may lead to a
possible singular behavior in proximity of the origin\footnote{Here, in proximity of the origin does not mean that the real troubles
will be exactly centered at the origin. The actual behavior will be clarified later in the exposition.}. Indeed, this is 
the eventuality we would like
to explore. For intense initial velocity fields and a very small diffusive parameter, there is the chance that smooth
solutions may, at some instant, lose regularity. 
\smallskip

To say the whole truth, we will not solve the just mentioned problem. Our rings will not move autonomously,
but they will be subject to external forces. This implies that ${\bf f}$ in (\ref{ns}) is going to be 
different from zero. By suitably manipulating the equations, we transfer part of the nonlinear term
on the right-hand side, so obtaining a forcing term ${\bf f}({\bf v})$ depending on the solution itself.
Assuming that the revised equation admits a unique solution, the field ${\bf f}({\bf v})$ (known a posteriori) is
interpreted as an external given force. Note that the new solution may not have any physical relevance.
These passages, that look like a trivial escamotage, have however
some hope to be useful. In fact, let us suppose that we are able to prove that ${\bf v}$ loses regularity
in a finite time, whereas ${\bf f}({\bf v})$ remains smooth (even if its knowledge is implicitly tied to
that of ${\bf v}$); this would mean that it is possible to generate singularities from regular data.
By `singularity' here we intend a degeneracy of some partial derivative of ${\bf v}$. 
It is known from the literature that a minimal degree of regularity for ${\bf v}$ is always preserved
during time. This means that we do not expect extraordinary explosions (we provide a more detailed 
explanation towards the end of section 8). It is important to remark that these mild forms of deterioration of the 
regularity are not clearly detected by standard numerical simulations. This makes our analysis a bit
uncertain.
\smallskip

We translate the full set of Navier-Stokes equations into a 3D nonlinear scalar differential equation,
where the unknown is a potential $\Psi$. Further simplifications in 2D and 1D, allow us to introduce
some toy problems aimed to provide a starting platform for possible theoretical advances.
By the way, we will not be able to prove rigorously the majority of the facts mentioned above. Some 
statements will be checked with the help of numerical experiments. Nevertheless, we believe that
the material given in the present paper establishes a strong foundation in view of more serious studies.
\smallskip

As far as the 3D incompressible fluid dynamics equations are concerned, the research on the regularity
of solutions has produced thousands of papers. A proof that the solutions maintain their smoothness 
during long-time evolution is at the moment not available. Indeed, the problem of describing the 
behavior in three space dimensions has always been borderline.  
Due to the viscosity term, smooth data are expected to produce solutions with an 
everlasting regular behavior. On the other hand, the lack of a conclusive 
theoretical analysis suggests the existence of possible counterexamples. 
%The subject has received a renewed impulse since  the institution of a famous prize (see \cite{fefferman}). 
The community supporting the idea that a blowup
may actually happen in a finite time, is growing, and numerous publications, both concerning the Euler
and the Navier-Stokes equations, are nowadays available. We cite here just a few titles, since an accurate review
would take too much time and effort. From the theoretical side, we mention: \cite{beale}, \cite{cannone}, \cite{chan}, 
\cite{foxall}, \cite{karch}. In \cite{gala} and \cite{tao}, possible scenarios regarding the development
of singularities are presented. From the numerical viewpoint, we quote \cite{grauer} and \cite{kerr}. 
Finally, sophisticated laboratory experiments on vortex rings at critical regimes are found for 
instance in \cite{keown} and \cite{lim}.

\par\medskip
\section{A suitable coordinates environment}

It is enough to study the Navier-Stokes problem on a single pyramidal subdomain and then
assemble the six pieces of solution (see the last section).
It is wise to work with a suitable system of of coordinates where the infinitesimal distance $ds$ is recovered by:
\begin{equation}\label{ds}
(ds)^2 = (dr)^2 + r^2 (d\theta )^2 + r^2 (d\phi )^2
\end{equation}
with $r$ denoting the radial variable, whereas $\theta$ and $\phi$ are angles.

\medskip
\begin{center}
\begin{figure}[h!]
\centerline{\hspace{-.3cm}\includegraphics[width=10.5cm,height=6.3cm]{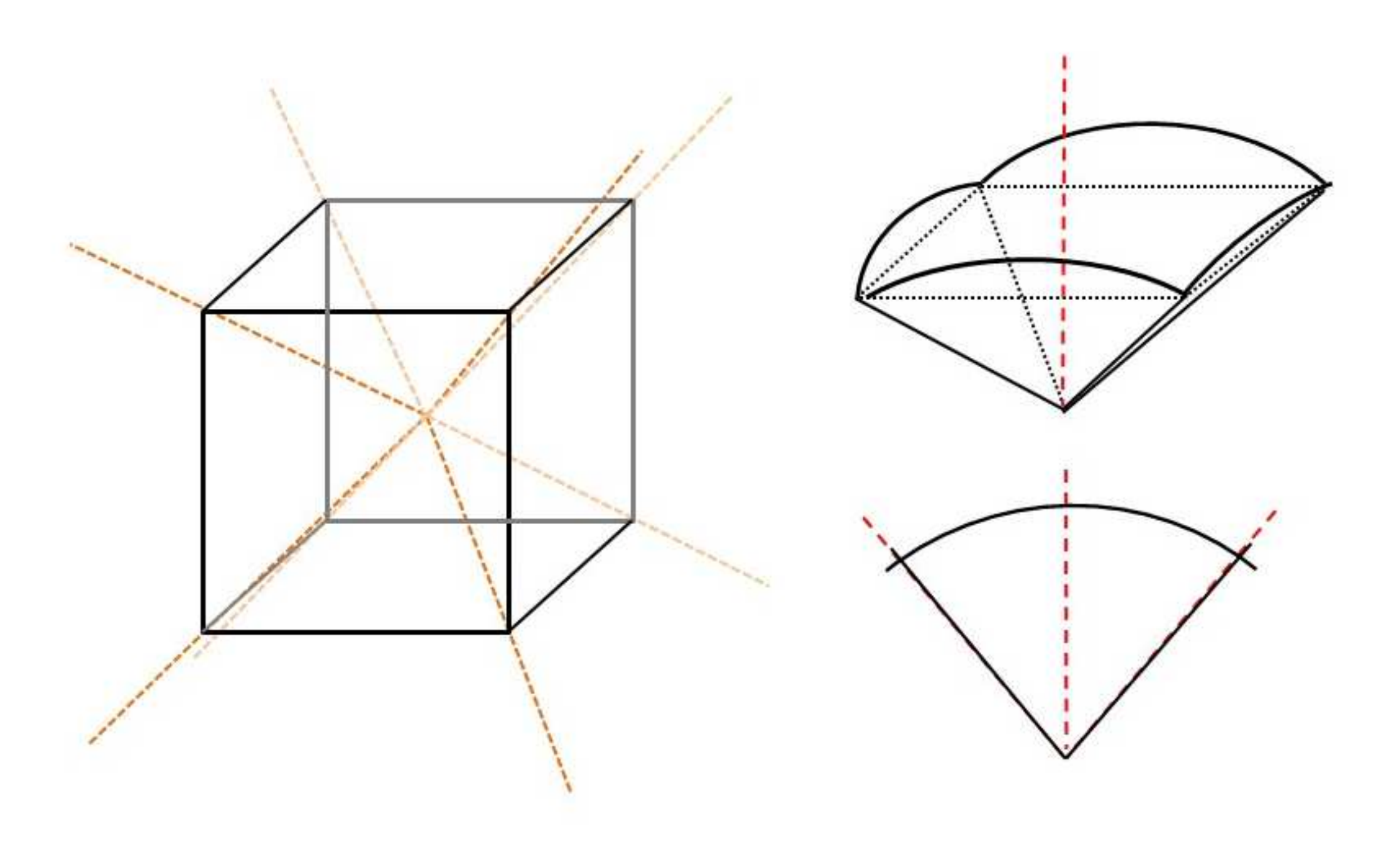}}
\vspace{-.5cm}
\begin{caption}{\small The whole three-dimensional space is virtually subdivided into six
pyramidal subdomains. These are separated by 12 triangular interfaces.}
\end{caption}
\end{figure}
\end{center}
\vspace{-.2cm}

Within this environment, the gradient of a potential $p$ is evaluated in the following way:
\begin{equation}\label{nabla}
\bar\nabla p = \left( \frac{\partial p}{\partial r}, \ \frac{1}{r}\frac{\partial p}{\partial \theta},
\ \frac{1}{r}\frac{\partial p}{\partial \phi}\right)
\end{equation}
For a given vector field ${\bf A}=(A_1, A_2, A_3)$,
we can compute some of the most classical differential operators:
%\begin{equation}\label{nablaa}
%\bar\nabla {\bf A}=\left( \frac{\partial A_1}{\partial r},\ \frac{1}{r}\frac{\partial A_2}{\partial \theta},
%\ \frac{1}{r}\frac{\partial A_3}{\partial \phi}\right)
%\end{equation}

\begin{equation}\label{diva}
{\rm div} {\bf A}= \frac{\partial A_1}{\partial r}+ \frac{2}{r}A_1+\frac{1}{r}\frac{\partial A_2}{\partial \theta}
+\frac{1}{r}\frac{\partial A_3}{\partial \phi}
\end{equation}

$$
{\rm curl} {\bf A}=\left( \frac{1}{r}\frac{\partial A_3}{\partial \theta}-
\frac{1}{r}\frac{\partial A_2}{\partial \phi}, \
-\frac{\partial A_3}{\partial r}- \frac{A_3}{r} +
\frac{1}{r}\frac{\partial A_1}{\partial \phi},\right.
$$
\begin{equation}\label{curla}
\qquad\qquad \left.
\frac{\partial A_2}{\partial r}+ \frac{A_2}{r} -
\frac{1}{r}\frac{\partial A_1}{\partial \theta}
\right)
\end{equation}

$$-\bar \Delta {\bf A} = {\rm curl}({\rm curl}{\bf A})= -\left(\frac{\partial^2 A_1}{\partial r^2} 
+\frac{4}{r}\frac{\partial A_1}{\partial r}  +\frac{2A_1}{r^2}+\frac{\Delta A_1}{r^2}, \qquad\qquad\ \right.$$
\begin{equation}\label{deltaav}
{\hspace{-.4cm}}\left. \frac{\partial^2 A_2}{\partial r^2} 
+\frac{2}{r}\frac{\partial A_2}{\partial r}  +\frac{\Delta A_2}{r^2} + \frac{2}{r^2}\frac{\partial A_1}{\partial \theta}, \
\frac{\partial^2 A_3}{\partial r^2} 
+\frac{2}{r}\frac{\partial A_3}{\partial r}  +\frac{\Delta A_3}{r^2}+ \frac{2}{r^2}\frac{\partial A_1}{\partial \phi}
\right)
\end{equation}
In the last expression we assumed that ${\rm div} {\bf A}=0$.
The symbol $\Delta$ without the upper bar denotes the usual Laplacian in the variables $\theta$ and $\phi$.
Applying $\Delta$ to the scalar functions $A_k$, $k=1,2,3$, leads us to the equality:
\begin{equation}\label{deltaa}
\Delta A_k =\frac{\partial^2 \hspace{-.05cm}A_k}{\partial \theta^2} + \frac{\partial^2 \hspace{-.05cm} A_k}{\partial \phi^2}
\end{equation}
Finally, we define the open set:
\begin{equation}\label{domain}
\Omega =\left\{ -{\textstyle\frac{\pi}{\omega}}< \theta < {\textstyle\frac{\pi}{\omega}}, \ 
-{\textstyle\frac{\pi}{\omega}}< \phi < {\textstyle\frac{\pi}{\omega}}\right\}
\end{equation}
In practice, we always choose $\omega =4$. A generic pyramidal domain corresponds to the set: $\Sigma = \{ (r,\theta ,\phi)\vert \
r>0, \ (\theta ,\phi )\in \Omega \}$

\par\medskip
\section{Stationary fields with singularity at the origin}

Before facing the general case (treated starting from section 7), we deal with some preliminary simplified examples.
We work in the reference frame $(r, \theta ,\phi)$ introduced in the previous section.
Our functions will be regular enough to allow for the exchange of the order of derivatives.
We start by discussing the case of a scalar potential $\Psi=\Psi (\theta , \phi)$ not depending on the variables 
$r$ and $t$. From this, we build the following vector potential:
\begin{equation}\label{campoas}
{\bf A}=\left( 0, \ -\frac{\partial \Psi}{\partial \phi}, \ \frac{\partial \Psi}{\partial \theta}\right)
\end{equation}
satisfying ${\rm div} {\bf A}=0$. Successively, we find the velocity field ${\bf v}=(v_1, v_2, v_3)$:
\begin{equation}\label{campovs}
{\bf v}={\rm curl} {\bf A}=\left( \frac{\Delta \Psi}{r}, \ -\frac{1}{r}\frac{\partial \Psi}{\partial \theta},
\ -\frac{1}{r}\frac{\partial \Psi}{\partial \phi}\right) =\frac{1}{r}(\Delta \Psi, 0 ,0) - \bar\nabla \Psi
\end{equation}
We recall that the symbol $\Delta$ is the scalar Laplacian in the variables
$\theta$ and $\phi$ (see (\ref{deltaa})). This means that: 
\begin{equation}\label{deltapsi}
\Delta \Psi =\frac{\partial^2 \Psi}{\partial \theta^2} + \frac{\partial^2 \Psi}{\partial \phi^2}
\end{equation}
By construction, the field ${\bf v}$ in (\ref{campovs}) satisfies the equation:
${\rm div}{\bf v}=0$. 
\smallskip

\noindent Going ahead, we compute:
\begin{equation}\label{curlvs}
{\rm curl}{\bf v}= -\bar \Delta {\bf A} =\frac{1}{r^2}\left( 0, 
\ \frac{\partial (\Delta\Psi)}{\partial \phi}, \ -\frac{\partial (\Delta\Psi)}{\partial \theta}
\right)
\end{equation}

$$-\bar \Delta {\bf v} = {\rm curl}({\rm curl}{\bf v})= -\frac{1}{r^3}\left( \Delta^2\Psi, 
\ \frac{\partial (\Delta\Psi)}{\partial \theta}, \ \frac{\partial (\Delta\Psi)}{\partial \phi}\right)
$$
\begin{equation}\label{deltavs}
= -\frac{1}{r^3}\left( \Delta^2\Psi +2\Delta \Psi, \ 0, \ 0\right) - \bar \nabla q_1
\end{equation}
where $\Delta^2\Psi =\Delta (\Delta \Psi )$. Here we find a first function $q_1 =\Delta\Psi /r^2$, playing the role of 
a scalar potential, in view of assembling the final pressure $p$ in (\ref{ns}).
A second function $q_2=-\frac12 \vert {\bf v}\vert^2$ takes part in the
vector relation:
\begin{equation}\label{vnablav}
({\bf v}\cdot \bar\nabla){\bf v}=-\bar\nabla q_2 -{\bf v}\times {\rm curl}{\bf v} 
\end{equation}
By making explicit the last member of the right-hand side in (\ref{vnablav}), we
get a third function $q_3$:
$$
{\bf v}\times {\rm curl}{\bf v}=\frac{1}{r^3}\left( \frac{\partial (\Delta\Psi)}{\partial \theta}
\frac{\partial \Psi}{\partial \theta} + \frac{\partial (\Delta\Psi)}{\partial \phi}
\frac{\partial \Psi}{\partial \phi}, \ \frac{\partial (\Delta\Psi)}{\partial \theta} \Delta\Psi, \
\frac{\partial (\Delta\Psi)}{\partial \phi} \Delta\Psi\right)$$
\begin{equation}\label{vcurlvs}
=\frac{1}{r^3}\left( \frac{\partial (\Delta\Psi)}{\partial \theta}
\frac{\partial \Psi}{\partial \theta} + \frac{\partial (\Delta\Psi)}{\partial \phi}
\frac{\partial \Psi}{\partial \phi} + (\Delta \Psi)^2, \ 0, \ 0\right) + \bar\nabla q_3
\end{equation}
with $q_3=\frac12 (\Delta\Psi)^2/r^2$. 

%$p_4 =K/2r^2$
\smallskip

By defining the global pressure $p=\nu q_1+q_2+q_3$ on the right-hand side of (\ref{ns})
and by setting ${\bf f}=0$, the whole Navier-Stokes system
is summarized in the fourth-order equation in the single scalar unknown $\Psi$:
\begin{equation}\label{deltapsis}
-\nu\Delta^2\Psi - 2\nu \Delta \Psi -\frac{\partial (\Delta\Psi)}{\partial \theta}
\frac{\partial \Psi}{\partial \theta} - \frac{\partial (\Delta\Psi)}{\partial \phi}
\frac{\partial \Psi}{\partial \phi} - (\Delta \Psi)^2=0
\end{equation}
It is convenient to write the above relation as a system of two second-order equations,
by introducing a new function $u$ such that:
\begin{equation}\label{eqpsinl}
u=\Delta \Psi 
\end{equation}
\begin{equation}\label{equnl}
\nu \Delta u + 2\nu u+ \nabla u\cdot \nabla \Psi + u^2 =0
\end{equation}
Here, the symbol $\nabla$ without the upper bar denotes the classical gradient
in two variables, i.e.:  $\nabla u =(\partial u/\partial \theta, \partial u/\partial \phi )$.
\smallskip

As far as boundary conditions are concerned, we first introduce the outward normal vector
$\bar n = (n_2, n_3)$ to the domain $\Omega$ defined in (\ref{domain}). At each one of the four corners, $\bar n$
is taken as the sum of the limits of the normal vectors along the two concurring sides (in this
case the norm of $\bar n$ is going to be equal to $\sqrt{2}$). We impose Neumann conditions
to both the unknowns $\Psi$ and $u$. This means that:
\begin{equation}\label{neumann}
\nabla\Psi \cdot \bar n=0 \qquad\qquad
\nabla u \cdot \bar n=0
\end{equation}
From (\ref{campovs}), the first relation implies that $v_2n_2+v_3n_3=0$
on $\partial\Omega$. This says that the velocity vector field is flattened on the
separation surfaces of the six pyramidal domains partitioning the whole three-dimensional space.
This construction holds with the exception of the point $r=0$, where our fields are singular.
\smallskip

By integrating the differential equation (\ref{equnl}) in $\Omega$, 
we discover the following compatibility condition for $u$:
$$
\nu \int_\Omega \Delta u \,d\theta d\phi + 2\nu \int_\Omega  u \,d\theta d\phi
+ \int_\Omega \nabla u\cdot \nabla \Psi \,d\theta d\phi
+\int_\Omega u^2 \,d\theta d\phi$$
$$= \nu \int_{\partial\Omega} \nabla u \cdot \bar n + 2\nu \int_\Omega  u \,d\theta d\phi+ 
\int_{\partial\Omega} u\nabla \Psi \cdot \bar n -\int_\Omega u^2 \,d\theta d\phi
+\int_\Omega u^2 \,d\theta d\phi =0$$
\begin{equation}\label{usebc}
\Rightarrow \int_\Omega u \,d\theta d\phi =0
\end{equation}

%$$r^2 \frac{\partial u}{\partial t} =\nu \Delta u + 2\nu u+ \nabla u\cdot \nabla \Psi + u^2 $$
\medskip
\section{Some preliminary numerical simulations}

In view of more sophisticated applications, we set up the computational machinery starting
from the one-dimensional version of the equations (\ref{eqpsinl}) and (\ref{equnl}).
Thus, we consider:
\begin{equation}\label{eqpsinl1}
u= \Psi^{\prime\prime}
\end{equation}
\begin{equation}\label{equnl1}
\nu (u^{\prime\prime}+2 u) +u^\prime\Psi^\prime +u^2=0
\end{equation}
where $u$ and $\Psi$ now depend exclusively on the variable $\phi$.
We then consider the Fourier expansions:
\begin{equation}\label{fexpa}
u(\phi)=c_0+\sum_{k=1}^\infty c_k \cos (\omega k \phi ) \qquad \qquad
\Psi (\phi)=d_0+\sum_{k=1}^\infty d_k \cos (\omega k \phi )
\end{equation}
where
\begin{equation}\label{cofur}
c_0=\frac{\omega}{2\pi}\int_{-\pi /\omega}^{\pi /\omega} {\hspace{-.2cm}}u(\phi ) \ d\phi
\qquad\quad
c_k=\frac{\omega}{\pi}\int_{-\pi /\omega}^{\pi /\omega} {\hspace{-.2cm}}u(\phi )
\cos (\omega k \phi ) \ d\phi, \ \ k\geq 1
\end{equation}
Analogous formulas hold for $d_k$, $k\geq 0$.
\smallskip

In this way, we are satisfying the boundary conditions $u^\prime=0$ and $\Psi^\prime=0$
at the endpoints $\phi =\pm \pi/\omega$.
As a consequence of (\ref{eqpsinl1}), for $k\geq 1$ the coefficients are connected
by the relation:
\begin{equation}\label{coeffco}
c_k = -\omega^2 k^2 d_k
\end{equation}
Moreover, the implication in (\ref{usebc}) suggests that $c_0=0$. Since $\Psi$ is
involved in the equations only through its derivatives, we can also set $d_0=0$.
\smallskip

From well-known trigonometric formulas, we get:
$$u^\prime\Psi^\prime =\sum_{{k=1}\atop{m=1}}^\infty km \omega^2 c_k d_m \sin (\omega k \phi )
\sin (\omega m \phi )$$
\begin{equation}\label{coeffco1}
=-\frac12 \sum_{{k=1}\atop{m=1}}^\infty km \omega^2 c_k d_m \Big[\cos (\omega (k+m) \phi )
-\cos (\omega (k-m) \phi )\Big]
\end{equation}

$$u^2 =\sum_{{k=1}\atop{m=1}}^\infty  c_k c_m \cos (\omega k \phi )
\cos (\omega m \phi )$$
\begin{equation}\label{coeffco2}
=\frac12 \sum_{{k=1}\atop{m=1}}^\infty c_k c_m \Big[\cos (\omega (k+m) \phi )
+\cos (\omega (k-m) \phi )\Big]
\end{equation}
\smallskip

For any fixed integer $n\geq 1$, by substituting (\ref{coeffco1}) and (\ref{coeffco2}) into the equation (\ref{equnl1}), 
we find out that, relatively to the mode $\cos (\omega n\phi)$, we must have:

$$\nu (2-n^2\omega^2)c_n+\frac12 \sum_{k+m=n} \Big[-km\omega^2 c_k d_m +c_kc_m\Big]$$
\begin{equation}\label{coeffren}
\qquad +~\frac12 \sum_{\vert k-m\vert =n}\Big[ km\omega^2 c_k d_m +c_kc_m\Big] =0
\end{equation}
All the indexes are greater or equal to one.
\smallskip

%sistema.m + funzsistemavar.m
\begin{center}
\begin{figure}[h!]
\centerline{{\includegraphics[width=6.cm,height=5.2cm]{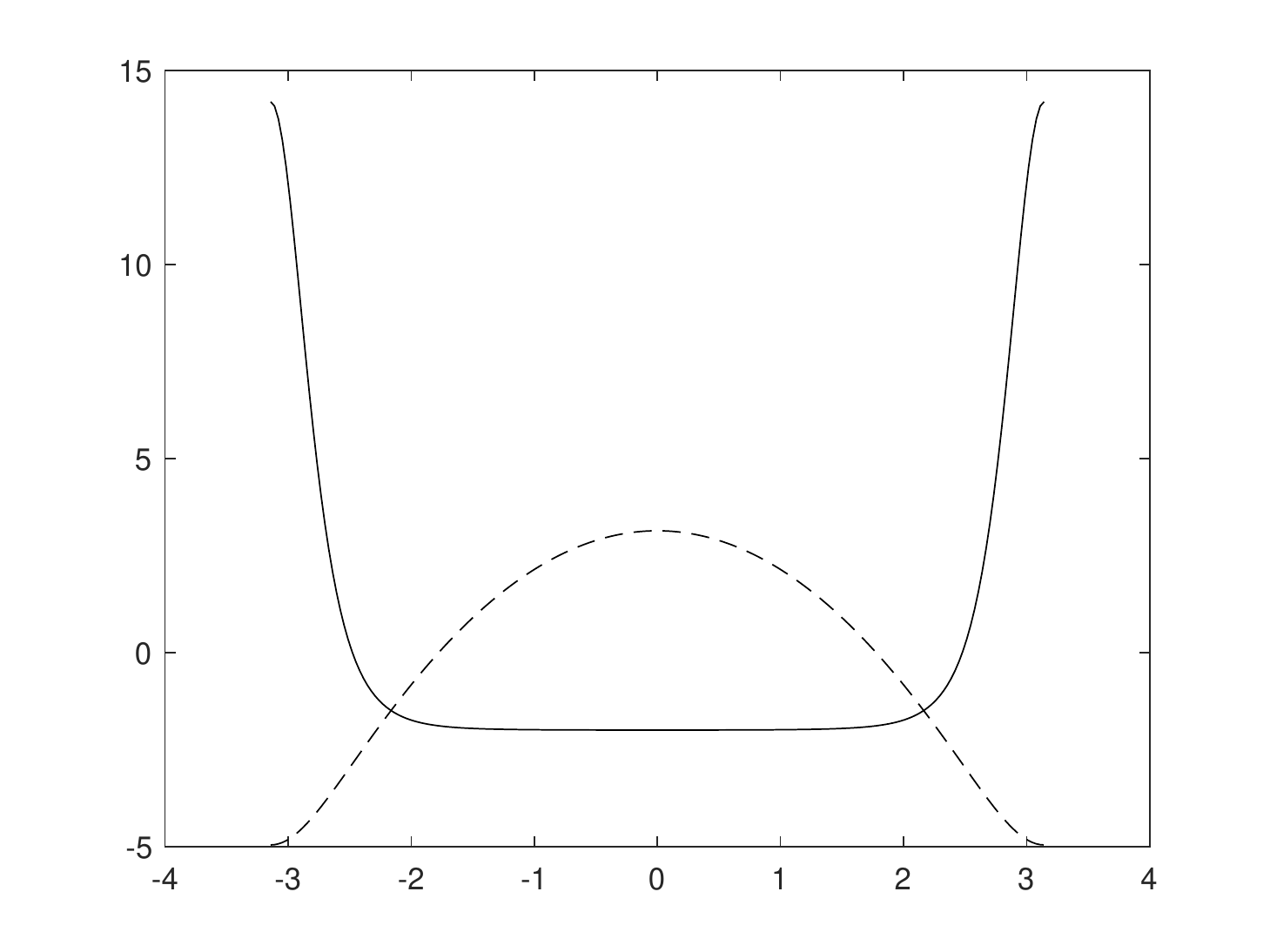}}
{\includegraphics[width=6.cm,height=5.2cm]{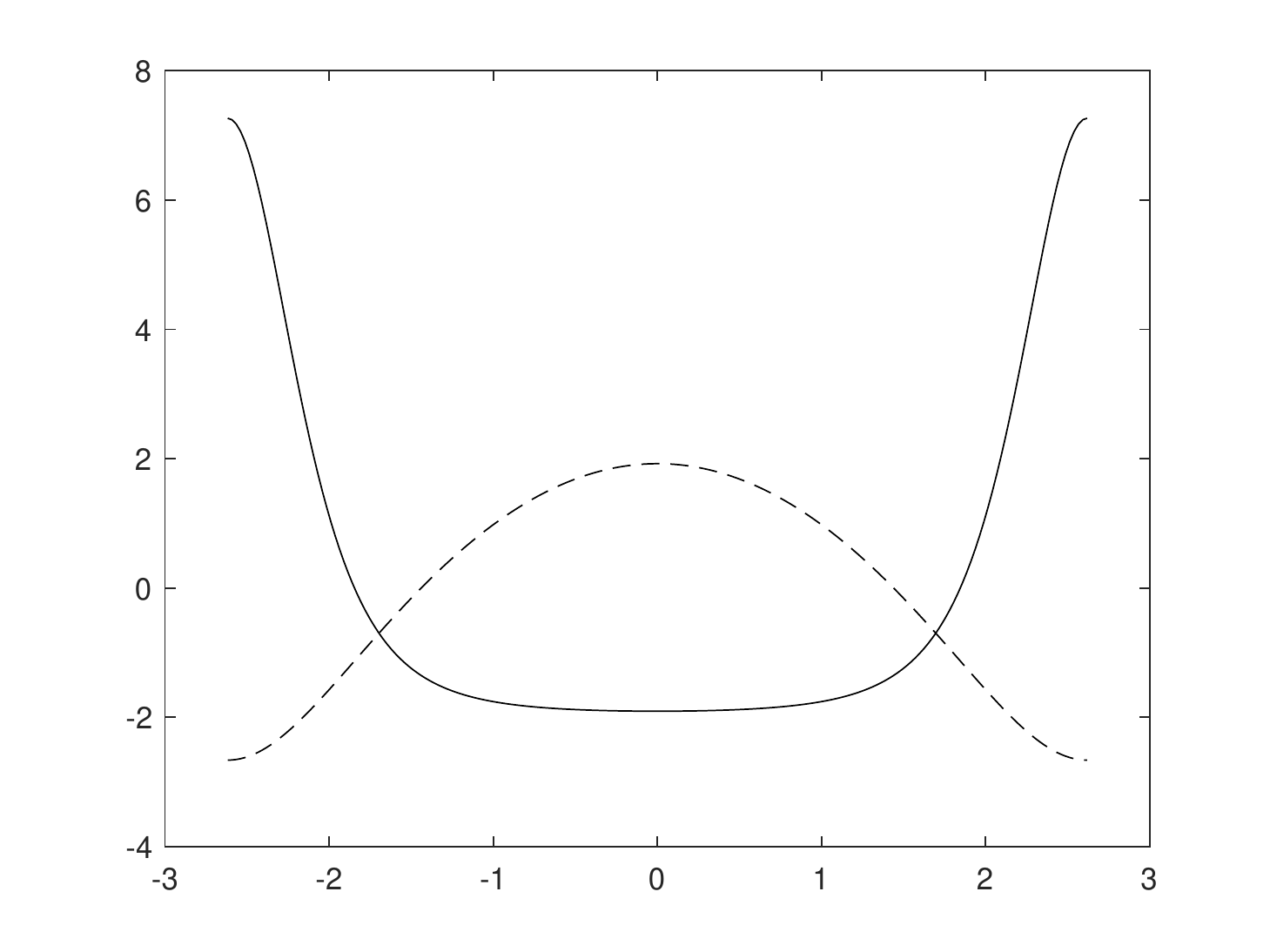}}}
\vspace{-.5cm}
\begin{caption}{\small Solutions $u$ (solid line) and $\Psi$ (dashed line) of the system (\ref{eqpsinl1}) and (\ref{equnl1})
for $\omega=1$ and $\omega=1.2$.}
\end{caption}
\end{figure}
\end{center}
\vspace{-.3cm}

Of course, the system (\ref{eqpsinl1}) and (\ref{equnl1}) always admits the trivial solutions $u=0$ and $\Psi=0$.
However, depending on the choice of $\omega$, another solution is available, that seems to be unique and rather stable.
The plots of Fig.3 show the couple of non-vanishing solutions obtained for $\omega=1$ and $\omega =1.2$. Their shape agrees with the one expected
from a rough preliminary theoretical inspection.
Nevertheless, by increasing $\omega$ (recall that we would like to have $\omega=4$), the corresponding
solutions display a certain number of oscillations, leading to a velocity field ${\bf v}$ that does not reflect
the behavior that we are trying to simulate. The first transition to the new typology of solutions happens when $\omega =\sqrt{2}$
(note that $u=\cos\sqrt{2}\phi$ is the first eigenfunction such that $u^{\prime\prime}+2u=0$).
\smallskip

The analysis of the set of equations (\ref{eqpsinl1}) and (\ref{equnl1}) has affinity with the study
of {\sl diffusive logistic models}, where the existence of non zero solutions depends on the location of a parameter
relatively to the distribution of the eigenvalues of the diffusive operator. The literature on the subject
is rather extensive. Since we did not find explicit references to our specific case, we limit 
our citations to the generic review paper \cite{taira}.
\smallskip

This first attempt to build a time-stationary solution presenting a singularity at the point $r=0$
has been a failure. Nevertheless, the construction is useful for further decisive improvements,
that are discussed in the following sections.

\par\medskip
\section{Evolutive field in the non viscous case}

In order to prepare the ground for the general case, other special solutions
may be proposed in the specific case when $\nu =0$. Within this new setting, we start from the vector
potential:
\begin{equation}\label{potaa}
{\bf A}=\left( 0, \ -r^2\frac{\partial \Psi}{\partial \phi}, \ r^2\frac{\partial \Psi}{\partial \theta}\right)
\end{equation}
The scalar potential $\Psi$ is function of the variables $\theta$, $\phi$ and $t$. 
Note that ${\rm div}{\bf A}=0$.
As before, we can determine the velocity field:
\begin{equation}\label{velo2}
{\bf v} ={\rm curl} {\bf A}=\left( r \Delta \Psi, \ -3r\frac{\partial \Psi}{\partial \theta},
\ -3r\frac{\partial \Psi}{\partial \phi}\right) =(ru, 0 ,0) - \bar\nabla q_0
\end{equation}
where, for convenience, we defined: $u=\Delta\Psi+6\Psi$ and $q_0=3 r^2\Psi$.
Of course, we still have: ${\rm div}{\bf v}=0$. The next step is to
evaluate the curl of ${\bf v}$ and its vector Laplacian:
\begin{equation}\label{curlv2}
{\rm curl}{\bf v}=\left( 0, 
\ \frac{\partial u}{\partial \phi}, \ -\frac{\partial u}{\partial \theta}
\right)
\end{equation}
\begin{equation}\label{deltav2}
-\bar \Delta {\bf v} = \frac{1}{r}\left( -\Delta u, 
\ \frac{\partial u}{\partial \theta}, \ \frac{\partial u}{\partial \phi}\right)
= \left( -\frac{\Delta u}{r}, 0, 0\right) - \bar \nabla q_1
\end{equation}
where now $q_1=-u$.
Concerning the nonlinear term, we obtain:
$${\bf v}\times {\rm curl}{\bf v}=\left( 3r\frac{\partial u}{\partial \theta}
\frac{\partial \Psi}{\partial \theta} + 3r\frac{\partial u}{\partial \phi}
\frac{\partial \Psi}{\partial \phi}, \ r\frac{\partial u}{\partial \theta} \Delta\Psi, \
r\frac{\partial u}{\partial \phi} \Delta\Psi\right)$$
\begin{equation}\label{nl2}
=\left( 3r\frac{\partial u}{\partial \theta}
\frac{\partial \Psi}{\partial \theta} + 3r\frac{\partial u}{\partial \phi}
\frac{\partial \Psi}{\partial \phi} - r(\Delta \Psi)^2, \ -f_2, \ -f_3\right) + \bar\nabla q_3
\end{equation}
where $q_3=\frac12 (r\Delta\Psi)^2$. The pressure in (\ref{ns}) can be now defined as
$p=(\partial q_0 /\partial t)+q_2+q_3$, where $q_2$ was introduced in (\ref{vnablav}).
This time, the forcing term ${\bf f}=(f_1,f_2,f_3)$
is not zero and we must have:
\begin{equation}\label{fterm}
f_1=0\qquad\qquad
f_2=-6r\frac{\partial \Psi}{\partial \theta} \Delta \Psi \qquad \qquad
f_3=-6r \frac{\partial \Psi}{\partial \phi}\Delta \Psi 
\end{equation}
Thus, ${\bf f}$ implicitly depends on the unknown itself. 
Alternatively, we can set $q_3=\frac12 (ru)^2$ and define $f_2$ and $f_3$ accordingly.
\smallskip

By putting together all the terms ($\bar \Delta {\bf v}$ excluded since $\nu =0$), the first 
component of the system yields the equation:
\begin{equation}\label{eq1com}
\frac{\partial u}{\partial t} -3 \nabla u\cdot \nabla \Psi + (\Delta\Psi)^2 =0
\end{equation}
with:
\begin{equation}\label{eq2com}
u=\Delta\Psi+6\Psi
\end{equation}
Neumann boundary conditions on $\partial\Omega$ will be assumed for both $u$ and
$\Psi$.
\smallskip

If we instead define $q_3=\frac12 (ru)^2$, the equation (\ref{eq1com}) takes the form:
\begin{equation}\label{eq1comalt}
\frac{\partial u}{\partial t} -3 \nabla u\cdot \nabla \Psi + u^2 =0
\end{equation}
There is no big difference concerning the behavior of the solutions for the two versions.
\smallskip

The second and the third components of the Navier-Stokes system are totally `absorbed'
by $f_2$, $f_3$ and by the gradient of pressure. The idea is that one can solve
(\ref{eq1com}) and (\ref{eq2com}) with zero right-hand side ($f_1=0$).
A posteriori, the couple $(f_2, f_3)$ is recovered from (\ref{fterm}) without solving 
any further equation. More comments about this procedure will be provided at the
end of section 7 (see, in particular, relation (\ref{nsrev})).
\smallskip

Here, the Laplacian $\Delta u$ is not taken into account
($\nu =0$) because the dependance with respect to $r$ in the expression 
(\ref{deltav2}) is not homogeneous with the other terms.
\smallskip

Some analysis can be carried out for the one-dimensional version of 
(\ref{eq1com}) and (\ref{eq2com}). In this case, we get the two equations:
\begin{equation}\label{eq2com2}
\frac{du}{dt}-3u^\prime\Psi^\prime +(\Psi^{\prime\prime} )^2=0
\end{equation}
\begin{equation}\label{eq1com1}
u= \Psi^{\prime\prime}+6\Psi
\end{equation}
According to (\ref{fexpa}), from (\ref{eq1com1}) a relation is soon established 
between the Fourier coefficients for $k\geq 0$:
\begin{equation}\label{eq2coef}
c_k = -\omega^2 k^2 d_k +6 d_k
\end{equation}
In particular, the coefficient $c_0=6d_0$ does not need to be zero.
Considering that:
$$(\Psi^{\prime\prime} )^2 =\sum_{{k=1}\atop{m=1}} k^2m^2 \omega^4 d_k d_m \cos (\omega k \phi )
\cos (\omega m \phi )$$
\begin{equation}\label{coefpp}
=\frac12 \sum_{{k=1}\atop{m=1}} k^2m^2 \omega^4 d_k d_m \Big[\cos (\omega (k+m) \phi )
+\cos (\omega (k-m) \phi )\Big]
\end{equation}
we can obtain the counterpart of (\ref{coeffren}) for a fixed integer $n\geq 1$, i.e.:
$$\frac{d c_n}{dt}+\frac12 \sum_{k+m=n} \Big[\mu_1 km\omega^2 c_k d_m +\mu_2 k^2 m^2 \omega^4 d_kd_m\Big]$$
\begin{equation}\label{coeffrenn}
\qquad +~\frac12 \sum_{\vert k-m\vert =n}\Big[ -\mu_1 km\omega^2 c_k d_m +\mu_2 k^2 m^2 \omega^4 d_kd_m\Big] =0
\end{equation}
with $\mu_1=3$ and $\mu_2=1$. For $n=0$, the first summation in (\ref{coeffrenn}) disappears.
Thus, we must have:
\begin{equation}\label{coeff0}
\frac{d c_0}{dt}
-\frac32 \sum_{j=1}^{\infty}j^2\omega^2 c_j d_j +\frac12 \sum_{j=1}^{\infty}j^4\omega^4 d_j^2=0
\end{equation}
By virtue of (\ref{eq2coef}), for $\omega =4$ we come out with the estimate:
\begin{equation}\label{coeff0ev}
\frac{d c_0}{dt}= \sum_{j=1}^{\infty}j^2\omega^2 (9-2j^2\omega^2) d_j^2 
=\sum_{j=1}^{\infty}\frac{j^2\omega^2 (9-2j^2\omega^2)}{(j^2\omega^2 -6)^2} \ c_j^2 
<-2\sum_{j=1}^{\infty} c_j^2 \leq 0
\end{equation}
\smallskip

Suppose that, for $t\rightarrow \hat t$ (where $\hat t$ may be finite or infinite), $u$ converges to a limit 
in $L^2(-\pi/4, \pi/4 )$. Let us also suppose that $\sum_{j=1}^{\infty} c_j^2$ tends to a positive constant.
Then (\ref{coeff0ev}) tells us that 
$\lim_{t\rightarrow \hat t} c_0(t)$ does not exist (i.e.: $c_0$ diverges negatively) and
this is against the hypothesis of convergence in $L^2(-\pi/4, \pi/4 )$.
The remaining possibility is that $\sum_{j=1}^{\infty} c_j^2$ tends to zero, which means that
$u$ converges to a constant function (i.e., $u$ minus its average tends to zero).
As a consequence, in the framework of functions with zero average, we expect 
$u$ and $\Psi$ to converge to zero, unless some compatibility conditions between the coefficients
($\mu_1$ and $\mu_2$) of differential systems of the type of (\ref{eq2com2})-(\ref{eq1com1}) are satisfied. 
We will be more precise in the coming section.

\par\medskip
\section{A simple 1D problem}

The results of the previous sections suggest to study more carefully the
system in the single variable $\phi$, involving the two unknowns $u$ and $\Psi$:
\begin{equation}\label{eq11d}
u= \Psi^{\prime\prime} +\lambda \Psi
\end{equation}
\begin{equation}\label{eq21d}
\frac{du}{dt}-\nu u^{\prime\prime}+ \Big[-\mu_1 u^\prime\Psi^\prime + 
\mu_2(\Psi^{\prime\prime} )^2\Big]=0
\end{equation}
where Neumann type boundary conditions are assumed at the endpoints, i.e.:
$u^\prime=0$ and $\Psi^\prime =0$ for $\phi=\pm \pi /\omega$. In (\ref{eq11d})-(\ref{eq21d}),
$\lambda$, $\mu_1$ and $\mu_2$ are real parameters.
\smallskip

After integration of (\ref{eq21d}) between $- \pi /\omega$ and $\pi /\omega$, one gets: 
$$
\frac{d}{dt}\int_{-\pi/\omega}^{\pi/\omega} u \ d\phi= -\mu_1 \int_{-\pi/\omega}^{\pi/\omega} u \Psi^{\prime\prime}\ d\phi
\ -\mu_2 \ \int_{-\pi/\omega}^{\pi/\omega} (\Psi^{\prime\prime})^2\ d\phi
$$
\begin{equation}\label{equint}
= -(\mu_1 +\mu_2) \int_{-\pi/\omega}^{\pi/\omega} (\Psi^{\prime\prime})^2 \ d\phi+\lambda \mu_1 
\int_{-\pi/\omega}^{\pi/\omega} (\Psi^{\prime})^2\ d\phi
\end{equation}
where we used the rule of summation by parts and imposed the boundary conditions.
\smallskip

We also recall the following Poincar\`e type inequality:
\begin{equation}\label{poinc}
\int_{-\pi/\omega}^{\pi/\omega} (\Psi^{\prime})^2 \ d\phi\ \leq \ \frac{1}{\omega^2}
\int_{-\pi/\omega}^{\pi/\omega} (\Psi^{\prime\prime})^2 \ d\phi
\end{equation}

We partly rediscover the system of section 4 by setting $\lambda =0$, $\mu_1=1$,
$\mu_2=-1$. In this case, the relation (\ref{equint}) is compatible with the fact that the first
Fourier coefficient $c_0$ of $u$ must remain zero during time evolution (see (\ref{usebc})). The system of section 5 
is instead recovered by setting $\nu =0$, $\lambda =6$, $\mu_1=3$,
$\mu_2=1$.
In the general case, relation (\ref{coeff0ev}) becomes:
\begin{equation}\label{coeff0evg}
\frac{d c_0}{dt}=
\frac12 \sum_{j=1}^{\infty}\frac{j^2\omega^2 [\lambda\mu_1- (\mu_1+\mu_2)j^2\omega^2]}{(j^2\omega^2-\lambda)^2} \ c_j^2 
\end{equation}
Note that we are in the peculiar situation where the right-hand side of (\ref{coeff0evg}) does not contain the coefficient
$c_0$.
If $c_0$ does not depend on $t$, the above formula may allow for non-vanishing Fourier coefficients $c_j$, $j\geq 1$, if
suitable compatibility conditions hold between the parameters $\lambda$, $\mu_1$ and $\mu_2$. Namely, it
is necessary that the generic quantity:
\begin{equation}\label{qua}
Q=\lambda \mu_1 -(\mu_1 +\mu_2)j^2\omega^2
\end{equation}
assumes both positive and negative values
depending on $j$. For $\lambda =6$, $\mu_1=3$, $\mu_2=1$, $\omega =4$, we have that $Q=18-64j^2$ is always negative,
which confirms that the projection of the system (\ref{eq2com2})-(\ref{eq1com1}) onto the space of zero average
functions does not admit solutions different from zero.
\smallskip

A numerical test has been made by truncating the Fourier sums at a given $N$ and the results are visible in
Fig.4. The diffusion parameter is $\nu=.01$.  The other parameters are: $\lambda=-3$, $\mu_1=.5$, $\mu_2=-1.5$. 
This choice ensures that $Q$ in  (\ref{qua}) may attain both positive and negative values, depending
on the frequency mode involved. In the computation we enforced the condition $c_0(t)=0, \forall t\in [0,T]$,
basically by not including the zero mode in the expansion of $u$ and noting that its knowledge is
not requested in the evaluation of the right-hand side of (\ref{eq21d}).
\smallskip

The explicit Euler scheme for $0\leq t\leq T=1.48$ has been implemented
with a sufficiently small time-step. The coefficients $c_n$ and $d_n$ are computed for $1\leq n\leq N$, with $N=50$. 
For $\omega=4$, the initial guess $u_0$ has been set in such a way that: $u_0(\phi )=\cos (\omega\phi)$.
Note that the sign of $u$ at time $t=0$ has a nontrivial impact on the branch of solution we would
like to follow. The coefficients of $\Psi$ are recovered at any iteration through relation  (\ref{eq11d}).
Very similar conclusions hold when $\mu_1=1$, $\mu_2=-4$ and $\lambda$ is negative. This particular case will be
rediscussed later in section 10. 

\vfill

%modello1d
\begin{center}
\begin{figure}[p!]
{\centerline{\includegraphics[width=9.4cm,height=8.2cm]{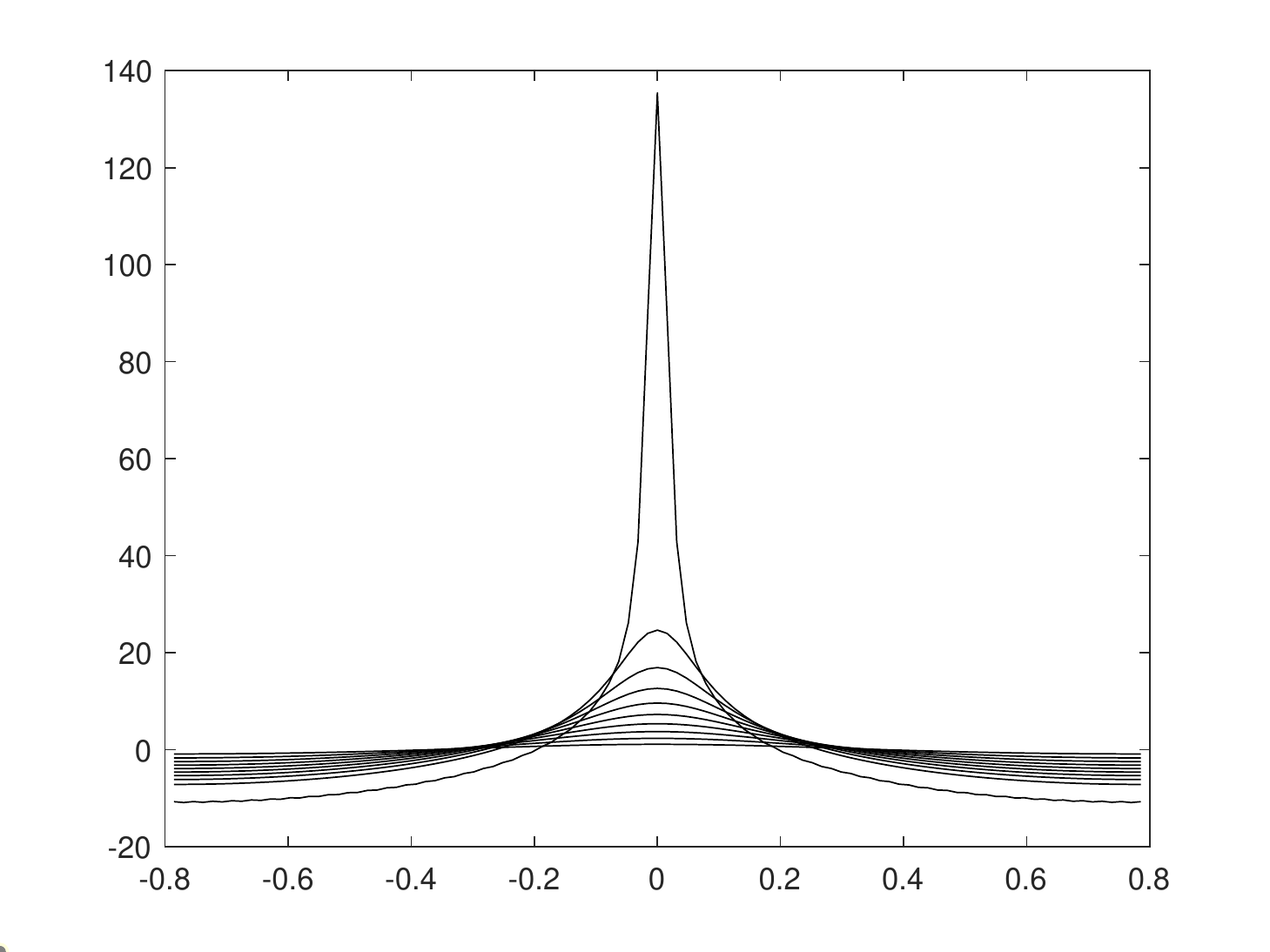}}
%\hspace{-.4cm}
\centerline{\includegraphics[width=9.4cm,height=8.2cm]{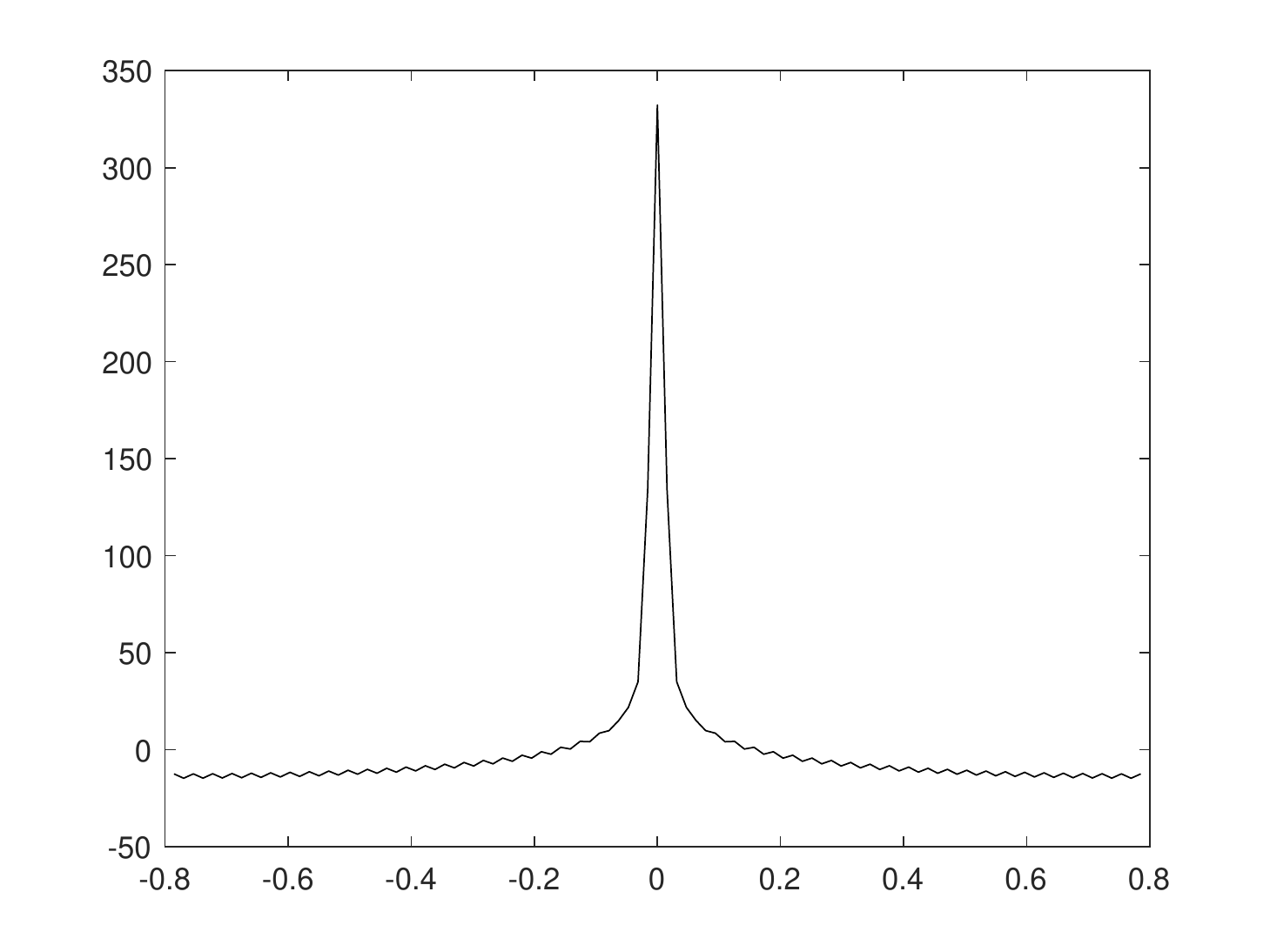}}}
\vspace{-.1cm}
\begin{caption}{\small Plots of the solution $u$ at equispaced times (top), for $\lambda=-3$, $\mu_1=.5$, $\mu_2=-1.5$,
with $0\leq t\leq 1.48$. At time $T=1.486$, the approximate solution starts producing oscillations (bottom).}
\end{caption}
\end{figure}
\end{center}
%\vspace{-.5cm}

The discrete solution is clearly trying to assume the shape of a very pronounced cusp at the center of the interval.
For times $t$ larger than $T=1.486$, the simulation first  produces oscillations and then overflow.
Without a theoretical analysis, we are however unable to decide if there is a real blowup of $u$ or just a bad
behavior of its derivatives.
With such a small value of $\nu$, the effects of diffusion are very mild, but still they may prevent the regularity
of the solution to degenerate, while the numerical instability may only be consequence of a non appropriate choice
of the discretization parameter $N$. Note that cosinus expansions are extremely easy to implement, especially in view to
enforce Neumann boundary conditions. Nevertheless, they give accurate results only in presence of high regularity,
which is not the prerogative of the functions we are examining here.
\smallskip

From our rough analysis, what we learned in this section is that, for certain values of the parameters, the model problem admits only the 
steady state solution identically zero. For other suitable choices of the parameters, non-vanishing stable
solutions emerge. They may display a degeneracy of the regularity after a certain time.
\smallskip

Nonlinear parabolic equations presenting a blowup of the solution in a finite time, are widely studied.
A classical example is:
\begin{equation}\label{nonlin}
\frac{\partial u}{\partial t} -\nu \Delta u= f(u)
\end{equation}
with Dirichlet boundary conditions. Assume that $f$ is convex and $f(u)>0$, for $u>0$.
If for some $a>0$, the integral $\int_a^\infty (1/f(u))du$ is finite, then the solution
of  (\ref{nonlin}) blows-up when the initial datum is sufficiently large.
This and similar other questions are reviewed for instance in \cite{hu}.
\smallskip

Our system may have affinities with other model equations deriving from the most disparate
applications. A prominent example is the Cahn-Hilliard equation (see \cite{cahn}).
The literature on this subject is quite extensive, so that we just limit ourselves to
mention the recent book \cite{miranville}. In its basic formulation, the Cahn-Hilliard equation
takes the form:
\begin{equation}\label{kh}
\frac{\partial u}{\partial t} =-\nu \Delta (f(u) -\gamma \Delta u)
\end{equation}
where $\nu>0$ and $\gamma >0$ are suitable parameters. It is often written as a system
after introducing the function $\mu = f(u) -\gamma \Delta u$. Typical boundary conditions
are of Neumann type, i.e.: $\partial u/\partial\bar n=0$ and $\partial \mu/\partial\bar n=0$.
Existence of nontrivial attractors is proven in several circumstances. A standard choice
for the nonlinear term is $f(u)=u^3-u$. Although there is no direct connection with
our model problem, it is not improbable that some similarities may emerge when approaching the 
study from the theoretical viewpoint.

\vfill

\par\medskip
\section{The most general case}

We start with the same vector potential as in (\ref{campoas}), but
$\Psi=\Psi (t, r, \theta , \phi)$ now also depends on the variable $r$.
We introduce $\Phi$ and $u$ such that:
\begin{equation}\label{dphi}
\Phi = \frac{\partial \Psi}{\partial r} +\frac{\Psi}{r}=
\frac{1}{r}\frac{\partial (r\Psi)}{\partial r} 
\end{equation}
\begin{equation}\label{defu}
 u=\Delta\Psi + r\frac{\partial (r\Phi)}{\partial r}
=\Delta\Psi + r\frac{\partial^2 (r\Psi)}{\partial r^2}
\end{equation}
Based on these assumptions, the velocity field ends up to be:
\begin{equation}\label{defv}
{\bf v} =(v_1, v_2, v_3)=\left( \frac{1}{r} \Delta \Psi, \ -\frac{\partial \Phi}{\partial \theta},
\ -\frac{\partial \Phi}{\partial \phi}\right) =\left(\frac{u}{r}, 0 ,0\right) - \bar\nabla q_0
\end{equation}
with $q_0=r\Phi$.
Through the use of standard calculus we also get:
\begin{equation}\label{curlvg}
 {\rm curl}{\bf v}= \frac{1}{r^2}\left( 0, 
\ \frac{\partial u}{\partial \phi}, \ -\frac{\partial u}{\partial \theta}
\right)
\end{equation}

$$-\bar \Delta {\bf v} = {\rm curl}({\rm curl}{\bf v})= \left( -\frac{1}{r^3}\Delta u, 
\ \frac{1}{r^2}\frac{\partial^2 u}{\partial r\partial \theta} -\frac{1}{r^3} 
\frac{\partial u}{\partial \theta},
\ \frac{1}{r^2}\frac{\partial^2 u}{\partial r\partial \phi} -\frac{1}{r^3} 
\frac{\partial u}{\partial \phi}\right)$$
\begin{equation}\label{delvg}
=\left(  -\frac{1}{r^3}\Delta u - \frac{\partial^2 }{\partial r^2}\left(\frac{u}{r}\right),
 \ 0, \ 0\right) - \bar \nabla q_1
\end{equation}
with $q_1=-(\partial /\partial r)(u/r)$. We recall that the symbols $\Delta$ and $\nabla$
(without the upper bars) do not contain partial derivatives with respect to $r$.
\smallskip

Regarding the nonlinear term, we have:
$$
{\bf v}\times {\rm curl}{\bf v}=\left( \frac{1}{r^2}\frac{\partial u}{\partial \theta}
\frac{\partial \Phi}{\partial \theta} + \frac{1}{r^2}\frac{\partial u}{\partial \phi}
\frac{\partial \Phi}{\partial \phi}, \ \frac{1}{r^3}\frac{\partial u}{\partial \theta} \Delta\Psi, \
\frac{1}{r^3}\frac{\partial u}{\partial \phi} \Delta\Psi\right)$$
\begin{equation}\label{nlvg}
=\left( \frac{1}{r^2}\frac{\partial u}{\partial \theta}
\frac{\partial \Phi}{\partial \theta} + \frac{1}{r^2}\frac{\partial u}{\partial \phi}
\frac{\partial \Phi}{\partial \phi}- \frac12 \frac{\partial}{\partial r}\hspace{-.1cm}
\left(\frac{\Delta\Psi}{r}\right)^{\hspace{-.13cm}2}, \ -f_2, \ -f_3\right) + \bar\nabla q_3
\end{equation}
In the above expression we introduced the following functions:
\begin{equation}\label{fnoto}
q_3=\frac12 \hspace{-.08cm}\left(\frac{\Delta\Psi}{r}\right)^{\hspace{-.13cm}2} \qquad f_2=-\frac{\Delta \Psi}{r^2}
\frac{\partial^2 (r\Phi)}{\partial r \partial \theta} \qquad f_3=-\frac{\Delta \Psi}{r^2}
\frac{\partial^2 (r\Phi)}{\partial r \partial \phi}
\end{equation}
In alternative, we can define $q_3=\frac12 (u/r)^2$ and adjust $f_2$ and $f_3$ accordingly.
\smallskip

After having defined the pressure $p=(\partial q_0 /\partial t)+\nu q_1+q_2+q_3$ 
(with $q_2$ given in (\ref{vnablav})) and the forcing term
${\bf f}=(0, f_2, f_3)$, the first component of the vector momentum equation (\ref{ns}) is synthetically
represented by the scalar equation:
$$\frac{1}{r}\frac{\partial u}{\partial t}- \nu \left( \frac{1}{r^3}\Delta u + 
\frac{\partial^2 }{\partial r^2}\left(\frac{u}{r}\right)\right)$$
\begin{equation}\label{fineq}
-\frac{1}{r^2}\frac{\partial u}{\partial \theta}
\frac{\partial \Phi}{\partial \theta} - \frac{1}{r^2}\frac{\partial u}{\partial \phi}
\frac{\partial \Phi}{\partial \phi}+ \frac12 \frac{\partial}{\partial r}{\hspace{-.1cm}}
\left(\frac{\Delta\Psi}{r}\right)^{\hspace{-.13cm}2} =0
\end{equation}
If $\Psi$ does not depend on $r$, we return to the case studied in section 3 by
setting $\Phi =\Psi /r$, $u=\Delta \Psi$. If $\Psi$ is function of $r$ only through the
factor $r^2$, we come back to the case studied in section 5.
\smallskip

With little manipulation, we finally arrive at the system of two second-order equations:
$$\frac{\partial u}{\partial t}- \nu \left( \frac{\Delta u}{r^2} + 
\frac{\partial^2 u}{\partial r^2}-\frac{2}{r} \frac{\partial u}{\partial r}+\frac{2u}{r^2}\right)$$
\begin{equation}\label{fineq2}
+\frac{1}{r}\left[-\nabla u\cdot\nabla {\hspace{-.11cm}} \left( \frac{\partial \Psi}{\partial r} 
+\frac{\Psi}{r}\right)+ \Delta\Psi \ \Delta {\hspace{-.11cm}} \left( \frac{\partial \Psi}{\partial r} 
-\frac{\Psi}{r}\right)\right] =0
\end{equation}
\begin{equation}\label{finequ}
u=\Delta\Psi +r^2\frac{\partial^2\Psi}{\partial r^2}
+2r\frac{\partial\Psi}{\partial r}
\end{equation}
For both the unknowns $u$ and $\Psi$, we will require Neumann type boundary conditions 
on $\partial\Omega$ (see (\ref{neumann})), for any value of $r>0$. For $r=0$, both $u$ and $\Psi$ must vanish.
A suitable decay for $r\rightarrow +\infty$ is also assumed.
\smallskip

We can make some heuristic considerations about the above system.
First of all, we introduce the two functionals:
\begin{equation}\label{elle1}
{\cal L}_1u = \frac{\Delta u}{r^2} + 
\frac{\partial^2 u}{\partial r^2}-\frac{2}{r} \frac{\partial u}{\partial r}+\frac{2u}{r^2}
\end{equation}
\begin{equation}\label{elle2}
{\cal L}_2\Psi = \Delta\Psi +r^2\frac{\partial^2\Psi}{\partial r^2}
+2r\frac{\partial\Psi}{\partial r}
\end{equation}
Afterwards, we take for instance the two low-order eigenmodes:
\begin{equation}\label{eig1}
u_0(r,\theta, \phi)= -\gamma^2 r^2\chi (r) \cos (\omega\theta )\cos (\omega\phi )
\end{equation}
\begin{equation}\label{eig2}
\Psi_0(r,\theta, \phi)= \chi (r) \cos (\omega\theta )\cos (\omega\phi )
\end{equation}
Here, for a given $\gamma >0$, the function $\chi$ is defined as:
\begin{equation}\label{chidef}
\chi (r) =\frac{1}{\sqrt{\gamma r}}\ J_{\sigma +\frac12} (\gamma r)
\end{equation}
where $J_{\sigma +\frac12}$ is the {\sl spherical Bessel's function} of the first kind.
This implies that $\chi$ solves the differential equation:
\begin{equation}\label{bessel}
\frac{d^2 \chi}{dr^2}+\frac{2}{r}\frac{d\chi}{dr} -\sigma (\sigma +1)\frac{\chi}{r^2}=-\gamma^2\chi
\end{equation}
In truth, the expression in (\ref{chidef}) is valid up to a multiplicative constant.
If $u_0$ in (\ref{eig1}) is taken as an initial guess, its sign is crucial for the successive evolution (see later on).
By choosing $\sigma$ in such a way that $\sigma (\sigma +1)=2\omega^2$, a straightforward computation 
passing through (\ref{bessel}) shows that:
\begin{equation}\label{elle10}
{\cal L}_1u_0 = -\gamma^2 u_0
\end{equation}
\begin{equation}\label{elle20}
{\cal L}_2\Psi_0 = -\gamma^2 r^2 \Psi_0= u_0
\end{equation}
By using again (\ref{bessel}), the last expression can be rewritten as:
\begin{equation}\label{elle20m}
u_0={\cal L}_2\Psi_0 = \Delta \Psi_0 +({\cal L}_2 -\Delta ) \Psi_0
= \Delta \Psi_0 +(2\omega^2 -\gamma^2 r^2 ) \Psi_0 = \Delta \Psi_0 +\lambda \Psi_0
\end{equation}
with $\lambda =2\omega^2 -\gamma^2 r^2$. This means that in first approximation,
one can suppose that: 
$\ u\approx \Delta \Psi +\lambda \Psi$ (although $\lambda$ depends on $r$).
\smallskip

We now proceed with further approximations. When $\omega =4$, we must have $\sigma (\sigma +1)=16$,
that provides: $\sigma \approx 5.18$. From classical estimates on Bessel's
functions, the behavior of $\chi$ in (\ref{chidef})  is like  $r^\sigma$ near the origin
(up to multiplicative constants). By denoting with $r_M>0$ the first nontrivial zero of $\chi$,
we can say that:
\begin{equation}\label{chiap}
\chi (r)\ \approx \ r^\sigma (r_M-r) \qquad {\rm for} \ 0\leq r\leq r_M
\end{equation}
The first nontrivial zero of the Bessel's function $J_{\sigma +\frac12}$ for $\sigma\approx 5.18$, is approximately
$z\approx 9.56$. Thus, we must have $\gamma =z/r_M$.
\smallskip

Relation (\ref{chiap}) specifies that $u$ and $\Psi$ decay to zero quite fast near the origin.
Thus, we will not expect any deterioration of the regularity in the neighborhood of $r=0$.
If something strange may happen, it will be at some place located at a distance from the
origin (see footnote 1).
\smallskip

%\begin{table}[h!]
%\begin{center}
%\begin{tabular}{lcccl}
%  $k$ & $\alpha_k$ & $z_k$  \\
%\hline
%1 & 3.53 & 7.62 \\
%2 & 7.51 & 12.23 \\
%3 & 11.51 & 16.70 \\
%4 & 15.50 & 21.08 \\
%5 & 19.50 & 25.41 \\
%6 & 23.50 & 29.71 \\
%\hline
%\end{tabular}
%\end{center}
%\end{table}

We continue this rough analysis by introducing a new parameter $\alpha \leq 1$. If $\hat r$ 
is a point such that:
\begin{equation}\label{chial}
\frac{d\chi}{dr}(\hat r ) +\frac{\chi (\hat r)}{\hat r}=\alpha\frac{\chi (\hat r)}{\hat r}
\end{equation}
by making use of $\chi$ in (\ref{chiap}), we obtain:
\begin{equation}\label{r0est}
\hat r~\approx \frac{\sigma +1-\alpha}{\sigma +2-\alpha}~r_M \ \approx \
\frac{6.18-\alpha}{7.18-\alpha}~r_M  < r_M
\end{equation}
Recalling the definition of $\lambda$, we also have:
\begin{equation}\label{lamest}
\lambda (\alpha ) = 32 -\gamma^2 \hat r^2 =32 -z^2 \left( \frac{\hat r}{r_M}\right)^{\hspace{-.1cm}2}
\approx \ 32- (9.56)^2 \left( \frac{6.18-\alpha}{7.18-\alpha}\right)^{\hspace{-.1cm}2}
\end{equation}
\smallskip

Going back to the equation (\ref{fineq2}), as
far as the initial guess $u_0=-\gamma^2\hat r^2\Psi_0$ is concerned, we can argue in a similar way. 
If $\hat r$ is such that:
\begin{equation}\label{relalfa}
\frac{\partial\Psi_0}{\partial r}(\hat r )+\frac{\Psi_0(\hat r )}{\hat r}=\alpha\frac{\Psi_0 (\hat r )}{\hat r}\qquad
\frac{\partial\Psi_0}{\partial r}(\hat r )-\frac{\Psi_0(\hat r )}{\hat r}=(\alpha -2)\frac{\Psi_0 (\hat r )}{\hat r}
\end{equation}
for small times $t$, the nonlinear term in square brackets, changes in accordance to what studied in section 6, i.e.:
\begin{equation}\label{fineq20}
\frac{\partial u}{\partial t}- \nu {\cal L}_1 u
+\ \frac{1}{\hat r^2}\left[-\mu_1\nabla u\cdot\nabla \Psi + \mu_2(\Delta\Psi)^2 \right] =0
\end{equation}
with $\mu_1=\alpha$ and $\mu_2=\alpha -2$.
As far as the equation (\ref{finequ}) is concerned, we are induced to write:
\begin{equation}\label{finequf}
u=\Delta\Psi +\lambda (\alpha) \Psi
\end{equation}
with $\lambda$ depending on $\alpha$ as in (\ref{lamest}). In the one-dimensional counterpart,
the quantity in (\ref{qua}) would take the value: $Q=\alpha \lambda (\alpha) -2(\alpha -1)j^2\omega^2$.
For $0<\alpha <1$, we get $\lambda (\alpha ) <0$ and $Q$ may actually change sign. As an example, we may
set $\alpha=.5$, so that $\lambda (\alpha ) \approx -34$ and $Q\approx -17+16 j^2$.
\smallskip

Roughly speaking, by fixing $\hat r$ in the interval $]0, r_M[$, we may encounter situations
similar to those examined in section 6, bringing to a (supposed) degeneracy of the regularity of the solutions.
This does not mean that such kind of troubles must actually manifest in the framework of the real
3D problem, especially because our preliminary analysis was oversimplified. We will better consolidate 
our knowledge in section 11, but unfortunately we will still remain far from
rigorous proofs. In the next section, we try some numerical simulations on the
global 3D problem. The aim is to check whether anomalous situations may effectively occur.
\smallskip

We think it is wise to better clarify the passages made in this section. We got a functional equation of the
type $G(\Psi )=0$, that can be obtained by replacing $u$ defined in (\ref{finequ}) into (\ref{fineq2}).
The aim was to solve the Navier-Stokes equation (\ref{ns}). Therefore, we can write:
$${\bf 0}=\frac{\partial {\bf v}}{\partial t} -\nu \bar\Delta {\bf v} + ({\bf v}\cdot \bar\nabla){\bf v}+\bar\nabla p -{\bf f}
$$
\begin{equation}\label{nsrev}
=\left[(G(\Psi),\ 0,\ 0)-\bar\nabla \left(\frac{\partial q_0}{\partial t}+\nu q_1+q_2+q_3\right)+ (0, f_2, f_3)
\right]+\bar\nabla p -{\bf f}
\end{equation}
After setting $p=(\partial q_0 /\partial t)+\nu q_1+q_2+q_3$ and ${\bf f}=(0, f_2, f_3)$, we actually arrive at the
relation $G(\Psi )=0$. In this way, the pressure is not an unknown of the system, since it can be built in
dependance of ${\bf v}$. Similarly, we have a forcing term ${\bf f}$ which is not given a priori, but still depends
on the unknown. At the end, we are not solving the autonomous movement of a fluid. Our vortex ring will
develop under the action of forces that depend on its dynamics. This evolution may have not physical interest and
we do not expect the results to be easily interpreted from the fluid mechanics viewpoint.
By the way, our interest here is mainly focused on the analytical viewpoint. Indeed, let us suppose that
the development of ${\bf v}$ presents some deterioration of smoothness in a finite time, then 
two eventualities may happen. If ${\bf f}$ also loses regularity, we end up with proving nothing, because
it is reasonable to assume that a bad forcing term may give raise to bad solutions. If we can show instead that ${\bf f}$
maintains a certain degree of regularity (even if it depends on the solution itself), then these results start becoming interesting.

\par\medskip
\section{Full 3D discretization}

In order to discretize the full system (\ref{fineq2})-(\ref{finequ}), we consider the series: 
\begin{equation}\label{fexpabid}
u=\sum_{{k=0}\atop{i=0}}^\infty c_{ki} \cos (\omega k \theta )\cos (\omega i \phi ) \quad\qquad 
\Psi =\sum_{{k=0}\atop{i=0}}^\infty d_{ki} \cos (\omega k \theta )\cos (\omega i \phi )
\end{equation}
where the Fourier coefficients depend on $r$ and $t$. 
In this fashion we are respecting the Neumann boundary constraints as prescribed in 
(\ref{neumann}). Here, we decided to set $c_{00}=d_{00}=0$.
For $n\geq 0$ and $l \geq 0$, the mode $\cos (\omega n \theta )\cos (\omega l \phi )$
is associated with the evolution of the corresponding coefficient $c_{nl}$:
$$\frac{\partial c_{nl}}{\partial t}-\nu \left( -(n^2+l^2)\omega^2\frac{c_{nl}}{r^2} +\frac{\partial^2 c_{nl}}{\partial r^2}
-\frac{2}{r} \frac{\partial c_{nl}}{\partial r}+\frac{2c_{nl}}{r^2}\right)$$
$$+\frac{1}{4r} \sum_{{k+m=n}\atop{i+j=l}} \left[(km+ij)\omega^2 c_{ki}\left(\frac{\partial d_{mj}}{\partial r}
+\frac{d_{mj}}{r}  \right)\right.$$
$$\qquad \qquad \left. +(k^2+i^2)(m^2+j^2)\omega^4d_{ki} \left(\frac{\partial d_{mj}}{\partial r}
-\frac{d_{mj}}{r}  \right)\right]$$
$$+\frac{1}{4r} \sum_{{\vert k-m\vert =n}\atop{i+j=l}} \left[(-km+ij)\omega^2 c_{ki}\left(\frac{\partial d_{mj}}{\partial r}
+\frac{d_{mj}}{r}  \right)\right.$$
$$\qquad \qquad \left. +(k^2+i^2)(m^2+j^2)\omega^4d_{ki} \left(\frac{\partial d_{mj}}{\partial r}
-\frac{d_{mj}}{r}  \right)\right]$$
$$+\frac{1}{4r} \sum_{{k+m=n}\atop{\vert i-j\vert =l}} \left[(km-ij)\omega^2 c_{ki}\left(\frac{\partial d_{mj}}{\partial r}
+\frac{d_{mj}}{r}  \right)\right.$$
$$\qquad \qquad \left. +(k^2+i^2)(m^2+j^2)\omega^4d_{ki} \left(\frac{\partial d_{mj}}{\partial r}
-\frac{d_{mj}}{r}  \right)\right]$$
$$+\frac{1}{4r} \sum_{{\vert k-m\vert =n}\atop{\vert i-j\vert =l}} \left[(-km-ij)\omega^2 c_{ki}\left(\frac{\partial d_{mj}}{\partial r}
+\frac{d_{mj}}{r}  \right)\right.$$
\begin{equation}\label{coeffrenbid}
\qquad {\hspace{2.6cm}} \left. +(k^2+i^2)(m^2+j^2)\omega^4d_{ki} \left(\frac{\partial d_{mj}}{\partial r}
-\frac{d_{mj}}{r}  \right)\right] =0
\end{equation}
where $c_{nl}$ and $d_{nl}$ are related via (\ref{finequ}) in the following way:
\begin{equation}\label{coeffrenbidc}
c_{nl}= -(n^2+l^2)\omega^2d_{nl} +r^2\frac{\partial^2 d_{nl}}{\partial r^2}
+2r \frac{\partial d_{nl}}{\partial r}
\end{equation}
The above formulas, based on simple trigonometric identities, generalize those proposed in
the previous sections. The two coefficients $c_{00}$ and $d_{00}$ will remain equal to zero, for all $r\geq 0$,
as time passes. Therefore, it is necessary to check whether a suitable integral of the nonlinear
term satisfies a compatibility condition (see section 11).
\smallskip

We compute approximate solutions where $r$ belongs to the interval $[0, r_M]$ for
some $r_M>0$. We impose homogeneous Dirichlet boundary conditions to $u$ and $\Phi$ at $r=0$ and $r=r_M$.
The final time is $T=0.11$.
The derivatives with respect to the variable $r$ are approximated by central finite-differences.
The discretization in time is performed by the explicit Euler scheme with a rather small time-step.
This allows us to easily update the coefficients $c_{nl}$ at each iteration. The coefficients $d_{nl}$ are
obtained at each step by solving an implicit 1D boundary-value problem which is recovered by a central finite-differences
discretization of (\ref{coeffrenbidc}).

%iniziale.m 
\vspace{.2cm}
\begin{center}
\begin{figure}[h!]
\centerline{\hspace{-.5cm}{\includegraphics[width=7.3cm,height=5.8cm]{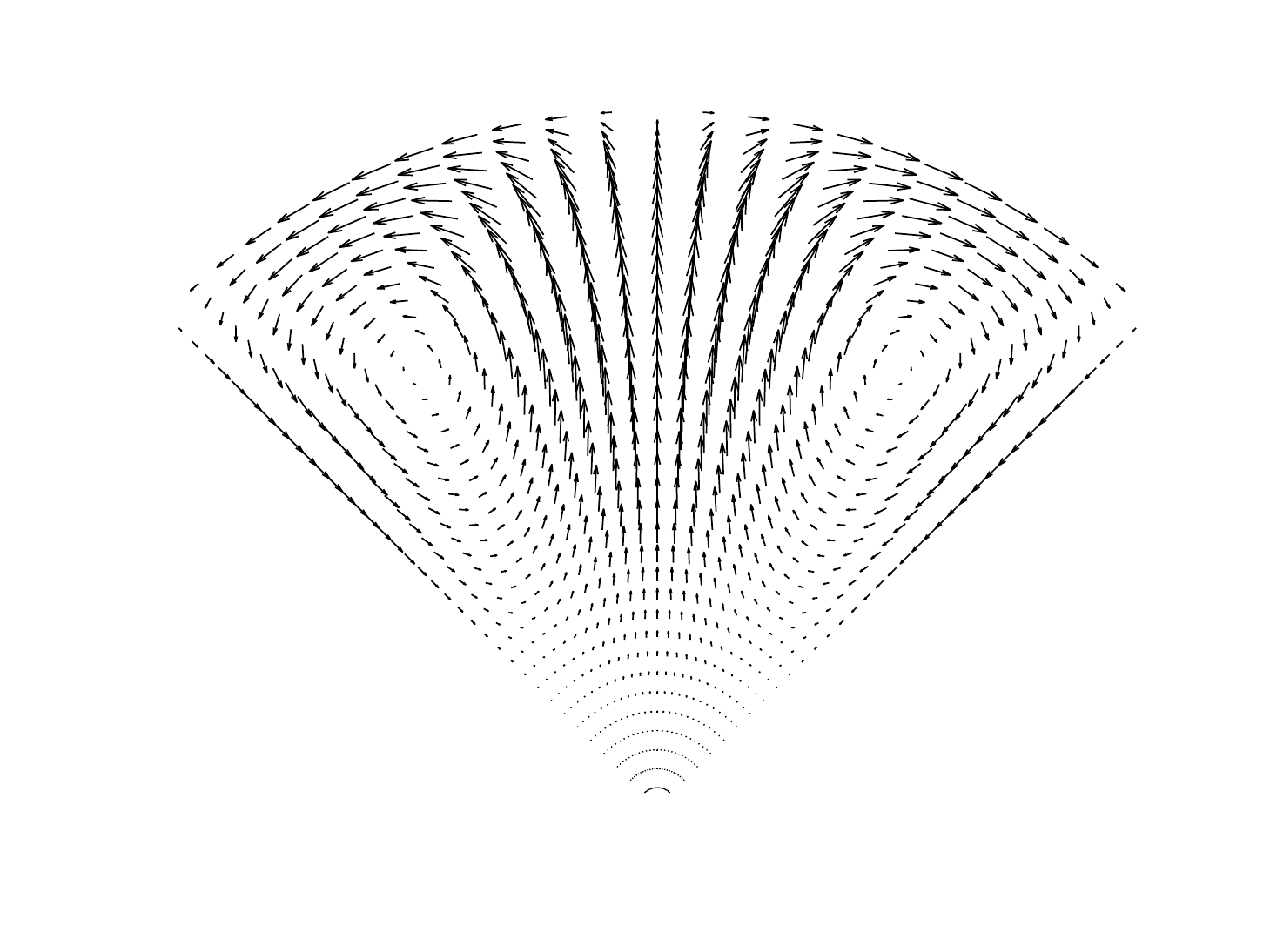}}
\hspace{-.9cm}{\includegraphics[width=7.3cm,height=5.6cm]{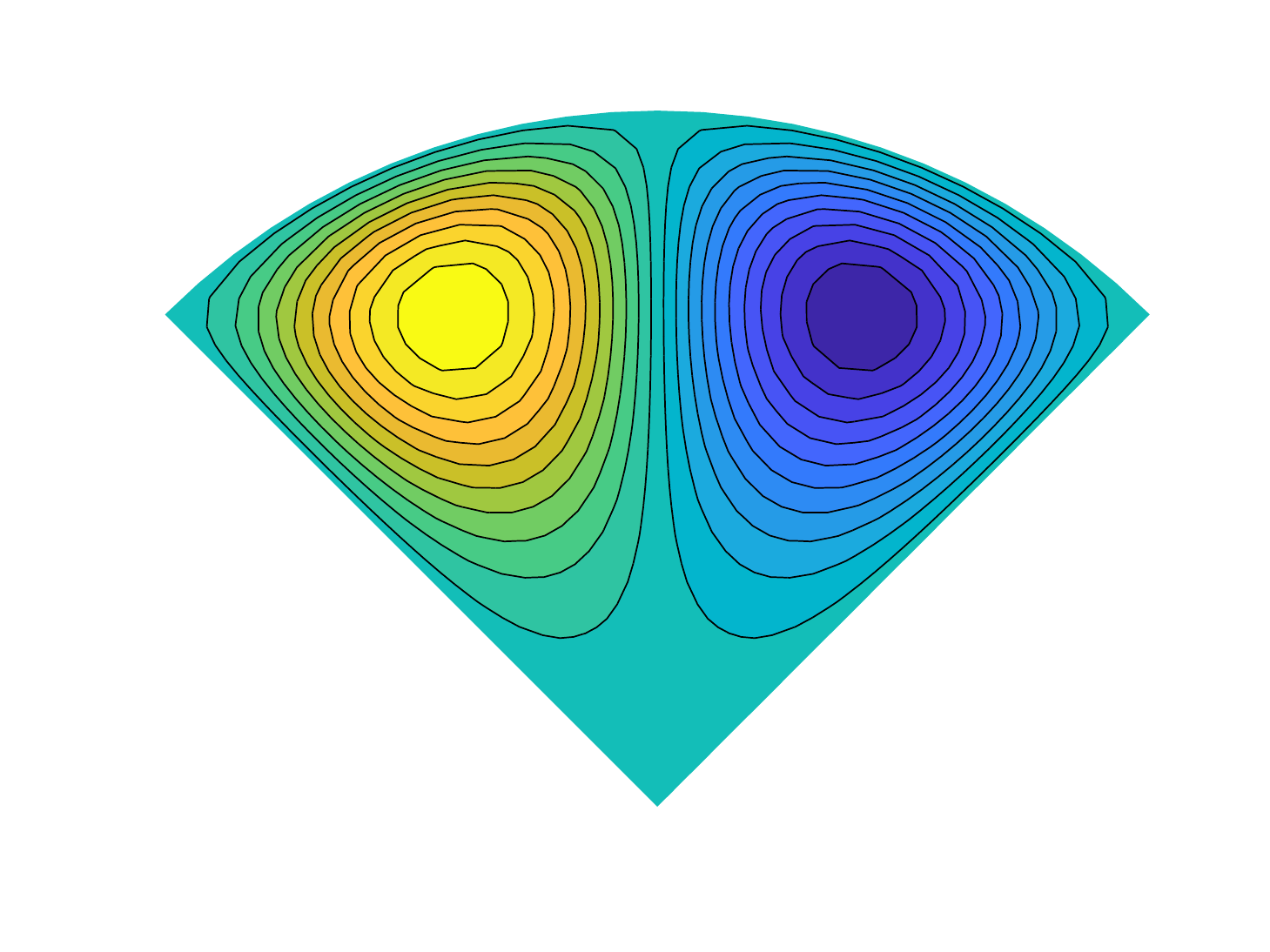}}}
\vspace{-.2cm}
\begin{caption}{\small Vector representation of the initial field ${\bf v}_0$ (left), corresponding to $u_0$ in 
(\ref{init}), together with the level lines of the third component of ${\bf A}$ (right). 
These pictures are referred to the section obtained for $\phi =0$.}
\end{caption}
\end{figure}
\end{center}
\vspace{.5cm}

In the experiments that follow, we set $\nu=.02$ and $r_M=10$. Inspired by (\ref{eig1}) and (\ref{chiap}),
at time $t=0$ we impose:
\begin{equation}\label{init}
u_0(r, \theta ,\phi )=\frac{r^7}{r^4_M} (r_M-r) \Big[\cos (\omega \theta ) \cos (\omega \phi )
+\cos (\omega \theta ) +\cos (\omega \phi )\Big]
\end{equation}
which means that $c_{00}=0$ and $c_{10}=c_{01}=c_{11}$.
\smallskip

We give in Fig.5 the section for $\phi =0$ of the initial velocity field ${\bf v}_0$ evaluated
according to (\ref{defv}). We also show the third component of ${\bf A}$ as prescribed in (\ref{campoas}). 
The level lines of $A_3$ do not exactly envelope the stream lines, but the give however a reasonable
idea of what is going on. The intensity of $u_0$ in (\ref{init}) has been calibrated
to guarantee stability for the time-advancing scheme, also in relation to the magnitude of $\nu$.
The sign of the initial datum influences the behavior of the evolution. With the sign as
in (\ref{init}), the corresponding ${\bf v}_0$ has the rotatory aspect visible in Fig.5.
Like in kind of {\sl driven cavity problem}, there is the tendency to form an internal layer
towards the center of the domain ($\theta =\phi =0$). By switching the sign of ${\bf v}_0$, 
the evolution tends to bring the fluid towards the pyramid vertex (see Fig.12). We prefer the first situation, and
the crucial question is whether this phenomenon may actually determine a deterioration of the regularity of
${\bf v}$ in a finite time.

\vfill

%nsfull3dsoloconti
%nsfull3dsolografica
\begin{center}
\begin{figure}[p!]
\centerline{\hspace{-.5cm}{\includegraphics[width=7.3cm,height=5.8cm]{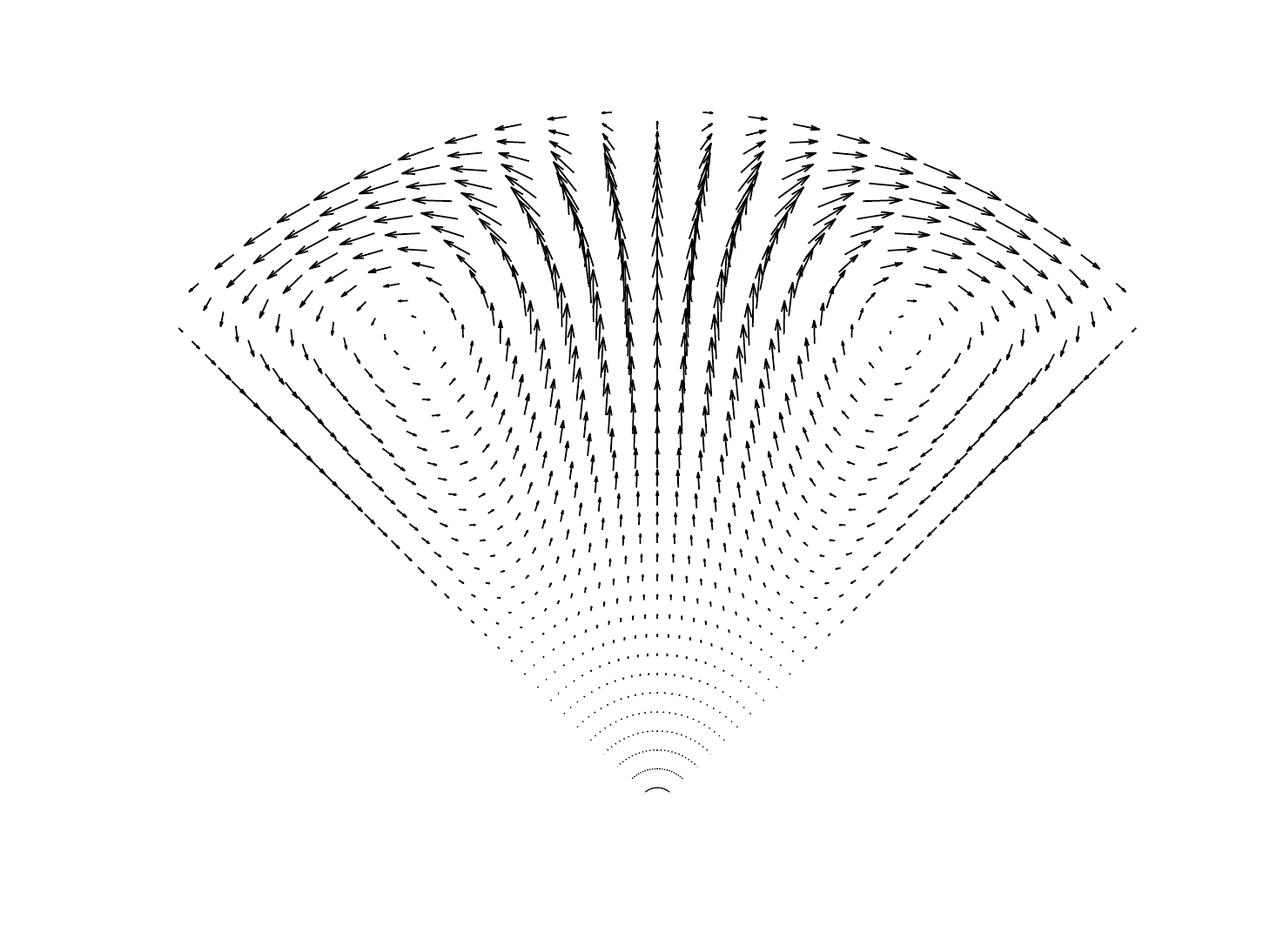}}
\hspace{-.9cm}{\includegraphics[width=7.3cm,height=5.6cm]{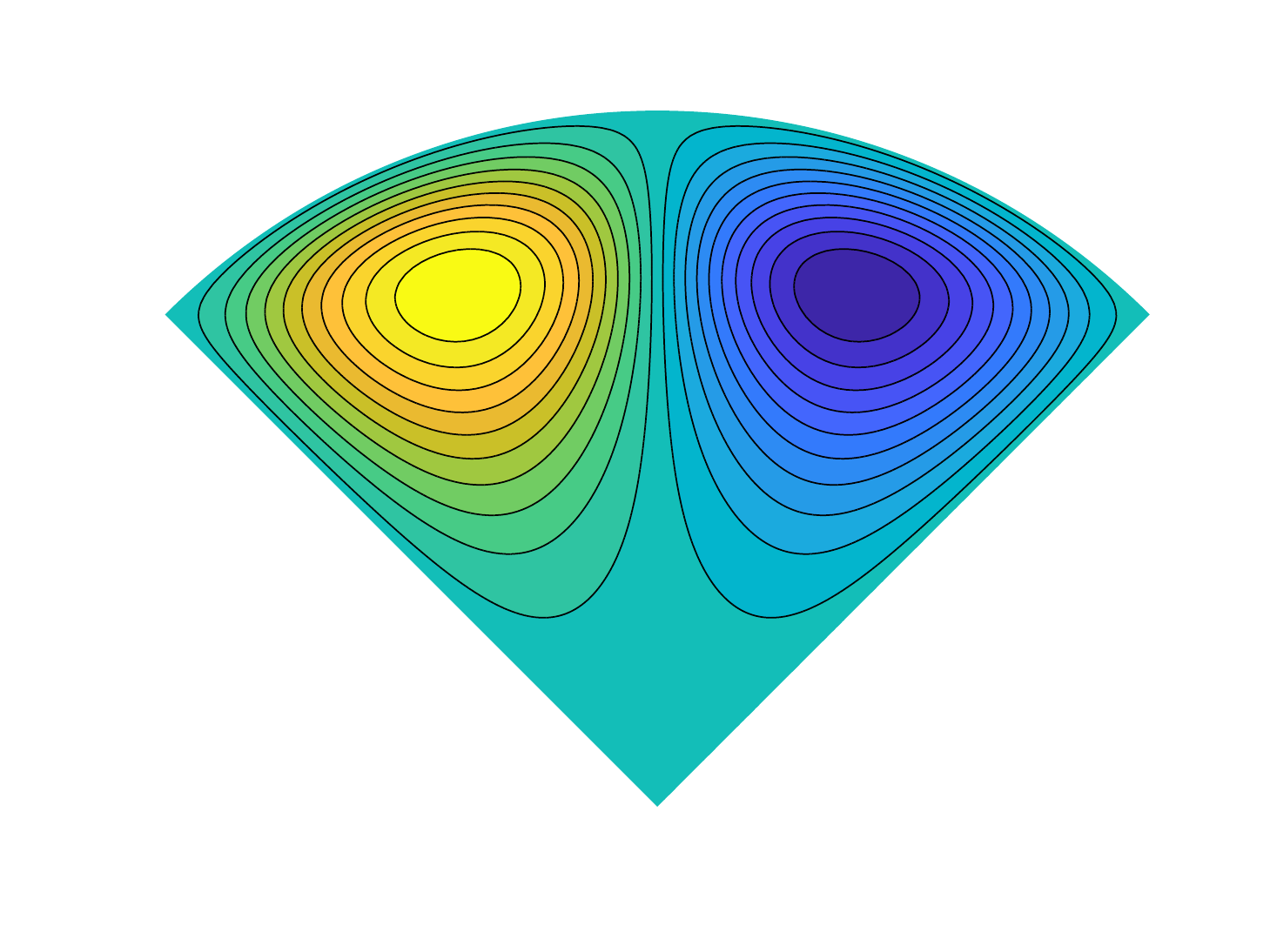}}}
\vspace{-.2cm}
\centerline{\hspace{-.5cm}{\includegraphics[width=7.3cm,height=5.8cm]{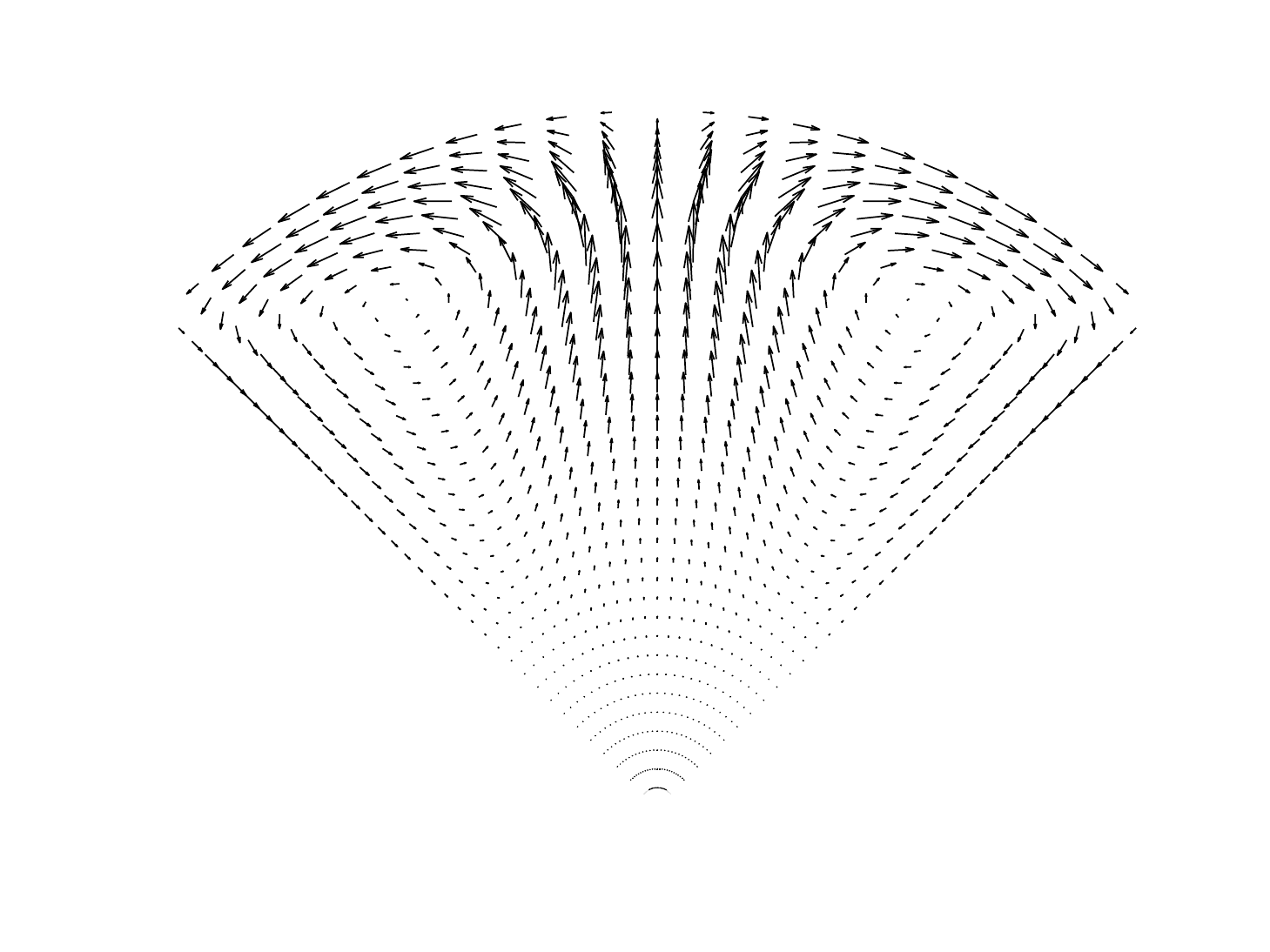}}
\hspace{-.9cm}{\includegraphics[width=7.3cm,height=5.6cm]{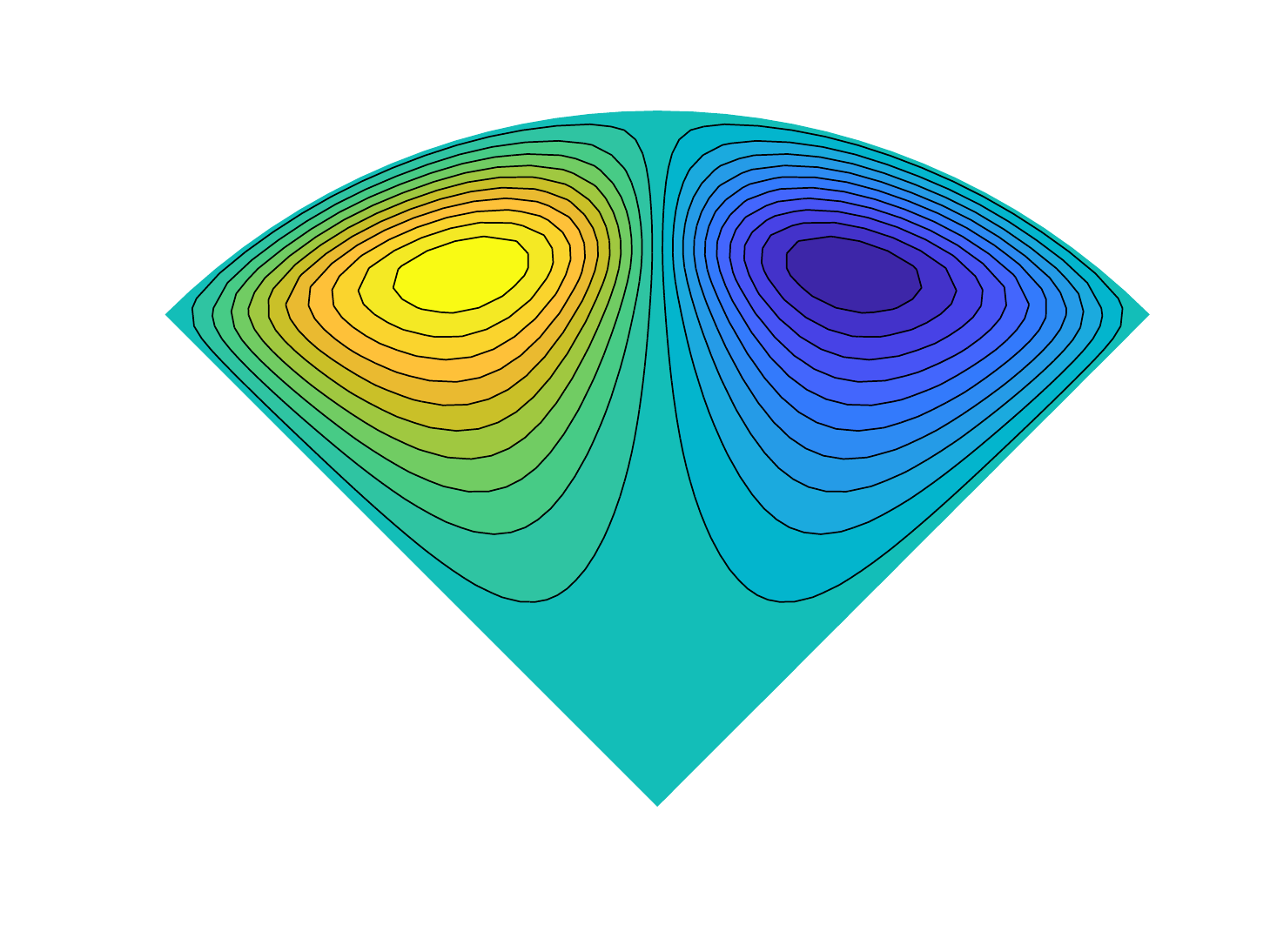}}}
\vspace{-.2cm}
\centerline{\hspace{-.5cm}{\includegraphics[width=7.3cm,height=5.8cm]{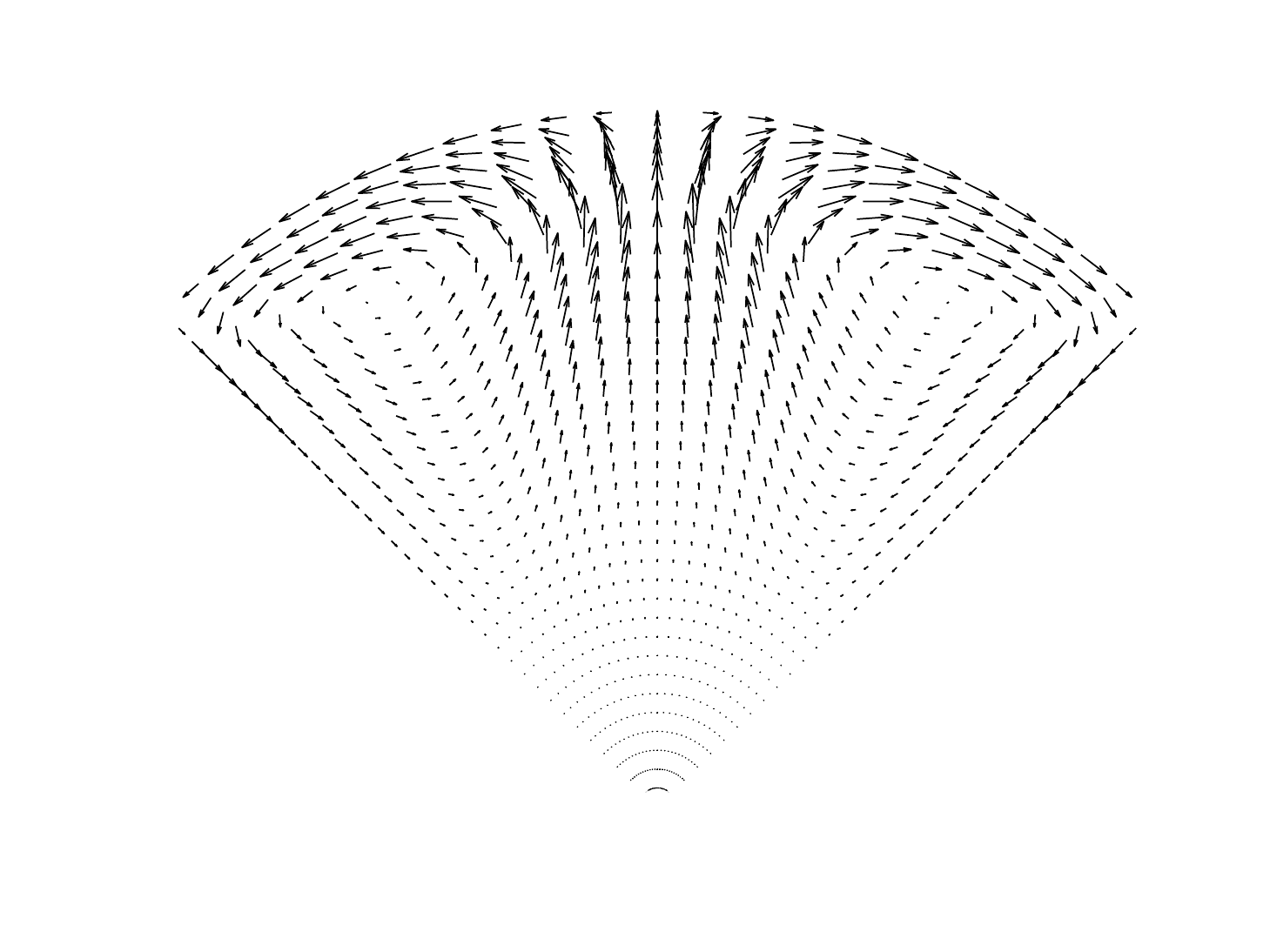}}
\hspace{-.9cm}{\includegraphics[width=7.3cm,height=5.6cm]{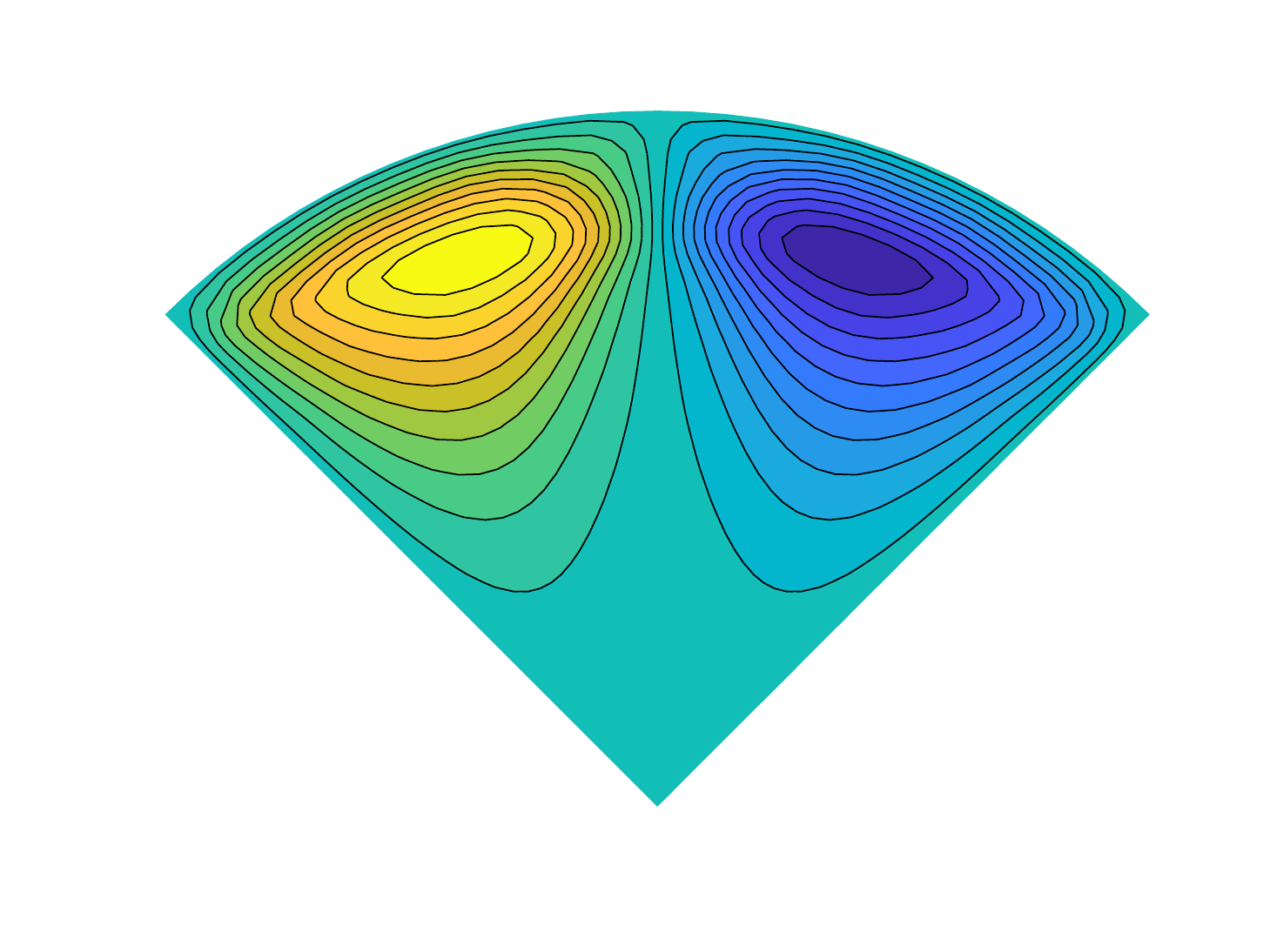}}}
\vspace{-.2cm}
\begin{caption}{\small Successive evolution of the ring sections, starting from the initial data of figure 5.
The pictures are referred to the cut corresponding to $\phi =0$, and the snapshots are
taken at times $t=0.044$, $t=0.077$, $t=T=0.110$, respectively. The Fourier series have been truncated at $N=7$.}
\end{caption}
\end{figure}
\end{center}
%\vspace{-.5cm}

\begin{center}
\begin{figure}[h!]
\centerline{\hspace{-.5cm}{\includegraphics[width=9.4cm,height=7.5cm]{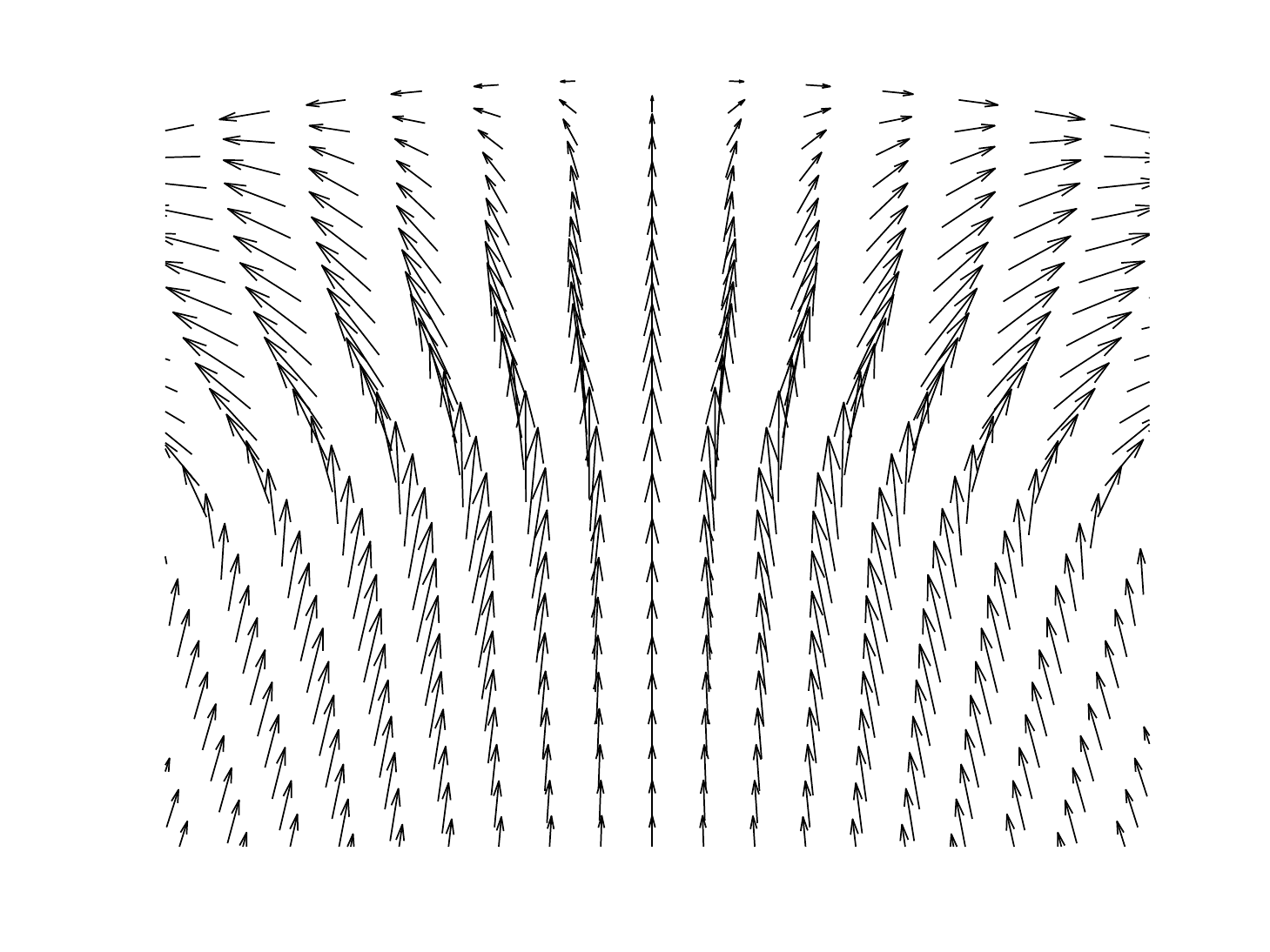}}}
\vspace{-.4cm}
\centerline{\hspace{-.5cm}{\includegraphics[width=9.4cm,height=7.3cm]{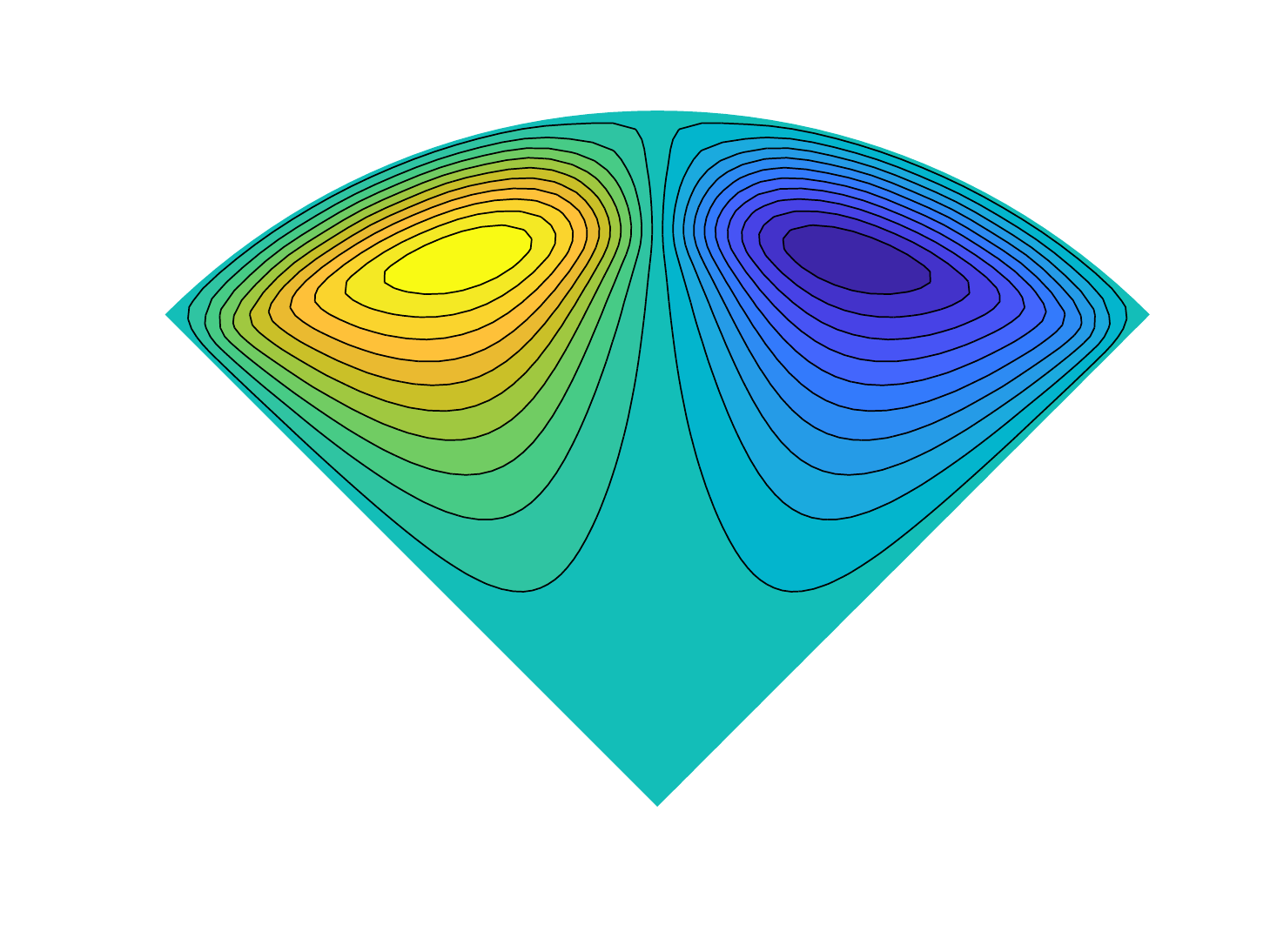}}}
\vspace{-.5cm}
\begin{caption}{\small Section at time $t=0.10$ for $\phi =0$, with an enlargement of the vector field.}
\end{caption}
\end{figure}
\end{center}
%\vspace{-.5cm}

We provide in Fig.6 some snapshots of the section (corresponding to $\phi =0$) of the evolving ring.
In truth, viewed from top (i.e., lying on the square $\Omega$ of the plane $(\theta ,\phi )$), the shape is not exactly
that of a classical rounded ring, but the body is a little elongated in proximity of the four corners.
The situation can be better examined in Fig.7, where an enlargement is provided for the solution at time $t=0.1$
(a bit earlier than the final time of computation). After that time, the evolution continues to be stable and the discrete solution 
remains bounded. The approximated solution has been obtained by truncating
the summations in (\ref{coeffrenbid})-(\ref{coeffrenbidc}) in correspondence to the indexes greater than $N=11$. 
The interval $[0, r_M]=[0,10]$ has been divided
into 73 parts. The $L^2$ norm of the velocity field shows very little variation during the evolution.
However, a decay should be normally observed due to the presence of the viscous term and the numerical diffusion
introduced by the discretization. In Fig.8 we can see the plot of the velocity component $v_1$ in the square 
$[0,r_M]\times [-\pi/4, \pi/4 ]$. Qualitatively, the pictures do not change too much by reducing or increasing
the degrees of freedom. That is true up to a critical time approximately equal to $t=0.1$.

\vspace{.1cm}
\begin{center}
\begin{figure}[h!]
\centerline{{\includegraphics[width=11.8cm,height=8.85cm]{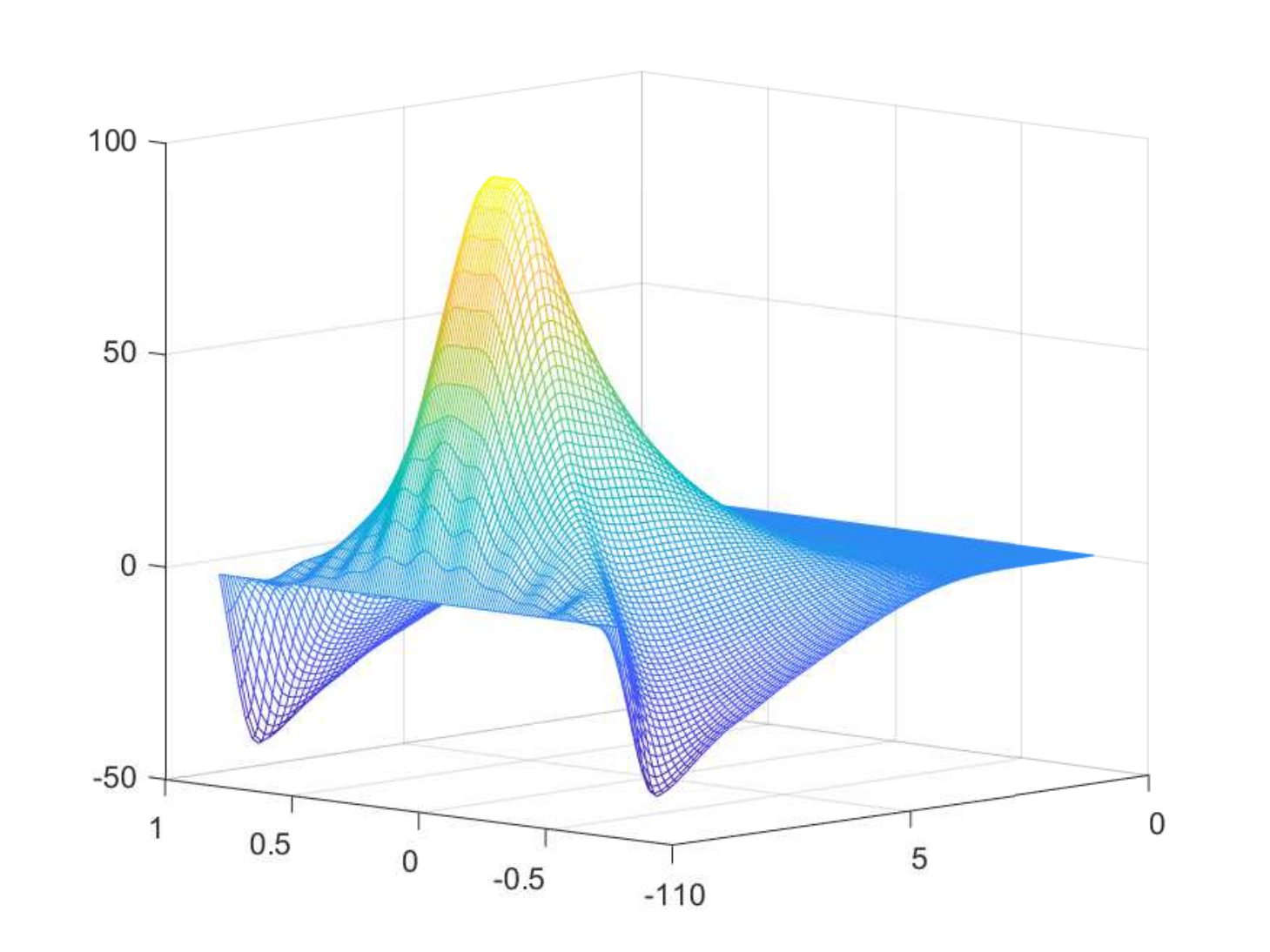}}}
\vspace{-.2cm}
\begin{caption}{\small Plot of the component $v_1$ at time $T=.10$, for $\phi = 0$ and 
$(r,\theta)\in [0,r_M]\times [-\pi/4, \pi/4 ]$. Some wiggles are present at the base. We suspect that they 
are due to the formation of layers at the corner points $(r_M, \pm\pi/4)$.}
\end{caption}
\end{figure}
\end{center}
\vspace{-.2cm}

The sections develop so that the main vortex moves upwards, trying to create a layer in proximity of the upper boundary.
We explain with some mathematical arguments why the vortex tends to be squeezed upwards as time evolves.
If in (\ref{eig1}) and (\ref{eig2}), we consider a higher mode, such as $\cos (k\omega\theta )\cos (k\omega\phi )$,
for $k\geq 1$, the corresponding $\sigma$ in (\ref{chidef}) is now required to satisfy the relation
$\sigma (\sigma +1)=2 k^2\omega^2$. This means that the approximated $\chi$ in (\ref{chiap}) shows a larger value of 
$\sigma$  at the exponent. The maximum of this function is reached at a point $\hat r$ given by:
\begin{equation}\label{chiapmax}
\hat r=\frac{\sigma}{\sigma +1}r_M 
\end{equation}
which approaches $r_M$ from below, as $\sigma$ tends to infinity. Thus, when the cosinus frequency increases,
the corresponding Bessel's function tends to reduce the distance between $\hat r$ and $r_M$. 
Of course, this justification, valid for the linear context, is not fully convincing in the case of the
nonlinear version. More insight comes from examining
Fig.9, where the radial component $v_1$ of the velocity field, as a function of the variable $r$, is shown 
for $\theta =\phi=0$ (the other two components $v_2$ and $v_3$ are zero).
The behavior seems to follow a kind of 1D Burgers equation, where the graph shifts from left to right. 
Up to $t=.09$ everything goes smooth, although the second derivatives tend to grow.
Between $t=.09$ and $t=.10$ there is a change of regime. 

\vspace{.5cm}
\vspace{-.4cm}
\begin{center}
\begin{figure}[h!]
\centerline{\hspace{-.5cm}{\includegraphics[width=9.4cm,height=7.5cm]{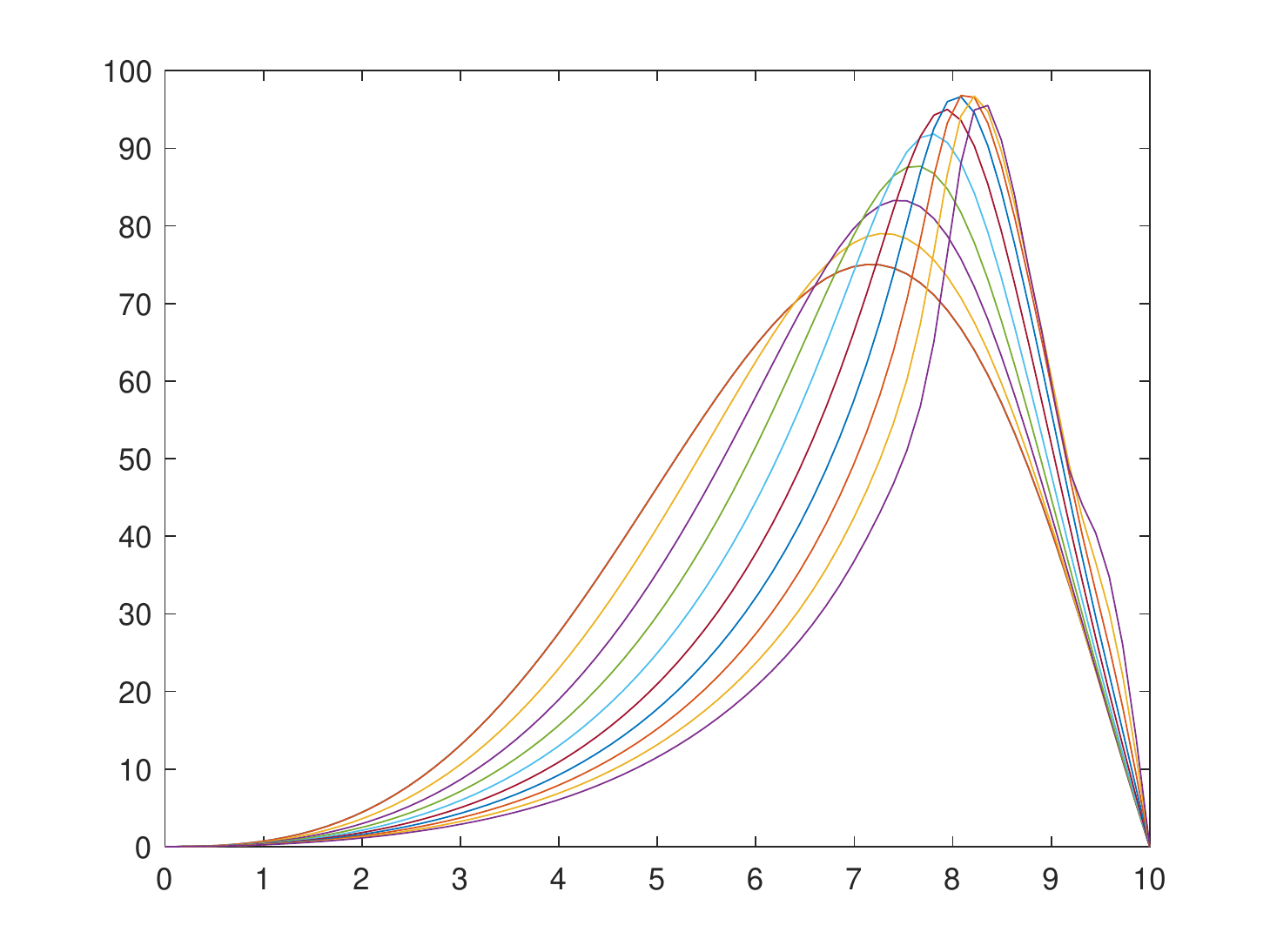}}}
\vspace{-.4cm}
\begin{caption}{\small Behavior of the component $v_1$ for $\theta =\phi =0$, at times $t=.01, .02, \cdots ,.11$.}
\end{caption}
\end{figure}
\end{center}
\vspace{-.2cm}

According to Fig.10, the vector field at the center, which is initially smooth, tends to generate a sort of
jump in the flux rate. This change is transmitted laterally, though one may argue that this is due either to a numerical
effect or to a consequence of the forcing term ${\bf f}({\bf v})$. 
Our guess is that too much fluid tends to accumulate at the center of the ring, and the presence of the upper
boundary cannot dissipate it.  Beyond $t=.10$, the numerical oscillations pollute the outcome (the anomaly is already visible
at the base of the last plot of Fig.8). Going ahead with time, we can reach situations as the one shown in Fig.10, obtained 
with more accurate expansions ($N=15$). These last computations are probably not trustworthy; some strange phenomenon 
is however detectable independently of the degrees of freedom used.

%\vspace{.5cm}
\begin{center}
\begin{figure}[h!]
\centerline{\hspace{-.5cm}{\includegraphics[width=9.4cm,height=7.5cm]{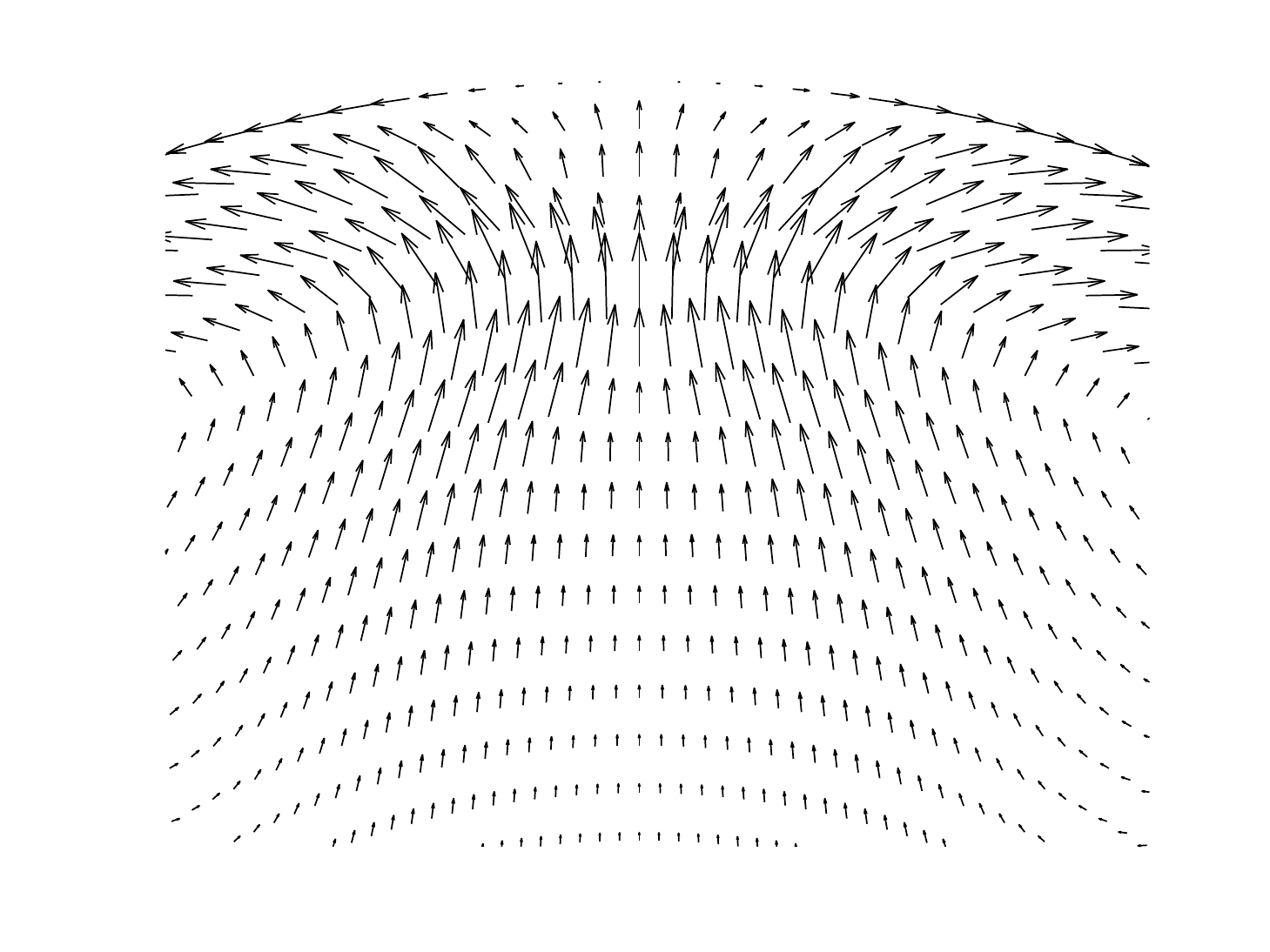}}}
\vspace{-.2cm}
\centerline{\hspace{-.5cm}{\includegraphics[width=9.4cm,height=7.3cm]{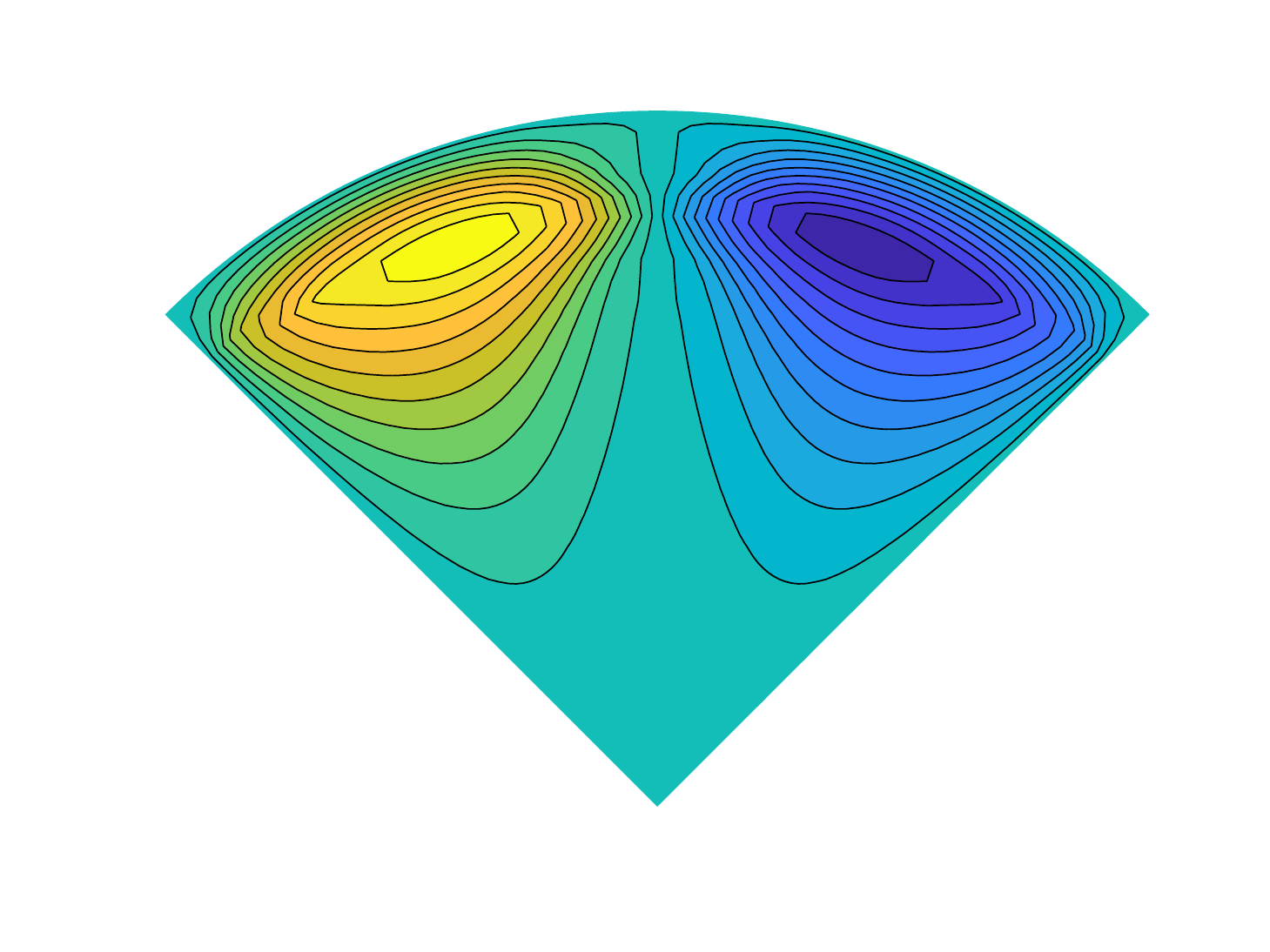}}}
\vspace{-.2cm}
\begin{caption}{\small Section at time $t=.11$ for $\phi =0$. Here the solution begins to be deteriorated,
though it has been calculated with more degrees of freedom than those relative to the previous figures.}
\end{caption}
\end{figure}
\end{center}

These computations are not massive, but rather intensive by the way. Thus, it is quite expensive to perform an
accurate analysis of the real behavior. It is also true that, confirming the presence of a jump of regularity
on the first derivative of the flux, may be practically impossible from the numerical point of view. 
At the critical time something different happens, the solution reaches a kind of steady state and the computation degenerates.  
For sure, we are not in presence of a blowup at infinity or a discontinuity of the field, but maybe of a lack of smoothness.
We suspect that a reliable verification of the facts is only achievable with rather large values of $N$, with
an abrupt growth of the costs for the numerical implementation. We address the reader to section 10 for further results based on
a simplified 2D version of the 3D originating problem. 
\smallskip

%\vspace{-.1cm}
%\begin{center}
%\begin{figure}[h!]
%\centerline{\hspace{-.5cm}{\includegraphics[width=9.4cm,height=7.5cm]{storiaN16Tlungo.pdf}}}
%\vspace{-.4cm}
%\begin{caption}{\small Behavior of the component $v_1$ for $\theta =\phi =0$, up to the critical time and beyond.}
%\end{caption}
%\end{figure}
%\end{center}
%\vspace{-.2cm}

\vspace{-.2cm}
\begin{center}
\begin{figure}[h!]
\centerline{\hspace{-.5cm}{\includegraphics[width=9.4cm,height=7.3cm]{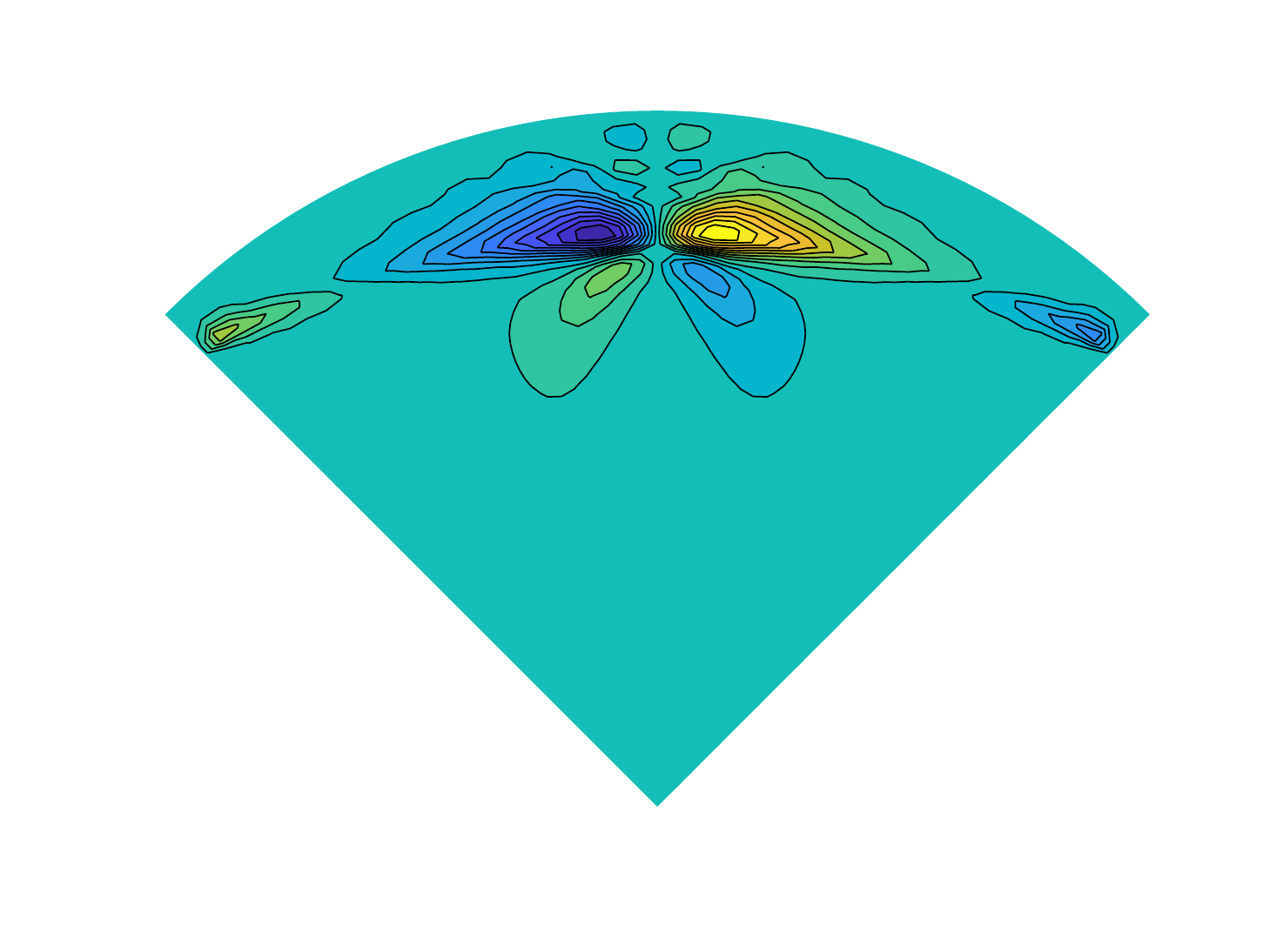}}}
\vspace{-.4cm}
\begin{caption}{\small Section of the function $f_2$ for $\phi =0$ at time $t=.10$.}
\end{caption}
\end{figure}
\end{center}
%\vspace{-.5cm}

As pointed out at the end of the previous section, the pictures presented so far do not reflect
the actual physical behavior of an autonomous velocity field ${\bf v}=(v_1,v_2,v_3)$ simulating a vortex ring. 
There is in fact a forcing term ${\bf f}$,
whose nature depends on the solution itself. According to (\ref{fnoto}) and (\ref{defv}), we have:
$$f_2=-\frac{\Delta \Psi}{r^2}\frac{\partial^2 (r\Phi)}{\partial r \partial \theta} =v_1\left(\frac{\partial v_2}{\partial r}
+\frac{v_2}{r}\right)
$$
\begin{equation}\label{fnotorev}
f_3=-\frac{\Delta \Psi}{r^2}\frac{\partial^2 (r\Phi)}{\partial r \partial \phi}=v_1\left(\frac{\partial v_3}{\partial r}
+\frac{v_3}{r}\right)
\end{equation}
We show in Fig.11 the plot of $f_2$ restricted to the plane $\phi=0$. The snapshot is taken at time $t=.10$ (the same
as in the pictures of Fig.7). Note that, relatively to the section $\theta =0$, $f_2$ is 
identically zero. 
\smallskip

The largest variations are manifested not too far from
the point (denoted by $P$) where $\vert v_1\vert $ reaches its maximum (see Fig.9). It has to be noticed, 
however, that both $f_2$ and
$f_3$ are the results of a multiplication of two terms and that $v_2$ and $v_3$ are identically 
zero for $\theta =\phi =0$. In Fig.11 there are regions where the function undergoes sharp changes, but things do not
seem to be so critical near $P$, where the worst variation should be expected.
\smallskip

We try to reach some heuristic conclusions by introducing the quantity $\rho=\sqrt{\theta^2 +\phi^2}$, and
assuming that at $P$ the function $\Psi$ behaves as $r^\alpha \rho^\beta$ (up to additive and multiplicative
constants), for appropriate values of the parameters
$\alpha$ and $\beta$. In this circumstance, we have the estimates:
$$
v_1 \approx r^{\alpha -1}\rho^{\beta -2} \qquad \quad v_2 \approx  r^{\alpha -1}\theta \rho^{\beta -2} \qquad
\quad v_3 \approx r^{\alpha -1}\phi\rho^{\beta -2}
$$ 
\begin{equation}\label{stime}
\qquad f_2 \approx r^{2\alpha -3}\theta\rho^{2\beta -4} \qquad \quad
f_3 \approx r^{2\alpha -3}\phi\rho^{2\beta -4} \qquad
\end{equation}
If for example we set $\alpha =5/2$ and $\beta =3$, the corresponding $v_1$ belongs
to $H^1(\Sigma )\cap C^0(\Sigma)$ but not to $H^2(\Sigma )$, where $\Sigma =]0,r_M[\times \Omega$.
On the other hand, we note that $f_2 \approx r^2 \theta (\theta^2 +\phi^2)$ and
$f_3 \approx r^2 \phi (\theta^2 +\phi^2)$ are locally smooth functions. This means that we have
room enough to suppose that a regular forcing term may produce a non regular solution, at least
for what concerns the integrability of certain derivatives. Note also that the second and the third
components of the smoothing term $\bar\Delta{\bf v}$ are entirely swallowed by the gradient of pressure.
\smallskip 

By looking for some old references relative to the regularity of Navier-Stokes solutions, we
come out for instance with the following papers: \cite{beirao}, \cite{caffa}, \cite{giga}, \cite{serrin}, \cite{struwe}.
Of course, much more material is available, as a consequence of an intense research activity. In our case, we have 
special type boundary conditions and an uncommon forcing term, therefore it is not easy to find pertinent
results. We leave this kind of analysis to the experts.
We guess that $v_1, v_2, v_3$ may comfortably stay into the space $H^1({\bf R}^3)$ during time evolution. 
The estimates above suggest a possible blowup at the interior of the functional space $H^2({\bf R}^3)$, 
which is just a bit more regular than $C^0({\bf R}^3)$.
Nevertheless, at the moment we have neither theoretical nor practical arguments to confirm this occurrence.
\smallskip

From our experiments it turns out that the role of the viscosity parameter $\nu$ is not
really crucial. It is true that, for relatively large values of $\nu$, the counterparts of
the plots of Fig.9 become smoother. Maybe, in those circumstances, it is just a matter of increasing the intensity
of the initial guess to restore the critical behavior. On the other hand, it is also possible to choose
$\nu=0$, without affecting the stability of the numerical scheme, and obtaining outputs very similar
to those of Fig.9. Perhaps, future theoretical studies may decree that our approach is fruitless
in the analysis of the possible blowup of the solutions of the Navier-Stokes equation.
However, the idea could still have chances to be applied successfully to the analysis of the non-viscous Euler equation.
\smallskip

We spend a few words regarding the possibility of switching the sign of the initial datum (i.e., by replacing 
$u_0$ by $-u_0$ in (\ref{init})). In Fig.12 we see two moments of this evolution. We are quite confident of the fact 
that a sort of singularity is going to be generated at the origin. For instance, it is reasonable to suppose that $v_1$ decays as $r$
when approaching the vertex of the pyramid. In the whole space ${\bf R}^3$, we would get $r=\sqrt{x^2+y^2+z^2}$, which
is not a regular function. On the other hand, by examining the functions $f_2$ and $f_3$ we find out a posteriori that they are affected
by the same pathology. Thus, we should be in the case where a bad forcing term ${\bf f}$ induces the creation of a bad field
${\bf v}$, and this not an interesting discovery.

\vspace{.2cm}
\begin{center}
\begin{figure}[h!]
\centerline{\hspace{-.5cm}{\includegraphics[width=7.3cm,height=5.8cm]{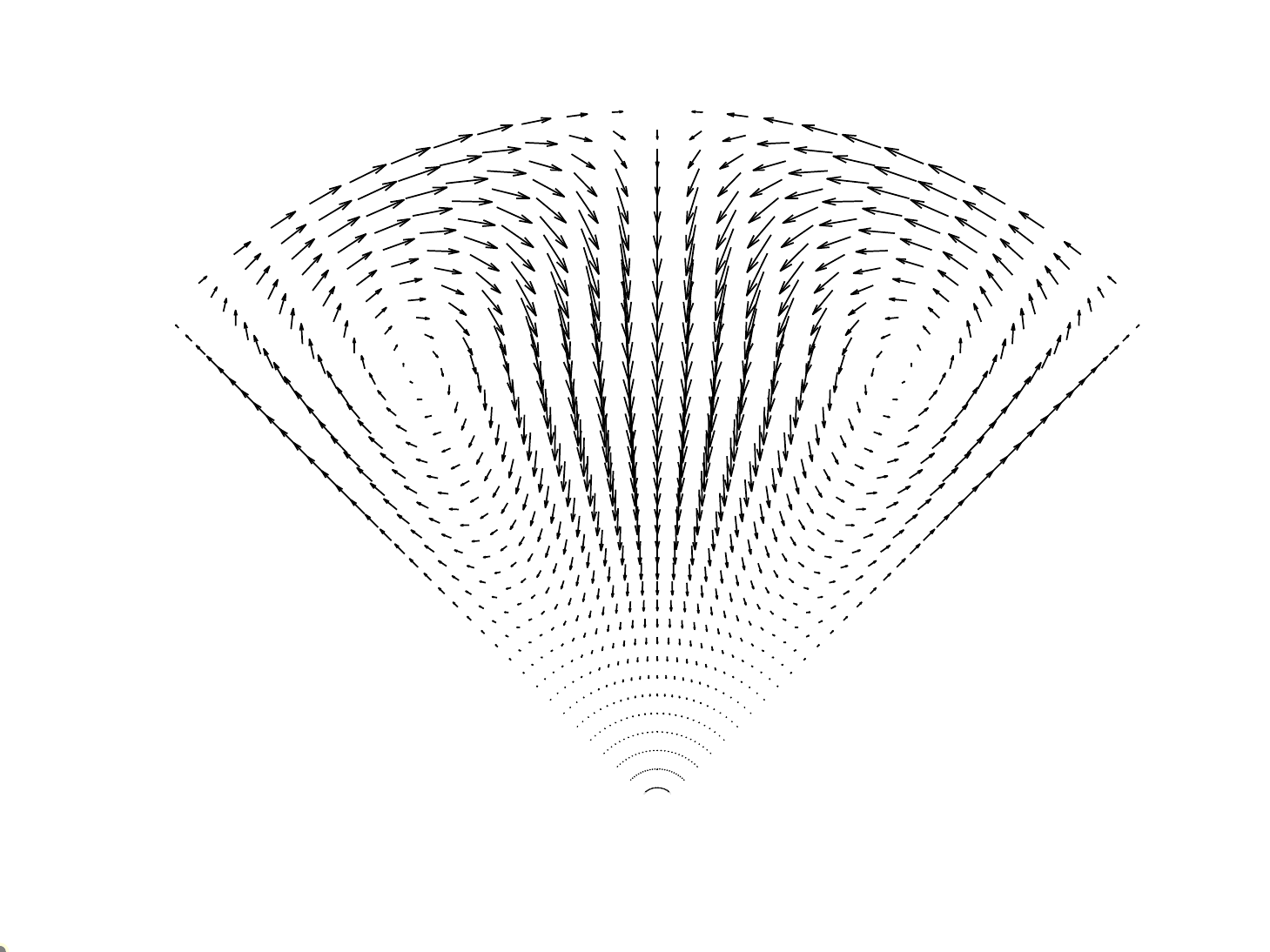}}
\hspace{-.9cm}{\includegraphics[width=7.3cm,height=5.6cm]{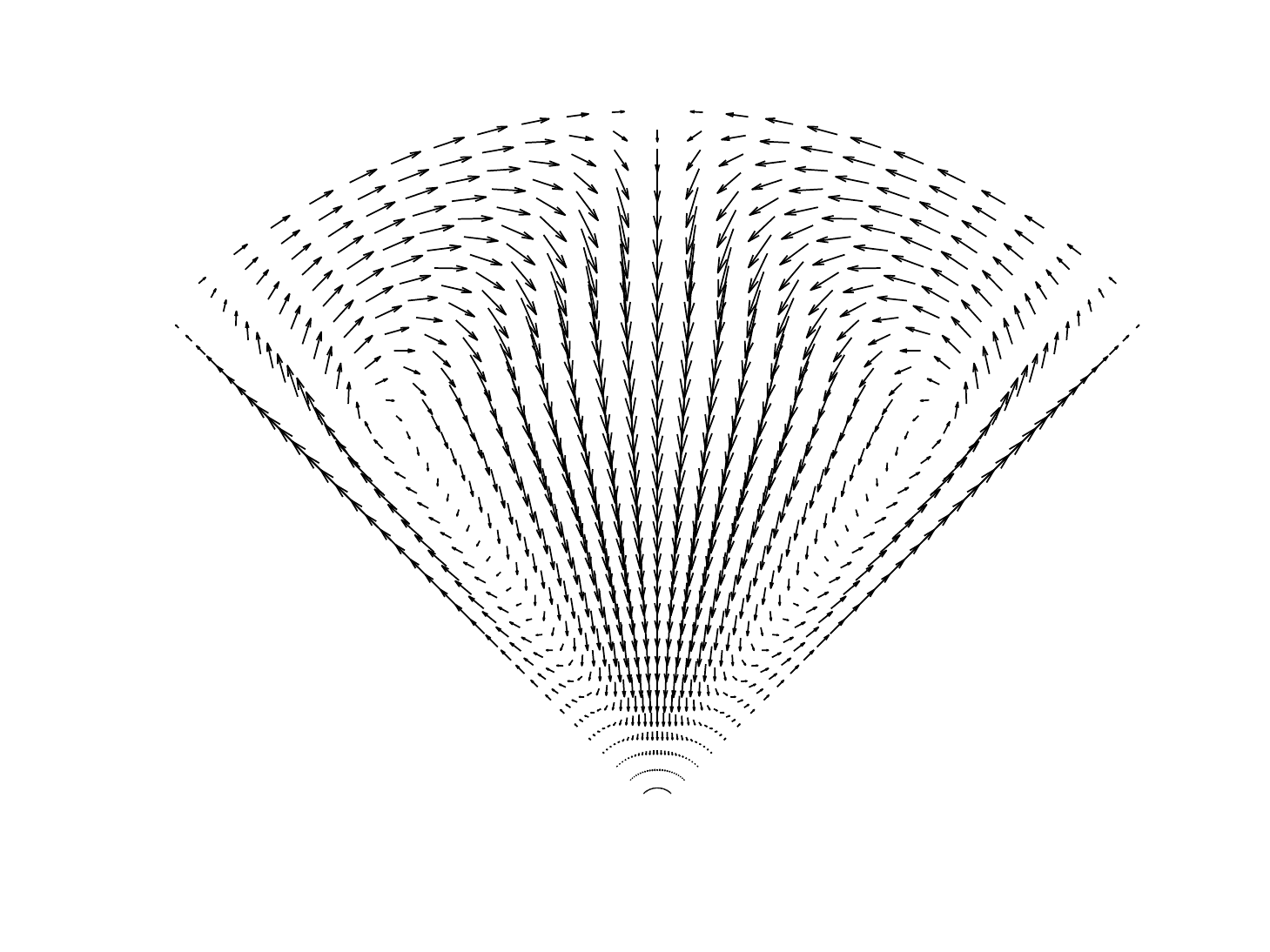}}}
\vspace{-.2cm}
\begin{caption}{\small Field distribution at time $t=0$  and time $t=.11$ when the initial field ${\bf v}_0$ corresponds to $-u_0$ in 
(\ref{init}).}
\end{caption}
\end{figure}
\end{center}
\vspace{.1cm}

As a final remark, we mention the possibility to substitute the pyramid $\Sigma$ with a
cone, and use spherical coordinates $(r,\theta ,\phi)$, where $\theta$ now denotes the
azimuthal angle. In this fashion we require that the expression of the ring does not involve 
the variable $\phi$, so obtaining a 2D problem. After the usual computations, we get:
$$\frac{\partial u}{\partial t}- \nu \left( \frac{1}{r^2\sin\theta } \ \frac{\partial}{\partial\theta}
{\hspace{-.1cm}}\left(\sin \theta \, \frac{\partial u}{\partial \theta} \right)
+\frac{\partial^2 u}{\partial r^2}-\frac{2}{r} \frac{\partial u}{\partial r}+\frac{2u}{r^2}\right)$$
\begin{equation}\label{fineq2sfe}
-\frac{1}{r} \frac{\partial u}{\partial \theta}\, \frac{\partial }{\partial \theta} {\hspace{-.1cm}} 
\left( \frac{\partial \Psi}{\partial r} +\frac{\Psi}{r}\right)+ 
\frac{r}{2}\, \frac{\partial}{\partial r}{\hspace{-.1cm}}\left[\left(\frac{1}{r\sin\theta }\, \frac{\partial}{\partial\theta}
{\hspace{-.1cm}}\left(\sin \theta \, \frac{\partial \Psi}{\partial \theta} \right)\right)^{\hspace{-.1cm}2}\right]=0
\end{equation}
with
\begin{equation}\label{finequsfe}
u=\frac{1}{\sin\theta } \, \frac{\partial}{\partial\theta}
{\hspace{-.1cm}}\left(\sin \theta \, \frac{\partial \Psi}{\partial \theta} \right) +r^2\frac{\partial^2\Psi}{\partial r^2}
+2r\frac{\partial\Psi}{\partial r}
\end{equation}
Unfortunately, if we approach the new set of equations by cosinus Fourier expansions (in order to preserve
Neumann boundary conditions) the formulas are not neat as in (\ref{coeffrenbid}), since there are spurious 
sinus components that cannot be easily handled. Thus, the computational
cost does not decrease significantly. Considering that we are not solving exactly the original
problem and that there are no numerical benefits, we decided not to proceed in this direction.
Nevertheless, in section 10, we examine a simplified version of (\ref{fineq2sfe})-(\ref{finequsfe}).
This surrogate problem will be more affordable from the numerical viewpoint, retaining however some of the main features.

\par\medskip
\section{Comparison with the 2D version}

It is known that the solutions of the 2D navier-Stokes equation preserve indefinitely
their regularity. The 2D version of the example examined so far, corresponds to
four flattened rings, built on triangular slices forming a partition of ${\bf R}^2$. In each single slice,
we work in polar coordinates $(r, \theta)$, or more appropriately in cylindrical coordinates 
$(r, \theta, z)$, where no dependence is assumed 
with respect to the variable $z$. In fact, the $z$-axis, orthogonal to the plane
${\bf R}^2$, is only introduced in order to use the operator {\sl curl}.
 We remind that, in this circumstance, the curl of a vector 
${\bf A}=(A_1, A_2, A_3)$ is determined as follows:
\begin{equation}\label{curlapol}
{\rm curl} {\bf A}=\left( \frac{1}{r}\frac{\partial A_3}{\partial \theta}, \
-\frac{\partial A_3}{\partial r}, \
\frac{\partial A_2}{\partial r}+ \frac{A_2}{r} -
\frac{1}{r}\frac{\partial A_1}{\partial \theta}
\right)
\end{equation}
For a scalar potential $\Psi$, which is function of $t$, $r$ and $\theta$, we define:
\begin{equation}\label{potapol}
{\bf A}=\left( 0, \ 0, \ \frac{\partial \Psi}{\partial \theta}\right)
\end{equation}
\smallskip
By going through the same passages followed for the 3D version, we get:
\begin{equation}\label{vpol}
{\bf v} ={\rm curl} {\bf A}=\left( \frac{1}{r} \frac{\partial^2 \Psi}{\partial \theta^2}, 
 \ -\frac{\partial^2 \Psi}{\partial r\partial \theta},
\ 0\right) =\left(\frac{u}{r}, 0 ,0\right) - \bar\nabla {\hspace{-.1cm}}\left(r \frac{\partial \Psi}{\partial r}
\right)
\end{equation}
where, the new function $u$ is introduced according to the expression:
\begin{equation}\label{upol}
u= \frac{\partial^2 \Psi}{\partial \theta^2}+r\frac{\partial}{\partial r}{\hspace{-.1cm}} 
\left(r \frac{\partial \Psi}{\partial r} \right) = \frac{\partial^2 \Psi}{\partial \theta^2}
+r^2\frac{\partial^2 \Psi}{\partial r^2}+r\frac{\partial \Psi}{\partial r}
\end{equation}

Proceeding with the computations, we have:
\begin{equation}\label{curlupol}
{\rm curl}{\bf v}= \left( 0, \ 0, \ -\frac{1}{r^2}\frac{\partial u}{\partial \theta}
\right)
\end{equation}

\begin{equation}\label{deltaupol}
-\bar \Delta {\bf v} = {\rm curl}({\rm curl}{\bf v})= \left(
-\frac{1}{r^3}\frac{\partial^2 u}{\partial \theta^2} - \frac{\partial}{\partial r}{\hspace{-.1cm}}
\left(r\frac{\partial}{\partial r}{\hspace{-.05cm}}\Big(\frac{u}{r^2}\Big)\right), \ 0,
\ 0 \right)- \bar\nabla q_1
\end{equation}
with 
\begin{equation}\label{q1pol}
q_1=-r\frac{\partial}{\partial r}{\hspace{-.05cm}}\Big(\frac{u}{r^2}\Big)
\end{equation}

Finally, we arrive at the nonlinear term:
\begin{equation}\label{nlpol}
{\bf v}\times {\rm curl}{\bf v}=\left( \frac{1}{r^2}\frac{\partial u}{\partial \theta}
\frac{\partial^2 \Psi}{\partial r\partial \theta} - \frac12 \frac{\partial}{\partial r}
\left(\frac{1}{r}\frac{\partial^2 \Psi}{\partial \theta^2}    
\right)^{\hspace{-.13cm}2}, \ -f_2, \ 0\right) + \bar\nabla q_3
\end{equation}
where
\begin{equation}\label{funcpol}
q_3=\frac12 \left(\frac{1}{r}\frac{\partial^2 \Psi}{\partial \theta^2}\right)^{\hspace{-.13cm}2} 
\qquad f_2=-\frac{1}{r^2}\frac{\partial^2 \Psi}{\partial \theta^2}{\hspace{.06cm}}
\frac{\partial}{\partial r}\hspace{-.1cm}\left(r
\frac{\partial^2 (r\Psi)}{\partial r \partial \theta}\right) 
\end{equation}
\smallskip

By putting all the pieces together, the 2D Navier-Stokes momentum equation
becomes a system in the scalar unknowns $\Psi$ and $u$:
$$\frac{\partial u}{\partial t}- \nu \left[ \frac{1}{r^2}\frac{\partial^2 u}{\partial \theta^2}
+r\frac{\partial}{\partial r}{\hspace{-.1cm}}\left(
r\frac{\partial}{\partial r}{\hspace{-.05cm}}\Big(\frac{u}{r^2}\Big)\right)\right]
$$
\begin{equation}\label{syspol}
-\frac{1}{r}\frac{\partial u}{\partial \theta}\, 
\frac{\partial^2 \Psi}{\partial r\partial \theta} + \frac{r}{2} \frac{\partial}{\partial r}
{\hspace{-.06cm}}\left(\frac{1}{r}\frac{\partial^2 \Psi}{\partial \theta^2}    
\right)^{\hspace{-.13cm}2} =0
\end{equation}
that can be equivalently written as:
$$\frac{\partial u}{\partial t}- \nu \left( \frac{1}{r^2}\frac{\partial^2 u}{\partial \theta^2}
+\frac{\partial^2 u}{\partial r^2}- \frac{3}{r}\frac{\partial u}{\partial r}+ \frac{4u}{r^2}\right)
$$
\begin{equation}\label{syspol2}
+\frac{1}{r}\left[ -\frac{\partial u}{\partial \theta}\, 
\frac{\partial }{\partial \theta}{\hspace{-.13cm}} \left(\frac{\partial \Psi}{\partial r}\right)+ 
\frac{\partial^2 \Psi}{\partial \theta^2}\, \frac{\partial^2 }{\partial \theta^2}{\hspace{-.13cm}}\left(
\frac{\partial \Psi}{\partial r}-\frac{\Psi}{r}\right)\right]=0
\end{equation}
in order to be similar to (\ref{fineq2}). We recall that $u$ is defined in (\ref{upol}).
\smallskip

Let us remark that the last equations have nothing in common with (\ref{fineq2sfe}) and (\ref{finequsfe}).
In the 3D version, defined on a cone, the flow comes from all directions and concentrates on the vertical 
axis. The section of the cone does not correspond to the slice of the 2D version, where the fluid only 
arrives from left or right. This is probably why the 3D version of the Navier-Stokes equations is more vulnerable to
an overcrowding of the fluid in certain areas, giving rise to an exceptional increase of pressure.  
\smallskip

As before, Neumann type boundary conditions are assumed for both $u$ and $\Psi$, i.e.:
$(\partial u /\partial\theta )(\pm \pi/\omega )= (\partial \Psi /\partial\theta )(\pm \pi/\omega )=0$,
for all $r>0$. The two functionals in (\ref{elle1}) an (\ref{elle2}) now become:
\begin{equation}\label{elle1p}
{\cal L}_1u =  \frac{1}{r^2}\frac{\partial^2 u}{\partial \theta^2}
+\frac{\partial^2 u}{\partial r^2}- \frac{3}{r}\frac{\partial u}{\partial r}+ \frac{4u}{r^2}
\end{equation}
\begin{equation}\label{elle2p}
{\cal L}_2\Psi = \frac{\partial^2 \Psi}{\partial \theta^2}
+r^2\frac{\partial^2 \Psi}{\partial r^2}+r\frac{\partial \Psi}{\partial r}
\end{equation}
By playing with the lowest order eigenmodes:
\begin{equation}\label{eig1p}
u_0(r,\theta, \phi)= -\gamma^2 r^2\chi (r) \cos (\omega\theta )
\qquad \qquad
\Psi_0(r,\theta, \phi)= \chi (r) \cos (\omega\theta )
\end{equation}
this time we discover that:
\begin{equation}\label{elle10p}
{\cal L}_1u_0 = -\gamma^2 u_0 \qquad \qquad {\cal L}_2\Psi_0 = u_0
\end{equation}
provided $\chi (r) =\ J_{\sigma } (\gamma r)$, with $\sigma=\omega$.
For $\omega =4$, the first nontrivial zero of the Bessel's function $J_4$ is 7.58.
\smallskip

We run some numerical experiments by setting $\nu=.02$, $r_M=8$ and $T=.21$. 
At time $t=0$ we impose $u_0=\pm(r^4 /r_M) (r_M-r) \cos (\omega \theta )$.
The Fourier expansions are truncated at $N=20$. The plots of Fig.13 show
the evolution of $v_1$ along the axes $\theta =0$ and $\theta =\pi /4$. Comparing with Fig.9, the
transition looks smoother and the effects of dissipation are more prominent. 
However, it has to be remembered that the role of the forcing term ${\bf f}$
(that implicitly depends on the solution itself) may alter the capacity to judge
what is really happening.

%nspol2dsoloconti
%nspol2dsolografica
\begin{center}
\begin{figure}[h!]
\centerline{{\includegraphics[width=6.8cm,height=5.2cm]{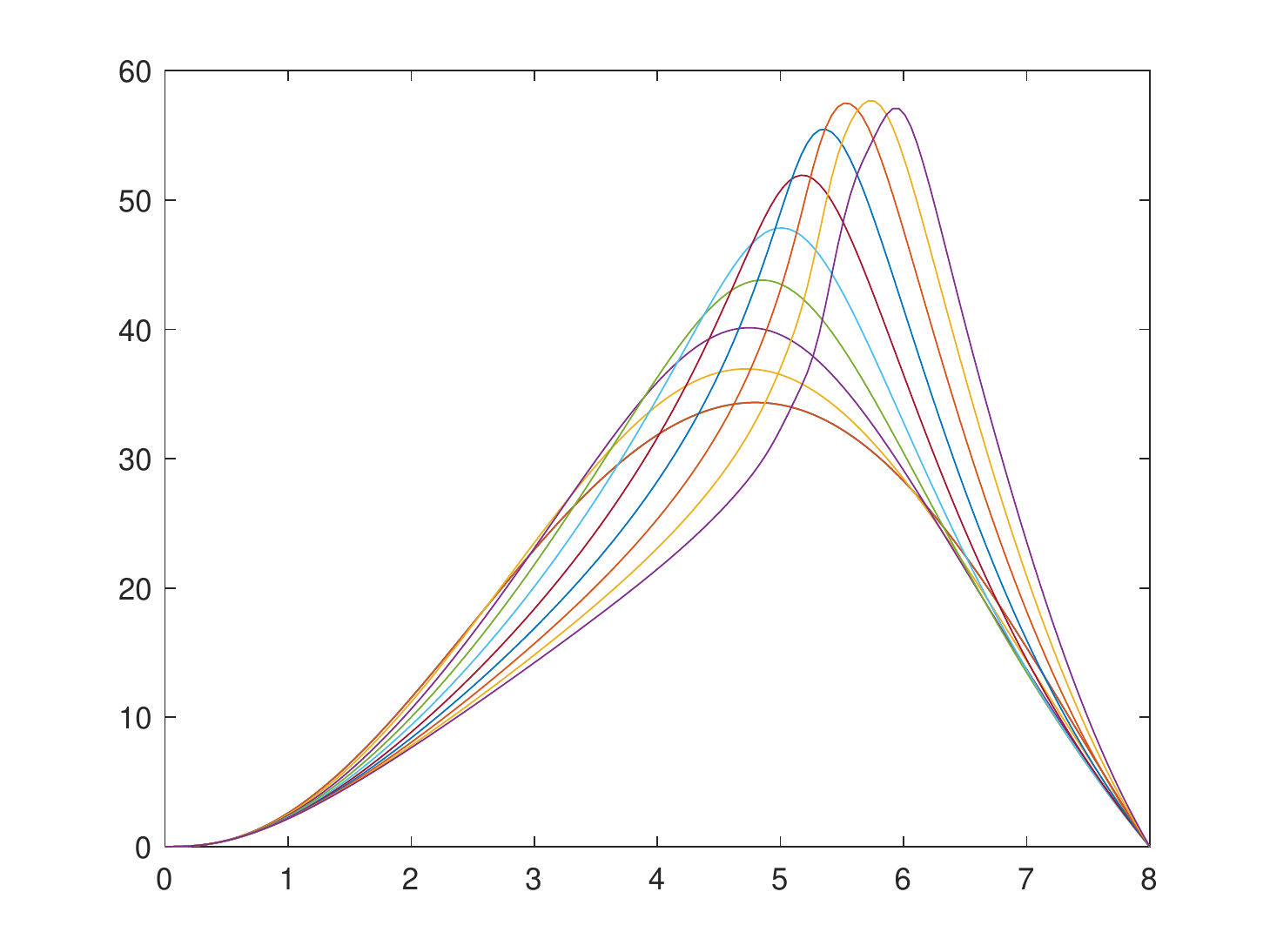}}
\hspace{-.6cm}{\includegraphics[width=6.8cm,height=5.2cm]{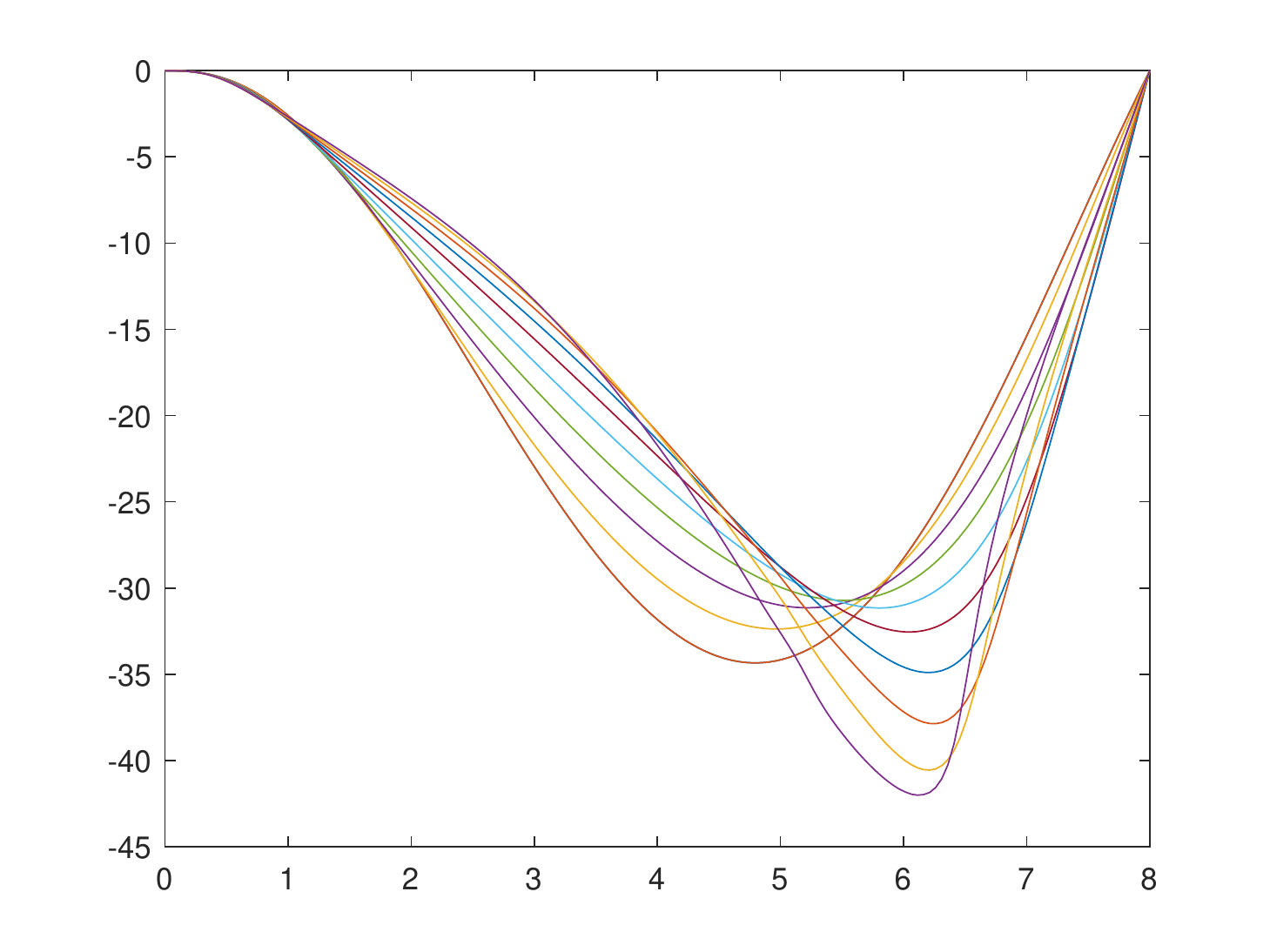}}}
%a centro, b lati
\vspace{-.2cm}
\begin{caption}{\small Behavior of the component $v_1$ for $\theta = 0$ (left) and $\theta =\pi /4$ (right), for 
equispaced time instants in the interval  $[0,T]=[0,.21]$.}
\end{caption}
\end{figure}
\end{center}
\vspace{-.2cm}

The plot of the velocity component $v_1$ on the rectangle $[0,r_M]\times [-\pi /4, \pi/4]$, at time $T=.21$ 
is provided in Fig.14. Everything looks pretty smooth. Once again, it is to be pointed out that the flattened 2D version
of the four rings has not at all the flavor of the original 3D counterpart. Indeed, referring to Fig.15, the amount
of fluid flowing outward along the segment $S_1$ corresponds to a shift of the vortexes towards the external boundary
(as also noticed in the 3D example). By inverting the sense of rotation, we observe a similar effect, as testified
(after a 45 degrees rotation) by the two vortexes separated by the segment $S_2$, where the fluid moves inward. 
Thus, a suitable twist of the whole apparatus amounts to an inversion of the velocity arrows. This is not true
in the three dimensional case, where the quantity of flow concentrating at the center of each ring is far more intense than
that coming from all around. In the 3D version, it is not possible to modify the polarity of the vortexes by a mere 45 
degrees rotation of the entire setting (as testified by the difference between the displacements of Fig.6 and Fig.12).

\vfill

\begin{center}
\begin{figure}[h!]
\centerline{{\includegraphics[width=11.8cm,height=8.85cm]{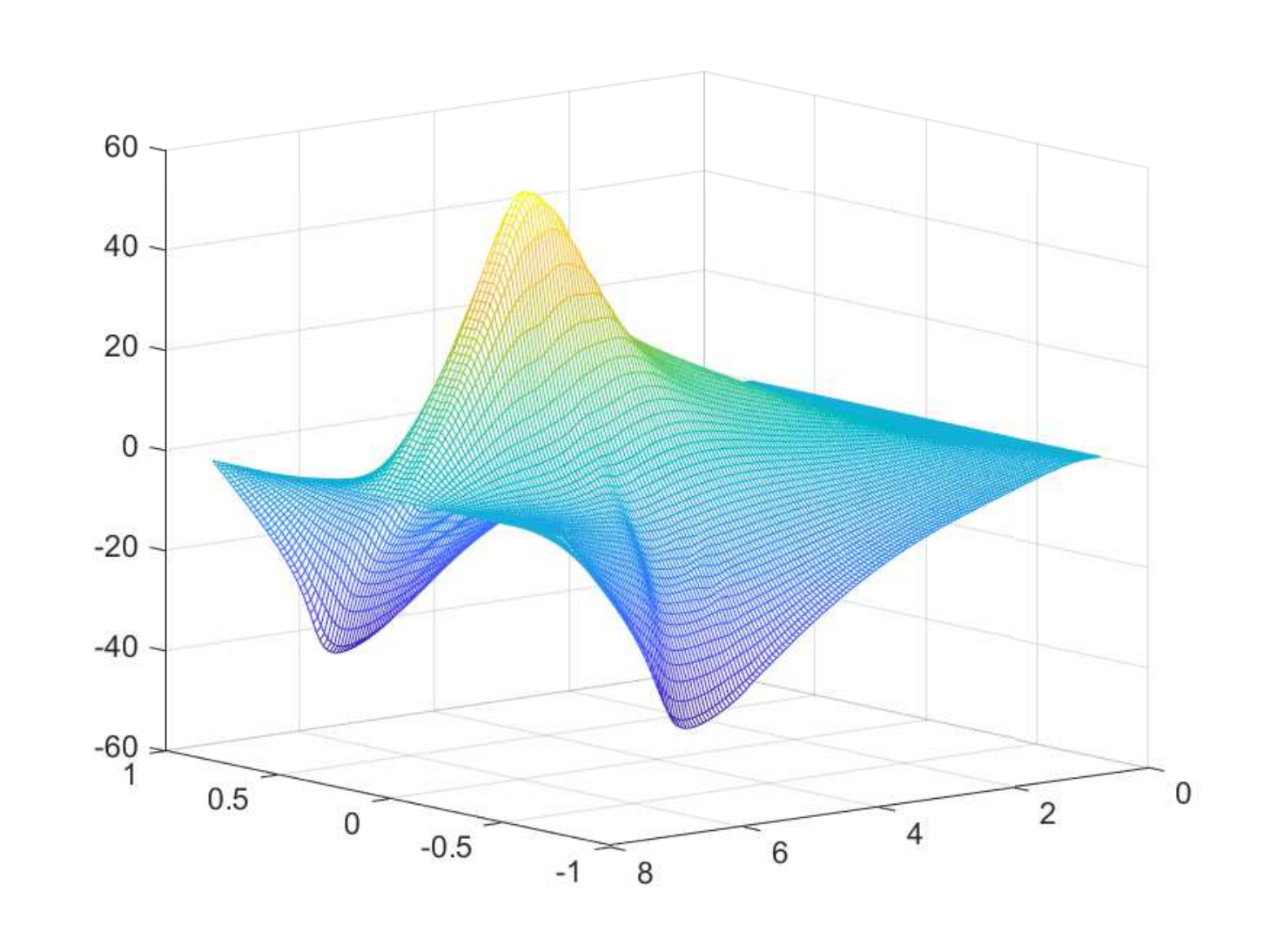}}}
\vspace{-.2cm}
\begin{caption}{\small Plot of the component $v_1$ at time $T=.21$, for  $(r,\theta)\in [0,8]\times [-\pi/4, \pi/4 ]$.}
\end{caption}
\end{figure}
\end{center}
\vspace{-1.cm}

\begin{center}
\begin{figure}[h!]
\centerline{{\includegraphics[width=5cm,height=7.4cm]{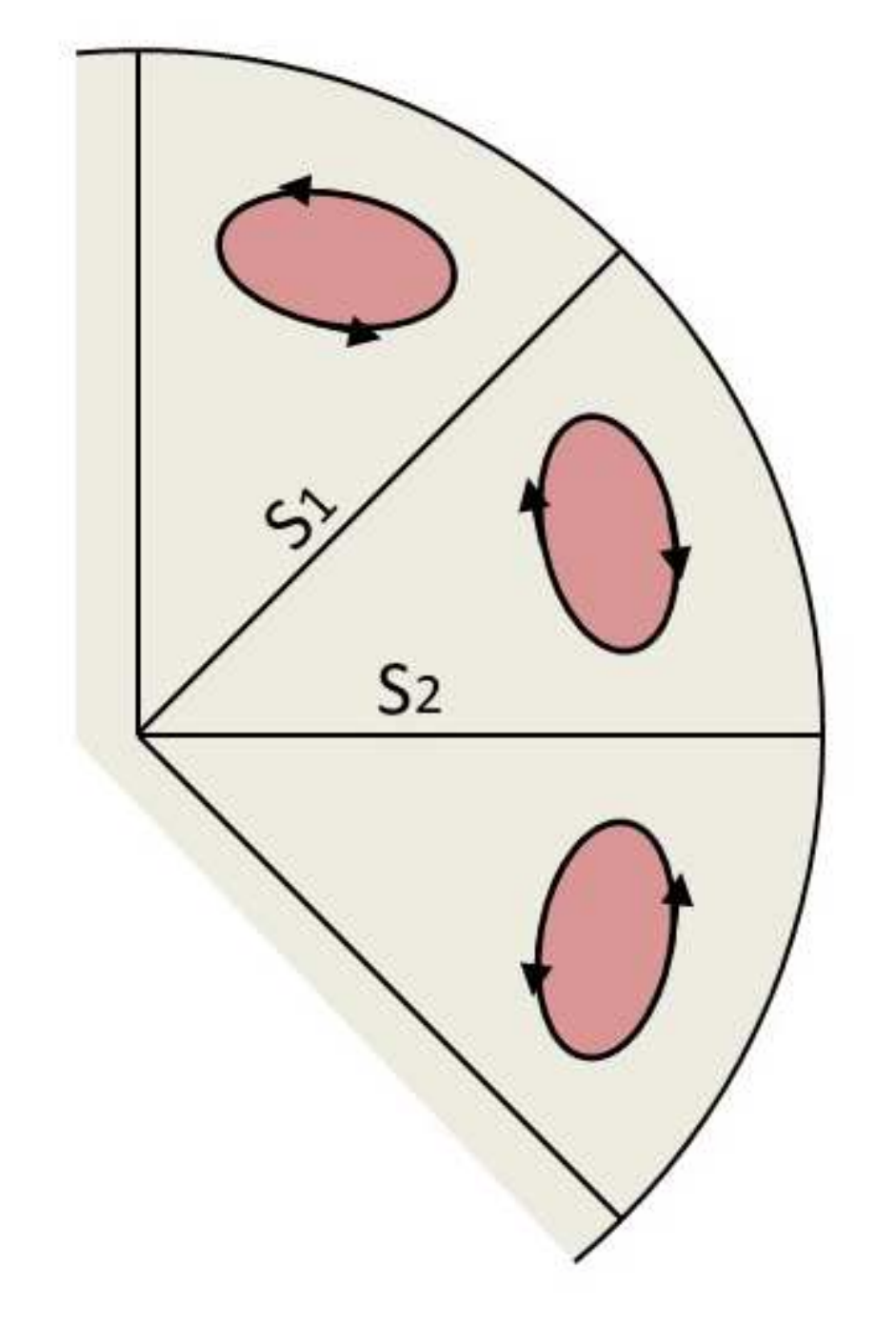}}}
\vspace{-.2cm}
\begin{caption}{\small In the 2D version, a 45 degrees rotation of the entire setting is equivalent to
switch the sense of rotation of the vortexes.}
\end{caption}
\end{figure}
\end{center}
%\vspace{-.2cm}

\vfill

\par\medskip
\section{A simplified model for the cone}

At the end of section 8, we introduced the equations (\ref{fineq2sfe})-(\ref{finequsfe}). 
Defined on a three-dimensional cone, they just make use of the two variables $r$ and $\theta$.
In order to develop a cheap numerical code for the calculation of their solutions, we
introduce the following approximation:
\begin{equation}\label{aptheta}
\frac{1}{\sin\theta } \, \frac{\partial}{\partial\theta}
{\hspace{-.1cm}}\left(\sin \theta \, \frac{\partial \Psi}{\partial \theta}\right)=
 \frac{\partial^2\Psi}{\partial\theta^2}+\frac{\cos\theta}{\sin\theta}\frac{\partial \Psi}{\partial\theta}
\ \approx \ 2  \frac{\partial^2\Psi}{\partial\theta^2}
\end{equation}
which is valid for small $\theta$. In this way we concentrate our attention on the central axis of the cone.
Meanwhile, we open the possibility of implementing Fourier cosinus expansions in an easy fashion.
\smallskip

First of all, the expression of the velocity field takes the form:
\begin{equation}\label{vpols}
{\bf v} =\left( \frac{2}{r} \frac{\partial^2 \Psi}{\partial \theta^2}, 
 \ -\frac{\partial^2 \Psi}{\partial r\partial \theta}-\frac{1}{r}\frac{\partial \Psi}{\partial\theta},
\ 0\right) 
\end{equation}
successively, the equations are modified as follows:
\begin{equation}\label{upols}
u= 2\frac{\partial^2 \Psi}{\partial \theta^2}
+r^2\frac{\partial^2 \Psi}{\partial r^2}+r\frac{\partial \Psi}{\partial r}
\end{equation}

$$\frac{\partial u}{\partial t}- \nu \left( \frac{2}{r^2}\frac{\partial^2 u}{\partial \theta^2}
+\frac{\partial^2 u}{\partial r^2}- \frac{2}{r}\frac{\partial u}{\partial r}+ \frac{2u}{r^2}\right)
$$
\begin{equation}\label{syspols}
+\frac{1}{r}\left[ -\frac{\partial u}{\partial \theta}\, 
\frac{\partial }{\partial \theta}{\hspace{-.13cm}} \left(\frac{\partial \Psi}{\partial r}+\frac{\Psi}{r}\right)+ 
4\frac{\partial^2 \Psi}{\partial \theta^2}\, \frac{\partial^2 }{\partial \theta^2}{\hspace{-.13cm}}\left(
\frac{\partial \Psi}{\partial r}-\frac{\Psi}{r}\right)\right]=0
\end{equation}

%nsproiet2dsoloconti
%nsproiet2dsolografica
\vspace{.2cm}
\begin{center}
\begin{figure}[h!]
\centerline{{\includegraphics[width=6.8cm,height=5.2cm]{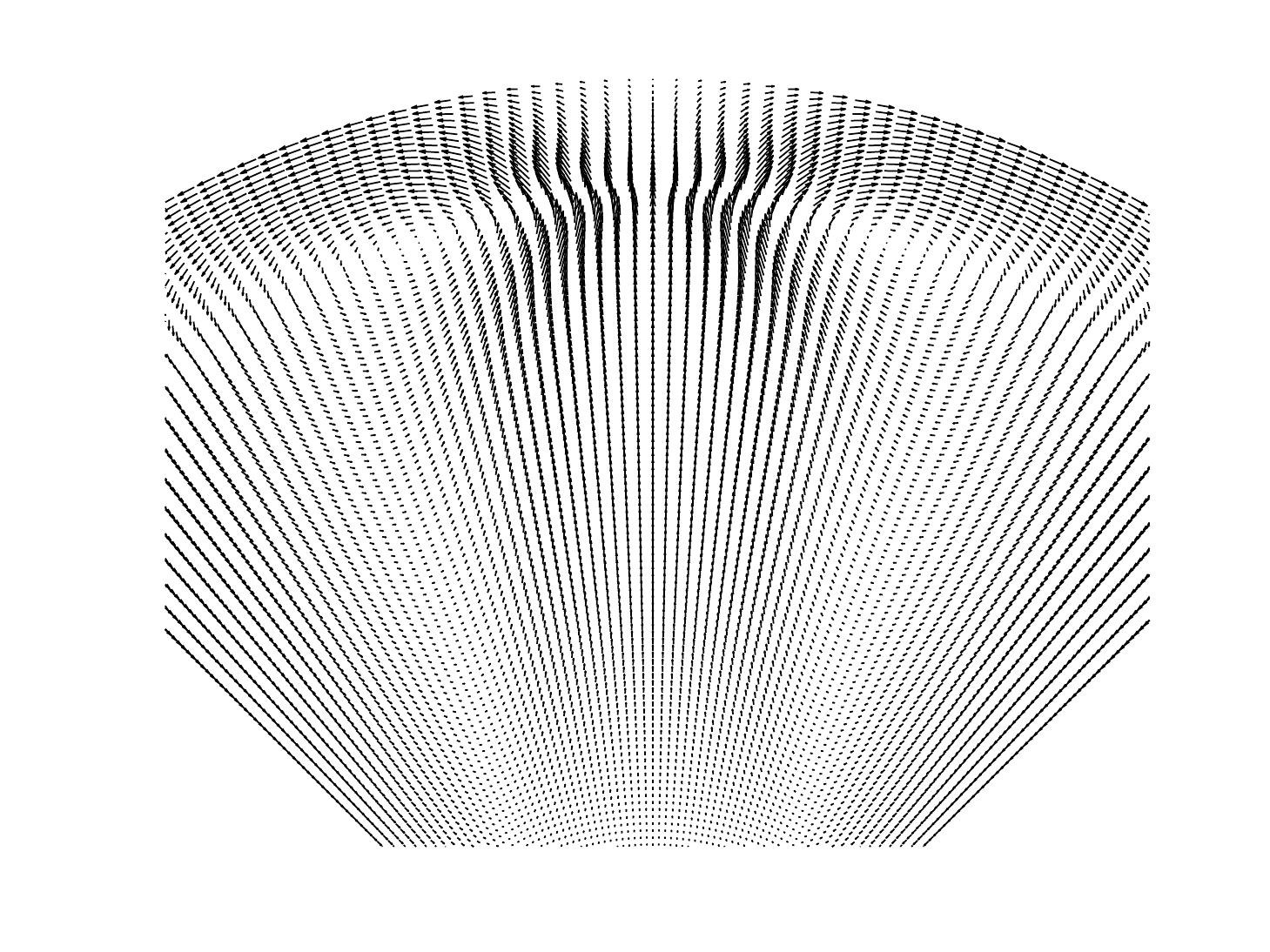}}
\hspace{-.6cm}{\includegraphics[width=6.8cm,height=5.2cm]{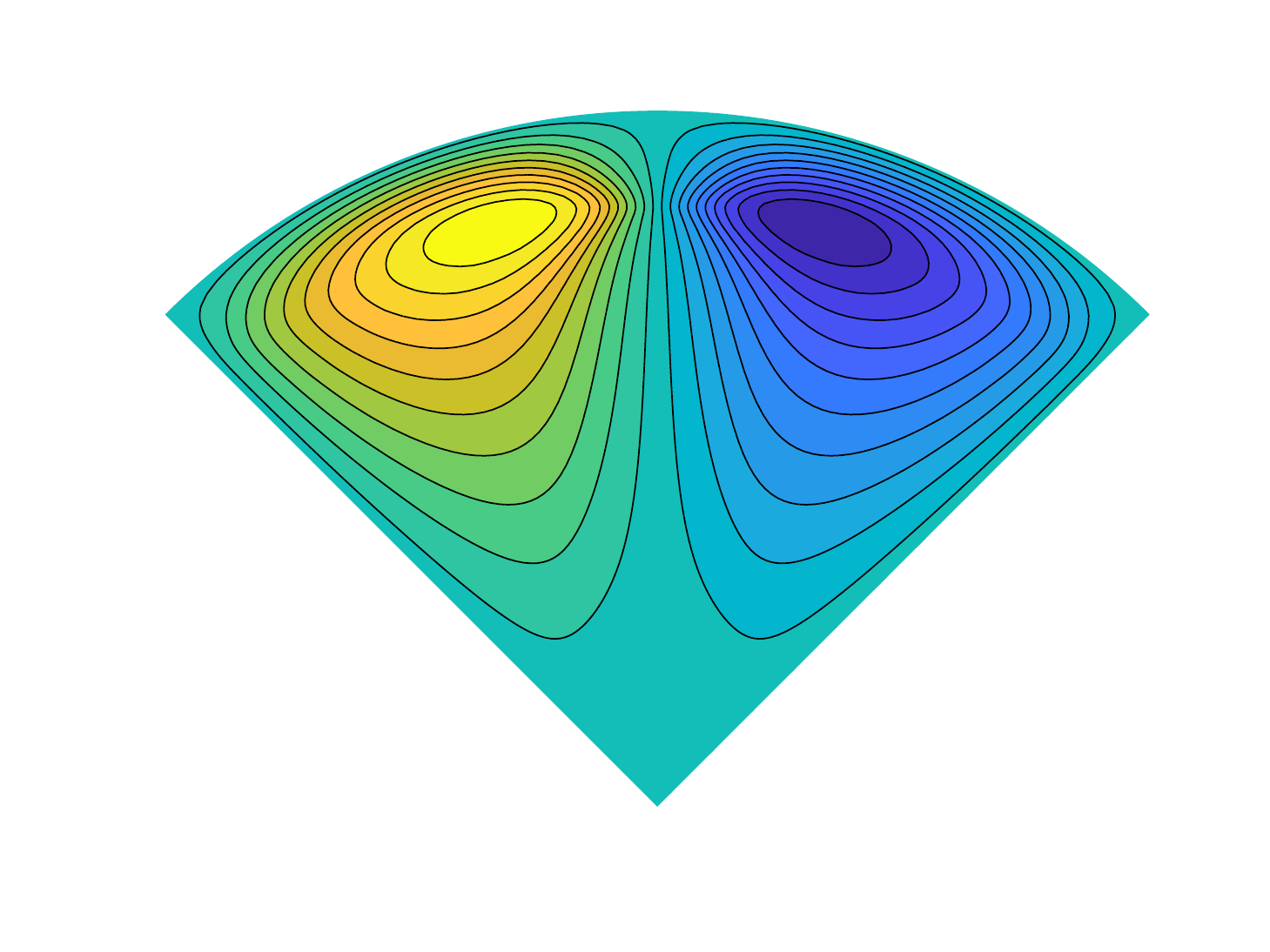}}}
%a centro, b lati
\vspace{-.2cm}
\begin{caption}{\small  Sections for $\theta = 0$ at time $T=.4$. We see an enlargement of the vector
field ${\bf v}$ (left) and the level lines of the function $\partial\Psi /\partial \theta$ (right).}
\end{caption}
\end{figure}
\end{center}
\vspace{.2cm}

The numerical code is the same as the one taken into account in the previous section. The results are however
rather different. We studied the behavior in the time interval $[0,T]=[0,.4]$, with $r_M=10$, $\nu=.02$ and the initial condition
$u_0=(r^7/r_M^4)(r_M-r)\cos \omega\theta$. Regarding the outcome, we refer to figures 16, 17, 18, where in the experiments the series have 
been truncated for $N>18$.

\vspace{-.2cm}
\begin{center}
\begin{figure}[h!]
\centerline{{\includegraphics[width=6.8cm,height=5.2cm]{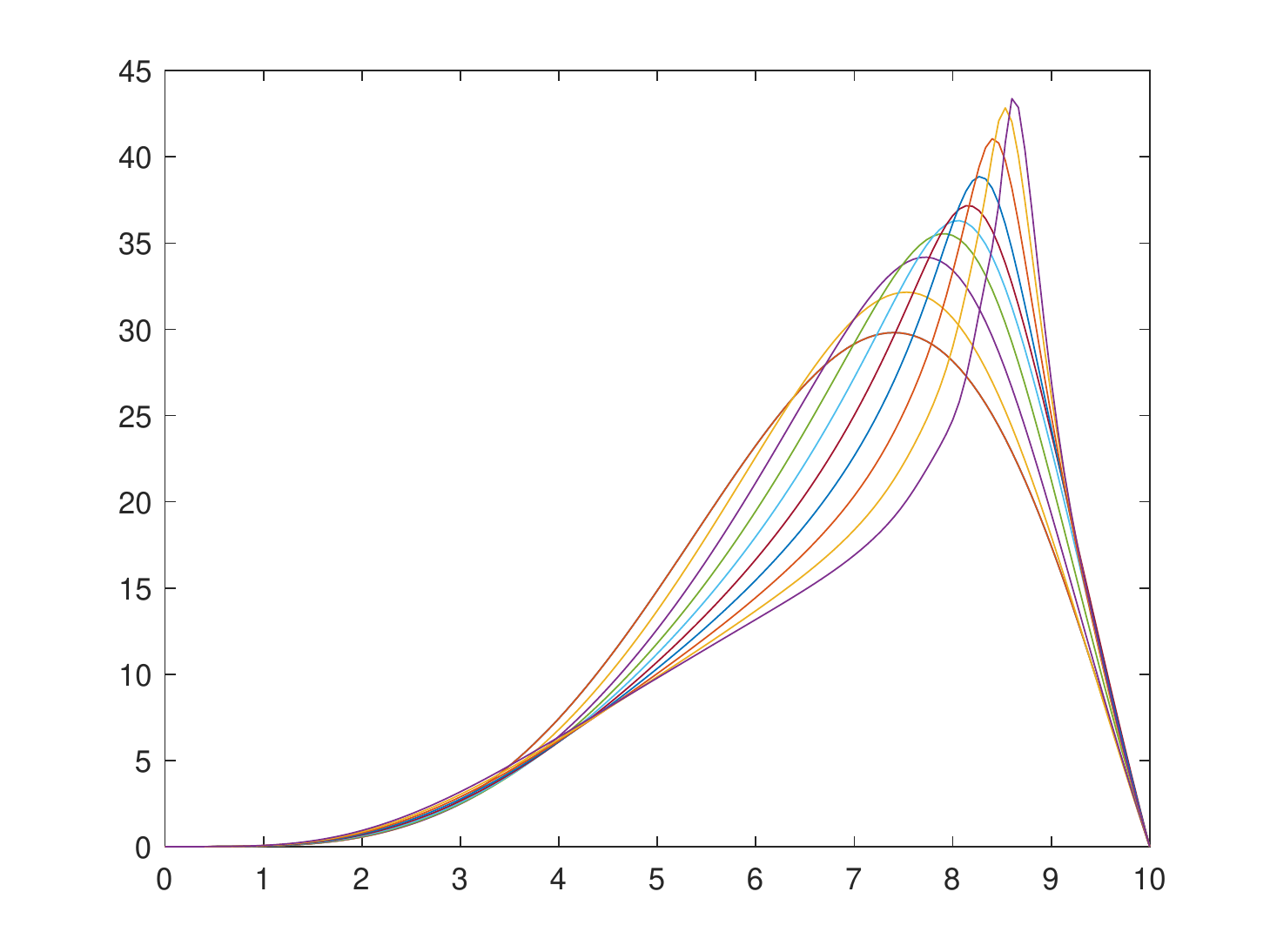}}
\hspace{-.6cm}{\includegraphics[width=6.8cm,height=5.2cm]{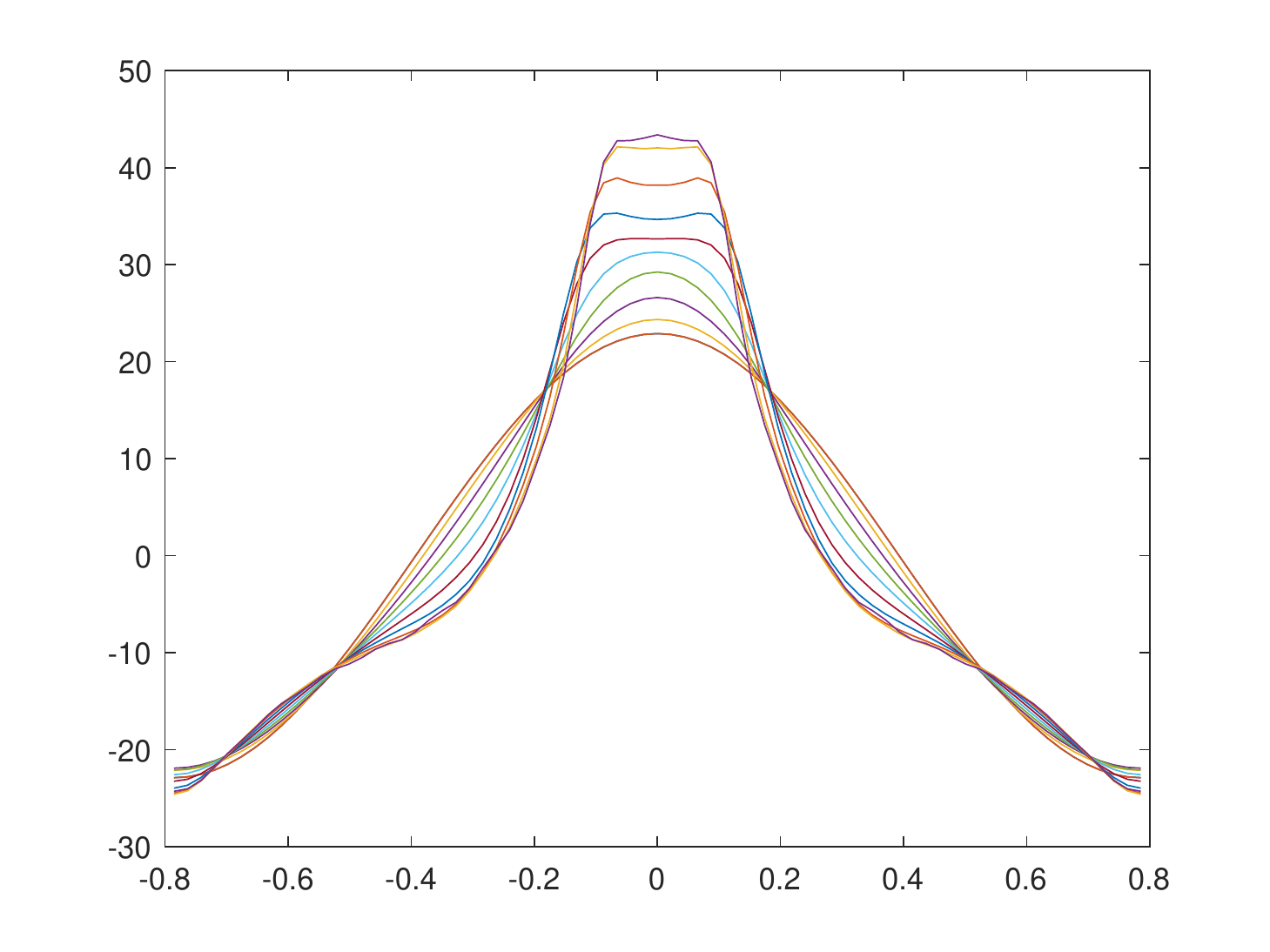}}}
%a centro, b lati
\vspace{-.2cm}
\begin{caption}{\small Behavior of the component $v_1$ at  equispaced time instants in the interval $[0,T]=[0,.4]$:
plots with respect to $r$ for $\theta = 0$ (left); plots with respect to $\theta$ (right) for a value of $r$ 
in the neighborhood of the maximum peak of the graphs on the left.}
\end{caption}
\end{figure}
\end{center}
\vspace{-1.3cm}

\begin{center}
\begin{figure}[h!]
\centerline{{\includegraphics[width=10.8cm,height=8.2cm]{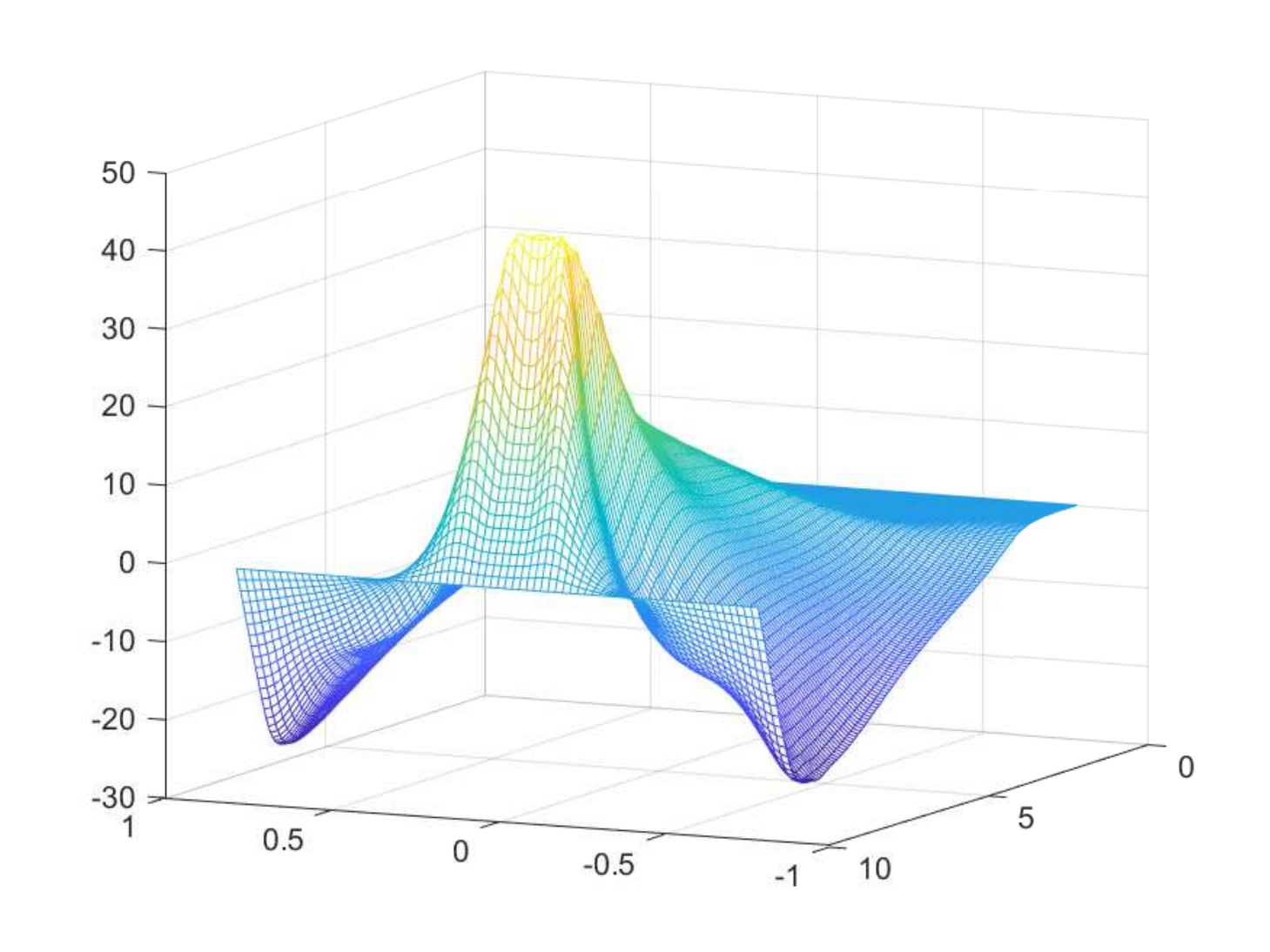}}}
%a centro, b lati
\vspace{-.2cm}
\begin{caption}{\small Plot of the component $v_1$ on the rectangle $[0, 10]\times [-\pi /4, \pi/4]$ at time $T=.4$.}
\end{caption}
\end{figure}
\end{center}
\vspace{-.2cm}

\vfill

It is interesting to observe that, at the points where $\partial \Psi /\partial r=0$, the coefficients of the
nonlinear term in (\ref{syspols}) correspond to the case $\mu_1=1$ and $\mu_2=-4$ for the 1D model problem
(\ref{eq21d}) introduced in section 6. This means that we are in the conditions such that the quantity $Q$ defined in 
(\ref{qua}) may attain different signs depending on the index $j$. In these circumstances, we made the guess that the solution 
of (\ref{eq21d}) blows up in a finite time. Here, we do not have an explosion. However, the behavior looks quite weird,
especially if we examine the picture on the right of Fig.17, in which a plateau is visible in the central part.
Other strange piece-wise like sections are obtained by weighting the terms of the nonlinear term in (\ref{syspols})
in a different manner. Again, we are not in the position to decide whether a break down of regularity is effectively
occurring, or the weirdness is just the consequence of the small diffusive term $\nu$ that allows for the development
of sharp layers without destroying the smoothness. We add further comments in the sections to follow.

\par\medskip
\section{Some theoretical considerations}

In the numerical simulations of section 8, we imposed that the functions $u$ and $\Psi$ had zero average
in $\Omega$, corresponding to the fact that $c_{00}=d_{00}=0$, for any $r$ and any $t$.
This property is compatible with (\ref{finequ}) and (\ref{coeffrenbidc}). Moreover, it is inspired by the fact that the 
nonlinear term in (\ref{fineq2}) is independent of $c_{00}$ and $d_{00}$ (see also (\ref{coeffrenbid})).
Thus, let us study more in detail this aspect. In order to do that, we integrate equation (\ref{fineq2}) in the domain
$\Omega$ and perform some integration by parts by taking into account the Neumann boundary constraints,
valid for any $r$.
Considering that $\int_\Omega u \, d\theta d\phi =0$, for any $r$ and $t$, we get:
$$
\int_\Omega \frac{1}{r}\left[\left(\Delta\Psi +r^2\frac{\partial^2\Psi}{\partial r^2}
+2r\frac{\partial\Psi}{\partial r}\right) \Delta {\hspace{-.11cm}} \left( \frac{\partial \Psi}{\partial r} 
+\frac{\Psi}{r}\right)+ \Delta\Psi \ \Delta {\hspace{-.11cm}} \left( \frac{\partial \Psi}{\partial r} 
-\frac{\Psi}{r}\right)\right] d\theta d\phi
$$
$$=\int_\Omega \frac{2}{r}\Delta\Psi\Delta {\hspace{-.11cm}} \left( \frac{\partial \Psi}{\partial r}\right)d\theta d\phi
-\int_\Omega \left[r \nabla {\hspace{-.11cm}}\left( \frac{\partial^2\Psi}{\partial r^2}\right)\cdot
\nabla {\hspace{-.11cm}}\left( \frac{\partial\Psi}{\partial r}\right)+2\left\vert \nabla {\hspace{-.11cm}}
\left( \frac{\partial\Psi}{\partial r}\right)
\right\vert^2\right]d\theta d\phi
$$
\begin{equation}\label{eqriel}
-\int_\Omega \left[\nabla {\hspace{-.11cm}}\left( \frac{\partial^2\Psi}{\partial r^2}\right)\cdot
\nabla \Psi +\frac{2}{r}\nabla\Psi\cdot \nabla {\hspace{-.11cm}}\left( \frac{\partial\Psi}{\partial r}\right)
\right]d\theta d\phi =0
\end{equation}
\smallskip 

The next step is to integrate the above expression with respect to $r>0$. We denote by $\Sigma$ the cartesian
product $\Omega \times ]0, r_M[$, where $r_M$ can be either finite or infinite. At $r=0$ and $r=r_M$ we impose 
vanishing boundary conditions, independently of $\theta$ and $\phi$. By integrating by parts when necessary, we must have:
\begin{equation}\label{eqriel2}
\int_\Sigma \left[\frac{1}{r^2}(\Delta\Psi)^2 -\frac12 \left\vert \nabla {\hspace{-.11cm}}
\left( \frac{\partial\Psi}{\partial r}\right)\right\vert^2 -\frac{1}{r^2}\vert\nabla\Psi\vert^2
\right]d\theta d\phi dr =0
\end{equation}
The above equality comes from the balance of positive and negative quantities. It does not say too much, except that 
is admissible with the existence of nontrivial functions $\Psi$ solving (\ref{fineq2})-(\ref{finequ}) 
and compatible with the constriction
$c_{00}=d_{00}=0$. If we instead multiply (\ref{eqriel}) by $r$ before the successive integration, the counterpart of
(\ref{eqriel2}) becomes $0=0$. If we finally multiply (\ref{eqriel}) by $r^2$ and integrate, the new version of
(\ref{eqriel2}) is:
\begin{equation}\label{eqriel3}
\int_\Sigma \left[-(\Delta\Psi)^2 +\frac{r^2}{2} \left\vert \nabla {\hspace{-.11cm}}
\left( \frac{\partial\Psi}{\partial r}\right)\right\vert^2 
\right]d\theta d\phi dr =0
\end{equation}
which also has an ambiguous sign. 
\smallskip

The same conclusions can be reached by arguing with the expansions (\ref{coeffrenbid})-(\ref{coeffrenbidc}).
We can substitute the generic coefficient $c_{nl}$, explicited in (\ref{coeffrenbidc}), into (\ref{coeffrenbid}).
Successively, by setting $n=l=0$, the first sum in (\ref{coeffrenbid}) disappears, the second one has $i=j=0$ and $k=m$,
the third one has $k=m=0$ and $i=j$, and the fourth one has $k=m$ and $i=j$. We can analyze the 
terms of the summation, after an integration with respect to the variable $r$. The conclusions are similar to those
of section 6, where, after introducing a suitable quantity $Q$, we distinguished between the case in which
$Q$ maintains the same sign (as a function of the indexes of the summation) or attains different signs.
Here we are in the second situation. 
\smallskip

Things change if we approach the two-dimensional Navier-Stokes problem.
Indeed, if we transfer the same kind of computations to the system
(\ref{upol})-(\ref{syspol2}), we first have:
$$
\int_\Omega \frac{1}{r}\left[\left(\frac{\partial^2\Psi}{\partial \theta^2}
 +r^2\frac{\partial^2\Psi}{\partial r^2}
+r\frac{\partial\Psi}{\partial r}\right) \frac{\partial^3\Psi}{\partial \theta^2\partial r}
+ \frac{\partial^2\Psi}{\partial \theta^2}
 \frac{\partial^2}{\partial \theta^2} {\hspace{-.11cm}} \left( \frac{\partial \Psi}{\partial r} 
-\frac{\Psi}{r}\right)\right] d\theta 
$$
$$=\int_\Omega \frac{2}{r}\frac{\partial^2\Psi}{\partial \theta^2}\, \frac{\partial^3\Psi}{\partial \theta^2\partial r}
d\theta - \int_\Omega \frac{1}{r^2}\left(\frac{\partial^2\Psi}{\partial \theta^2}\right)^{\hspace{-.11cm}2}
d\theta 
$$
\begin{equation}\label{eqrielpo}
-\int_\Omega \left[r\frac{\partial^3\Psi}{\partial \theta\partial r^2}\, \frac{\partial^2\Psi}{\partial \theta\partial r}
+\left(\frac{\partial^2\Psi}{\partial \theta\partial r}\right)^{\hspace{-.11cm}2}\right]
d\theta  =0
\end{equation}
where $\Omega =]-\pi/4, \pi /4[$.
A further integration with respect to $r$, produces:
\begin{equation}\label{eqriel2po}
-\frac12\int_\Sigma \left(\frac{\partial^2\Psi}{\partial \theta\partial r}\right)^{\hspace{-.11cm}2}
d\theta dr =0
\end{equation}
This situation is rather different from that of the three-dimensional case, since the right-hand side
in (\ref{eqriel2po}) is negative and the compatibility 
with $c_0=0$ now only happens for $\Psi =0$. The outcome does not change if we multiply
(\ref{eqrielpo}) by $r$ before integration, so obtaining:
\begin{equation}\label{eqriel3po}
-\int_\Sigma \frac{1}{r}\left(\frac{\partial^2\Psi}{\partial \theta^2}\right)^{\hspace{-.11cm}2}
d\theta dr =0
\end{equation}
\smallskip

The considerations made in section 6 were supported by some numerical tests and suggested as a rule of thumb that, 
when $Q$ has constant sign, the evolutive nonlinear model problem (projected into the subspace of functions with zero 
average) has a unique attractor consisting of the zero function.
On the other hand, when $Q$ attains different signs, there are stable singular solutions that are
reached in a finite time. Can we deduce similar conclusions for the set of Navier-Stokes equations?
Is the behavior of some indicator $Q$ the discriminant factor between the two and the three-dimensional 
cases? We have no answers at the moment, but we hope that the results here discussed
may serve as starting point to advance in this investigation. We also point out that the model problem introduced in 
section 6 might be of interest by itself, both for its mathematical elegance and for possible applications in other contexts.

\par\medskip
\section{Discussion}

There are a few things still to be fixed before concluding this paper. First of all, we
need to say something about the assemblage of the six pyramidal domains representing a partition
of the whole space ${\bf R}^3$ (see Fig.2). The Neumann conditions imposed to $\Psi$
(and consequently to $\Phi$) guarantee that ${\bf v}$ is flattened on each triangular boundary,
for any $r$ and $t$ (see (\ref{defv})). Due to the Neumann conditions imposed on $u$,
from an inspection of (\ref{nlvg}), the above property is also true for the
nonlinear term ${\bf v}\times{\rm curl}{\bf v}$. Thus, the transfer of
information between the domains only takes place through the diffusive term $\nu\bar\Delta {\bf v}$.
After integration over $\Omega$, 
the Laplacian $\Delta u$ can be expressed in weak form and the Neumann boundary conditions
allow for a good match across the interfaces, if we also take into account all the symmetries
involved. As a matter of fact, each normal derivative cancels out the corresponding
normal derivative of the contiguous domain, since the two normal vectors are opposite.
This property is not only true for the 12 triangles
dividing the domains, but also for the 8 straight-lines constituting the boundary of the boundary.
These last are made of the so called {\sl cross-points}. A reasonable initial condition,
such as for instance the one given in (\ref{init}), may ensure a $C^1$ matching across the interfaces.
Of course, global initial data can be chosen as smooth as we please. In the event that some loss of regularity
occurs during the evolution, we expect it to happen at some points in the middle of the pyramids. 
If a deterioration of the regularity appears before at some other places (for instance at the origin or at the interfaces),
it will be anyway a confirmation of the possibility to generate singularities in a finite time.
\smallskip

We did not talk too much about the pressure $p$ in the whole paper. This is also strictly depending
on ${\bf v}$. It is actually defined as the sum of all the potentials than can be plugged in form of a gradient
on the right-hand side of the Navier-Stokes momentum equation. Whatever the expression of $p$ is,
as far as the velocity field remains smooth, we expect the same to 
happen to the pressure. Otherwise, as ${\bf v}$ starts showing a bad behavior, so it will be that of $p$.
\smallskip

A further question concerns with the radial type boundary constraints. At $r=0$ we assume everything
to be zero. Indeed, as seen in our experiments, we expect a reasonably fast decay of the solutions near the origin. 
Nothing interesting will develop there, so that (\ref{init}) seems again a practicable choice.
For the other extreme, i.e. for $r=r_M$, the examples here considered are equivalent to force homogeneous Dirichlet boundary 
conditions on the surface of a sphere.
If we want our problem to be defined in the whole space ${\bf R}^3$ (i.e.: $r_M=+\infty $), we may
require either an appropriate monotone decay at infinity, or an oscillating behavior.
The Bessel's function in (\ref{chidef}) can be an option, since it oscillates remaining bounded for all $r\geq 0$, though it
has not a rapid decay at infinity ($\approx 1/\sqrt{r}$). It is also to be reminded that the
sign of the initial guess influences in different ways the successive development (compare figures 6 and 12).
Presumably, without the Dirichlet type constraint at $r_M$, the vortexes will try to escape outbound, so
we suggest the adoption of an initial function with alternate signs. Unfortunately, our computational capabilities
are not enough to handle these types of experiments.
\smallskip

As a final remark we say that the idea of the six collapsing rings described in section 1 can
be approached as it is, i.e. without resorting to the trick of simplifying the equations through the help
of a fictitious force ${\bf f}$. In alternative, an on purpose attractive radial force (i.e.: $f_1\not =0$), may be added
to speed up the collapsing process. This 3D fluid dynamics exercise can be tackled by a numerical 
code with a certain amount of computational effort. It would be worthwhile to have a try; unforeseen
surprises may come out.

\par

\end{document}